\documentclass[12pt]{elsarticle}
\usepackage[left=1in,right=1in,top=1in,bottom=1in]{geometry}%

\usepackage{amsrefs}
\usepackage{amsmath}
\usepackage{enumerate}
\usepackage{amssymb}                
\usepackage{amsmath}                
\usepackage{amsfonts}
\usepackage{amsthm}
\usepackage{bbm}
\usepackage{float}
\usepackage{mathtools}
\usepackage{graphicx,epsfig}
\usepackage{hyperref}

\RequirePackage{xcolor}[2.11]
\colorlet{siaminlinkcolor}{green!50!black}
\colorlet{siamexlinkcolor}{red!50!black}
\colorlet{siamreviewcolor}{black!50}
\if@hidelinks
\hypersetup{ hidelinks = true }
\else
\hypersetup{
  hidelinks = true,
  colorlinks = false,
  allcolors = siaminlinkcolor,
  urlcolor = siamexlinkcolor,
}
\fi
\makeatletter
\def\ps@pprintTitle{%
   \let\@oddhead\@empty
   \let\@evenhead\@empty
   \def\@oddfoot{\reset@font\hfil\thepage\hfil}
   \let\@evenfoot\@oddfoot
}
\makeatother

\usepackage[capitalize,nameinlink]{cleveref}
\crefname{section}{section}{sections}
\crefname{subsection}{subsection}{subsections}
\Crefname{section}{Section}{Sections}
\Crefname{subsection}{Subsection}{Subsections}

\Crefname{figure}{Figure}{Figures}

\crefformat{equation}{\textup{#2(#1)#3}}
\crefrangeformat{equation}{\textup{#3(#1)#4--#5(#2)#6}}
\crefmultiformat{equation}{\textup{#2(#1)#3}}{ and \textup{#2(#1)#3}}
{, \textup{#2(#1)#3}}{, and \textup{#2(#1)#3}}
\crefrangemultiformat{equation}{\textup{#3(#1)#4--#5(#2)#6}}%
{ and \textup{#3(#1)#4--#5(#2)#6}}{, \textup{#3(#1)#4--#5(#2)#6}}{, and \textup{#3(#1)#4--#5(#2)#6}}

\Crefformat{equation}{#2Equation~\textup{(#1)}#3}
\Crefrangeformat{equation}{Equations~\textup{#3(#1)#4--#5(#2)#6}}
\Crefmultiformat{equation}{Equations~\textup{#2(#1)#3}}{ and \textup{#2(#1)#3}}
{, \textup{#2(#1)#3}}{, and \textup{#2(#1)#3}}
\Crefrangemultiformat{equation}{Equations~\textup{#3(#1)#4--#5(#2)#6}}%
{ and \textup{#3(#1)#4--#5(#2)#6}}{, \textup{#3(#1)#4--#5(#2)#6}}{, and \textup{#3(#1)#4--#5(#2)#6}}

\crefdefaultlabelformat{#2\textup{#1}#3}

\def\R{{\mathbb R}}

\def\C{{\mathbb C}}
\def\Z{{\mathbb Z}}

\DeclareMathOperator{\ran}{ran}

\newtheorem{lemma}{Lemma}
\newtheorem{theorem}{Theorem}
\newtheorem{corollary}{Corollary}

\newtheorem{hypothesis}{Hypothesis}
\newtheorem{remark}{Remark}

\graphicspath{ {images/} }

\begin{document}

\begin{frontmatter}

\title{Existence and spectral stability of multi-pulses in discrete Hamiltonian lattice systems}

\author[1]{Ross Parker}
    \ead{ross\_parker@brown.edu}
\author[2,3]{P.\,G. Kevrekidis} 
    \ead{kevrekid@math.umass.edu}
\author[1]{Bj\"{o}rn Sandstede}
    \ead{bjorn\_sandstede@brown.edu}

\address[1]{Division of Applied Mathematics, Brown University, Providence, RI 02912, USA}
\address[2]{Department of Mathematics and Statistics, University of Massachusetts, Amherst MA 01003, USA}
\address[3]{Mathematical Institute, University of Oxford, Oxford, OX2 6GG, UK}

\begin{abstract}
    In the present work, we consider the existence and spectral stability of multi-pulse solutions in Hamiltonian lattice systems.
    We provide a general framework for the study
    of such wave patterns based on a discrete
    analogue of Lin's method, previously used in
    the continuum realm. We develop explicit conditions for the existence of $m$-pulse structures and
    subsequently develop a reduced matrix allowing us to
    address their spectral stability. As a prototypical
    example for the manifestation of the details of
    the formulation, we consider the discrete nonlinear
    Schr\"{o}dinger equation. Different families
    of $2$- and $3$-pulse solitary waves are discussed,
    and analytical expressions for the corresponding 
    stability eigenvalues are obtained which are in very good agreement
    with numerical results.
\end{abstract}


\end{frontmatter}

\section{Introduction and motivation}

The study of multi-pulse wave structures has a time
honored history in continuum systems. Attempts at
a systematic formulation have taken place both at
a more phenomenological, asymptotic level~\cite{elphick}
and at a more rigorous level~\cite{Sandstede1998}.
The development in the latter work of the so-called
Lin's method for such wave patterns offered a 
systematic view into a reduced formulation where
the characteristics of the pulses (such as their centers,
or possibly also their widths) could constitute effective
dynamical variables for which simpler dynamical equations,
i.e. ordinary differential equations, could be derived. While Lin's method for discrete dynamical systems has been developed in~\cite{Knobloch2000}, it has not so far been applied to the discrete multi-pulse problem.
Over the following decade, methods were sought to 
isolate and freeze the dynamics of individual pulses
within the patterns~\cite{beyn1,beyn2}. More recently,
such freezing techniques have also been extended to other
structures including rotating waves~\cite{beyn3}. 

Despite the intense interest in such multiple coherent
structure patterns at the continuum limit, similar 
techniques have not been systematically developed
at the discrete level. Parts of the relevant efforts
have involved an attempt at adapting the
asymptotic methodology of~\cite{elphick} 
(in the work of~\cite{kevold}) and also the consideration
of structures systematically in the vicinity of the
so-called anti-continuum limit~\cite{Pelinovsky2005}.
The latter setting involves as a starting point the
limit of vanishing coupling between the discrete sites,
whereby suitable Lyapunov-Schmidt conditions can 
be brought to bear to identify persistent configurations
for finite coupling strengths between the adjacent
lattice sites. While works such as~\cite{Kapitula2001a}
have emerged that develop instability criteria, it would
be useful to have a systematic toolbox to study the spectrum of multi-pulses in the spatially discrete setting. This would
serve to both quantify the persistence conditions of
the multi-structure states, and also to offer specific
predictions on their spectral stability and nonlinear
dynamics. 

It is this void that it is the aim of the present work
to fill. We start from a general formulation of Lin's method
for the discrete multi-pulse problem in Hamiltonian systems. (For non-Hamiltonian systems, an adaptation of the results in \cite{Sandstede1998} is also possible). Assuming that a homoclinic orbit exists (the single pulse), we systematically develop conditions for the persistence of multi-pulse states. We then provide estimates of their relevant stability eigenvalues for the low dimensional (reduced) system of the pulses. These eigenvalues are close to 0, and we call them interaction eigenvalues, since they result from nonlinear interaction between neighboring copies of the primary pulse.

As a concrete example for the implementation of the method, we revisit the
discrete nonlinear Schr{\"o}dinger (DNLS) system for which many of the methods of the previous paragraph have been developed~\cite{Kevrekidis2009} (see also~\cite{pelinovsky_2011}). In particular, we give a 
systematic description especially of 2- and 3-pulse solutions and explain how the relevant conclusions can be generalized
to arbitrary multi-pulse structures. 
Our presentation will be structured as follows.
In section 2, we will present the mathematical
setup of the problem and of the special case (DNLS)
example of interest. In section 3, we will develop
Lin's method providing the main results but deferring
the proof details to later sections. In section 4,
we apply the method to the DNLS, comparing the theoretical findings to systematic computations of multi-pulse solutions. Our results are then summarized and some possible directions for future work are offered. Details of the proofs are presented in sections 6-8.

\section{Mathematical setup}

A lattice dynamical system is an infinite system of ordinary differential equations which are indexed by points (nodes)
on a lattice. For the purposes of this work, we will only consider dynamical systems on the integer lattice $\Z$, where the differential equation for each point on the lattice is identical, and the equations are coupled by a centered, second order difference operator.

As a specific example, we will look at the discrete nonlinear Schr{\"o}dinger equation (DNLS)
\begin{equation}\label{DNLS}
i\dot{\psi}_n + d(\psi_{n+1} - 2 \psi_n + \psi_{n-1}) + |\psi_n|^2 \psi_n = 0,
\end{equation}
which is (2.12) in \cite{Kevrekidis2009}, where we have taken $\beta = -1$ and $\sigma = 1$. The parameter $d$ represents the coupling between nodes; $d > 0$ is the focusing case, and $d < 0$ the defocusing case~\cite{Kevrekidis2009}. Equation \cref{DNLS} is Hamiltonian, with energy given by (2.17) in \cite{Kevrekidis2009,pelinovsky_2011}. Of general interest in this type of lattice is the existence and stability of standing waves, which are bound state solutions of the form $\psi_n(t) = e^{i \omega t}\phi_n$~\cite{alfimov}. Making this substitution in \cref{DNLS} and simplifying, a standing wave solves the steady state equation
\begin{equation}\label{DNLSequilib}
d(\phi_{n+1} - 2 \phi_n + \phi_{n-1}) - \omega \phi_n + |\phi_n|^2 u_n = 0.
\end{equation}
From~\cite{herrmann_2011}, a symmetric, real-valued, on-site soliton solution $q_n$ exists to \cref{DNLSequilib} for all $\omega \neq 0$ and $d \geq 0$. This solution $q_n$ furthermore is differentiable in $\omega$. 

We will write DNLS as a system of two real variables $u = (v, w) \in \ell^2(\Z, \R^2)$, where $v = \text{Re }\psi$ and $w = \text{Im }\psi$. In this fashion, we can write \cref{DNLS} in Hamiltonian form as
\begin{equation}\label{DNLSrealHam}
\dot{u}_n + J [\mathcal{H}'(u)]_n = 0,
\end{equation}
where $J$ is the standard skew-symmetric symplectic matrix
\[
J = \begin{pmatrix}0 & 1 \\ -1 & 0\end{pmatrix}
\]
and the Hamiltonian $\mathcal{H}: \ell^2(\Z,\R^2) \rightarrow \R$ is
\begin{equation}\label{DNLSrealH}
\mathcal{H}(v, w) = -\sum_{n = -\infty}^\infty 
\left( \frac{d}{2}\left(v_n - v_{n-1}\right)^2 + \frac{d}{2}\left(w_n - w_{n-1}\right)^2 - \frac{1}{4}\left( v_n^2 + w_n^2 \right)^2 \right).
\end{equation}
The Hamiltonian $\mathcal{H}$ is invariant under the standard rotation group $R(\theta)$, given by
\begin{equation}\label{Rtheta}
R(\theta) = \begin{pmatrix}
\cos(\theta) & \sin(\theta) \\
-\sin(\theta)& \cos(\theta)
\end{pmatrix},
\end{equation}
which has infinitesimal generator $R'(0) = J$. In addition, there is another conserved quantity, often called the norm or the power of the solution, which is given by
\begin{equation}\label{DNLSQ}
\mathcal{Q}(v, w) = \frac{1}{2} \sum_{n = -\infty}^\infty 
\left( v_n^2 + w_n ^2\right).
\end{equation}

Standing waves are solutions of \cref{DNLSrealHam} of the form $R(\omega t) u$, where $u$ is independent of $t$. Substituting this into \cref{DNLSrealHam}, we obtain the equivalent system of equations
\begin{equation}\label{DNLSequilib1}
\mathcal{H}'(u) - \omega u = 0,
\end{equation}
which for DNLS is given by
\begin{equation}\label{DNLSequilib2}
\begin{aligned}
d (v_{n+1} - 2 v_n + v_{n-1}) + v_n w_n^2 + v_n^3 - \omega v_n &= 0 \\
d (w_{n+1} - 2 w_n + w_{n-1}) + v_n^2 w_n + w_n^3 - \omega w_n &= 0 \:.
\end{aligned} 
\end{equation}
If $u$ is a standing wave solution, then $R(\theta) u$ is also a standing wave by symmetry. We note that the steady state system has the form
\begin{equation}
\mathcal{H}'(u) - \omega \mathcal{Q}'(u) = 0,
\end{equation}
which is the stationary equation \cite[(2.15)]{Grillakis1987}. The steady state equation \cref{DNLSequilib2} also has a conserved quantity $E$ \cite{Johansson2000}, which is given by
\begin{equation}\label{DNLSE}
E = 2d(v_n w_{n-1} - v_{n-1} w_n) = 2d \langle u_n, J u_{n-1} \rangle.
\end{equation}
By a conserved quantity in this setting we mean
that this quantity is independent of the lattice
index $n$.

For stability analysis, the linearization of \cref{DNLSrealHam} about a standing wave solution $u^*$ 
yields the linear operator $L(u^*)$, given by 
\begin{equation}\label{DNLSeigproblem}
L(u^*) = J \mathcal{H}''(u^*)  - \omega J.
\end{equation}
Let $u^*_n = (q_n, 0)$, where $q_n$ is the real-valued, on-site standing wave solution to \cref{DNLS}. It is straightforward to verify that
\begin{equation}\label{DNLSkernel1}
\begin{aligned}
L(u^*) R'(0) u^* &= 0 \\
L(u^*) \partial_\omega u^* &= R'(0) u^* \:.
\end{aligned}
\end{equation}
Based on these statements, we have that $0$ is an
eigenvalue with algebraic multiplicity $2$ and
geometric multiplicity $1$ in the DNLS problem. 

\section{Main theorems}

\subsection{Setup}

With DNLS as our principal motivation, we will consider the following more general setting. Consider the Hamiltonian lattice differential equation 
\begin{equation}\label{lattice1}
\dot{u}_n = J [\mathcal{H}'(u)]_n,
\end{equation}
where $u(t) \in \ell^2(\Z, \R^{2k})$, $\mathcal{H}: \ell^2(\Z, \R^{2k}) \rightarrow \R$ is smooth with $\mathcal{H}(0) = 0$ and $\mathcal{H}'(0) = 0$, and $J$ is a $2k \times 2k$ symplectic matrix. For simplicity, and again using DNLS as motivation, we will assume that $\mathcal{H}'(u)$ takes the form
\begin{equation}\label{latticeform}
[\mathcal{H}'(u)]_n = d (\Delta_2 u)_n + f(u_n),
\end{equation}
where $\Delta_2$ is the second difference operator $(\Delta_2 u)_n = u_{n+1} - 2 u_n + u_{n-1}$, $d$ is the coupling constant, and $f: \R^{2k} \rightarrow \R^{2k}$ is smooth with $f(0) = 0$ and $Df(0) = 0$. This implies that, other than the terms from $\Delta_2 u$, the RHS of \cref{lattice1} only involves the lattice site $u_n$. We note that $Df(u(n))$ is self-adjoint since $\mathcal{H}''(u)$ is self-adjoint.

We make the following hypothesis concerning symmetries of the system.
\begin{hypothesis}\label{symmetryhyp}
There is unitary group of symmetries $\{ R(\theta) : \theta \in \R \}$ on $\R^{2k}$ such that 
\begin{enumerate}[(i)]
\item The Hamiltonian $\mathcal{H}$ is invariant under $R(\theta)$, i.e. 
\begin{equation}\label{Hinvariance}
\mathcal{H}(R(\theta)u) = \mathcal{H}(u).
\end{equation}
\item $R'(0) = J$, where $R'(0)$ is the infinitesimal generator of $R(\theta)$.
\end{enumerate}
\end{hypothesis}
\noindent For DNLS, $R(\theta)$ is the rotation group \cref{Rtheta}.

Equilibrium solutions to \cref{lattice1} satisfy 
\begin{equation}\label{equilib1}
\mathcal{H}'(u) = 0.
\end{equation}
Differentiating the symmetry invariance \cref{Hinvariance} as in \cite{Grillakis1987}, we obtain the symmetry relations
\begin{equation}\label{symmetryrel}
\begin{aligned}
\mathcal{H}'(R(\theta)u)) &= R(\theta) \mathcal{H}'(u) \\
\mathcal{H}''(R(\theta)u)) &= R(\theta) \mathcal{H}'(u) R(\theta)^* \:,
\end{aligned}
\end{equation}
from which it follows that $u$ is a solution to \cref{equilib1} if and only if $R(\theta)u$ is a solution. We also note that $f(R(\theta)u) = f(u)R(\theta)$.

We are interested in bound states (referred to also
as standing waves), which are solutions to \cref{lattice1} of the form $u(t) = R(\omega t)u$, where $u \in \ell^2(\Z, \R^{2k})$ is independent of $t$. Bound states satisfy the equilibrium equation
\begin{equation}\label{latticestat}
\mathcal{H}'(u) - \omega u = 0,
\end{equation}
and we note that if $q$ is a bound state, $R(\theta)q$ is also a bound state. Let $q$ be a bound state solution to \cref{latticestat}. The linearization of \cref{lattice1} about a bound state $q$ is the linear operator
\begin{equation}\label{latticeL}
L(q) = J \mathcal{H}''(q) - \omega J.
\end{equation}
By substituting $R(\theta)q$ into \cref{latticestat} and differentiating with respect to $\theta$ at $\theta = 0$, we can verify that 
\begin{equation}\label{Lkernel1}
L(q) R'(0) q = 0.
\end{equation}

\noindent As in \cite{Grillakis1987}, we take the following hypothesis about the existence of bound state solutions.
\begin{hypothesis}\label{boundstatehyp}
For $\omega \in (\omega_1, \omega_2)$, there exists a $C^1$ map $\omega \mapsto q$ such that $q \in \ell^2(\Z, \R^{2k})$ is a bound state solution to \cref{latticestat}.
\end{hypothesis}

By \cref{boundstatehyp}, $\partial_\omega q$ exists for $\omega \in (\omega_1, \omega_2)$. Differentiating \cref{latticeL} with respect to $\omega$, $\partial_\omega q$ satisfies 
\begin{equation}\label{Lkernel2}
L(q)\partial_\omega q = R'(0) q.
\end{equation}
We note that this requires $R'(0) = J$.

\subsection{Spatial dynamics formulation}

We write the bound state equation \cref{latticestat} as the first order difference equation
\begin{equation}\label{diffeq}
U(n+1) = F(U(n)),
\end{equation}
where $U(n) = (u(n), \tilde{u}(n)) = (u_n, u_{n-1}) \in \R^{4k}$ and $F: \R^{4k}\rightarrow \R^{4k}$ is smooth and defined by
\begin{equation}\label{latticeF}
F\begin{pmatrix}u \\ \tilde{u} \end{pmatrix} =
\begin{pmatrix}
\left( 2 + \frac{\omega}{d} \right)u - \frac{1}{d}f(u) - \tilde{u} \\
u
\end{pmatrix}.
\end{equation}
We note that $F(0) = 0$. It is straightforward to verify the symmetry relation 
\begin{equation}
F(T(\theta)U) = F(U)T(\theta),
\end{equation}
where
\begin{equation}
T(\theta) = \begin{pmatrix}
R(\theta) & 0 \\ 0 & R(\theta)
\end{pmatrix}.
\end{equation}

Let $Q(n) = (q_n, q_{n-1})$ be an equilibrium solution to \cref{diffeq}. We can similarly write the eigenvalue problem $(L(q) - \lambda I)v = 0$ as the first order difference equation 
\begin{equation}\label{latticeEVP}
V(n+1) = DF(q_n) V(n) + \lambda B V(n),
\end{equation}
where
\begin{equation}\label{latticeDF}
DF(q(n)) = 
\begin{pmatrix}
2 + \frac{\omega}{d} - \frac{1}{d}Df(q_n) & -I  \\
I & 0
\end{pmatrix} 
\end{equation}
and $B$ is the constant-coefficient block matrix
\begin{equation}\label{latticeB}
B = \frac{1}{d} 
\begin{pmatrix}
J & 0 \\
0 & 0 
\end{pmatrix}.
\end{equation}
It follows from \cref{Lkernel1} that
\begin{equation}\label{DFkernel1}
T'(0) Q(n+1) = DF(Q(n)) T'(0) Q(n).
\end{equation}
We also note that since $J = R'(0)$, $T(\theta)$ commutes with $B$.

Since $F(0) = 0$, 0 is an equilibrium point for the dynamical system \cref{diffeq}. Fix $\omega \in \omega_1, \omega_2)$, and let $q$ be the bound state from \cref{boundstatehyp} corresponding to $\omega$. Let $Q(n) = ( q_n, q_{n-1} )$. Since $q \in \ell^2(\Z, \R^{2k})$, $q_n \rightarrow 0$ as $n \rightarrow \pm \infty$, thus $Q(n)$ is a homoclinic orbit solution to $\cref{diffeq}$ connecting the equilibrium at 0 to itself. We will refer to this as the primary pulse solution. It follows from \cref{Lkernel2} that
\begin{equation}\label{DFkernel2}
\partial_\omega Q(n+1) = DF(Q(n)) \partial_\omega Q(n) + B T'(0)Q(n).
\end{equation}

Since $f(0) = Df(0) = 0$, for the equilibrium at 0 we have
\begin{equation}\label{DF0}
DF(0) = \begin{pmatrix}u \\ \tilde{u} \end{pmatrix} =
\begin{pmatrix}
2 + \frac{\omega}{d} & -I  \\
I & 0
\end{pmatrix} ,
\end{equation}
which has eigenvalues $\mu = \{r, 1/r\}$, each with multiplicity $k$, where
\begin{align}\label{eigr}
r = 1 + \frac{\omega}{2 d} \left( 1 + \sqrt{1 + \frac{4 d}{\omega}} \right).
\end{align}
For $\omega, d > 0$, we have $r > 1$. Thus the equilibrium at 0 is hyperbolic with $k$-dimensional stable and unstable manifolds. The homoclinic orbit $Q(n)$ lies in the intersection of the stable and unstable manifolds. We take the following additional hypothesis regarding their intersection.

\begin{hypothesis}\label{intersectionhyp}
The tangent spaces of the stable and unstable manifolds $W^s(0)$ and $W^u(0)$ have a one-dimensional intersection at $Q(n)$.
\end{hypothesis}

\noindent By \cref{Lkernel1}, this intersection is spanned by $T'(0)Q(n)$. By the stable manifold theorem, we have the decay rate
\begin{equation}\label{Qdecay}
|Q(n)| \leq C r^{-|n|}.
\end{equation}

By \cref{intersectionhyp}, $T'(0) Q(n)$ is the unique bounded solution to the variational equation
\begin{align*}
V(n+1) &= D_U F(q(n)) V(n),
\end{align*}
where
\begin{equation}\label{varsol}
T'(0) Q(n) = \begin{pmatrix} R'(0) q(n) \\ R'(0) q(n-1) \end{pmatrix} = \begin{pmatrix} J q(n) \\ J q(n-1) \end{pmatrix}.
\end{equation}
It follows that there exists a unique bounded solution $Z_1(n)$ to the adjoint variational equation
\begin{align*}
Z(n) &= D_U F(q(n))^* Z(n+1).
\end{align*}
We can verify directly that
\begin{equation}\label{adjvarsol}
Z_1(n) = \begin{pmatrix} R'(0) q(n-1) \\ -R'(0) q(n) \end{pmatrix} = \begin{pmatrix} J q(n-1) \\ -J q(n) \end{pmatrix}.
\end{equation}
In both of these cases, uniqueness is up to scalar multiples.

\subsection{Existence of multi-pulses}

We are interested in multi-pulses, which are bound states that resemble multiple, well separated copies of the primary pulse $Q(n)$. In this section, we give criteria for the existence of multi-pulses. We will characterize a multi-pulse solution in the following way. Let $m > 1$ be the number of copies of $Q(n)$; $N_i$ ($i = 1, \dots, m-1$) be the distances (in lattice points) between consecutive copies; and $\theta_i = G$ ($i = 1, \dots, m$) be symmetry parameters associated with each copy of $Q(n)$. We seek a solution which can be written piecewise in the form 
\begin{equation}\label{Upiecewise}
\begin{aligned}
U_i^-(n) &= T(\theta_i) Q(n) + \tilde{Q}_i^-(n) && n \in [-N_{i-1}^-, 0] \\
U_i^+(n) &= T(\theta_i) Q(n) + \tilde{Q}_i^+(n) && n \in [0, N_i^+] \:,
\end{aligned}
\end{equation}
where $N_i^+ = \lfloor \frac{N_i}{2} \rfloor$, $N_i^- = N_i - N_i^+$, $N_0^- = N_m^+ = \infty$, and
\begin{equation}\label{defN}
N = \frac{1}{2} \min\{ N_i \}.
\end{equation}
The individual pieces are joined together end-to-end as in \cite{Sandstede1998}. The functions $\tilde{Q}_i^\pm(n)$ are remainder terms, which we expect to be small; see the estimates in \cref{transversemulti} below.

In addition to satisfying \cref{diffeq}, the pieces $U_i^\pm(n)$ must match at endpoints of consecutive intervals. Thus, in order to have a multi-pulse solution, $U_i^\pm(n)$ must satisfy the system of equations
\begin{equation}\label{Usystem}
\begin{aligned}
(U_i^\pm)(n+1) &= F(U_i^\pm(n))  \\
U_i^+(N_i^+) - U_{i+1}^-(-N_i^-) &= 0 \\
U_i^+(0) - U_i^-(0) &= 0
\end{aligned}
\end{equation}
for $i = 1, \dots, m$. The first equation in \cref{Usystem} states that the individual pulses are solutions to the difference equation \cref{diffeq} on the appropriate domains; the second equation glues together the individual pulses at their tails; and the third equation is a matching condition at the centers of the pulses.

We will solve \cref{Usystem} using Lin's method. Lin's method yields a solution which has $m$ jumps in the direction of $Z_1(0)$. An $m-$pulse solution exists if and only if all $m$ jumps are 0. These jump conditions are given in the next theorem.

\begin{theorem}\label{ntmulti}
Assume \cref{symmetryhyp}, \cref{boundstatehyp}, and \cref{intersectionhyp}, and let $Q(n)$ be the primary pulse solution to \cref{diffeq}. Then there exists a positive  integer $N_0$ with the following property. For all $m > 1$, pulse distances $N_i \geq N_0$ and symmetry parameters $\theta_i$, there exists a unique $m-$pulse solution $Q_m(n)$ to \cref{diffeq} if and only if the $m$ jump conditions 
\begin{equation}\label{jumpcondexist}
\begin{aligned}
\xi_1 &= \langle T(\theta_1) Z_1(N_1^+), T(\theta_{2}) Q(-N_1^-) \rangle + R_1 = 0 \\
\xi_i &= \langle T(\theta_i) Z_1(N_i^+), T(\theta_{i+1}) Q(-N_i^-) \rangle \\
&\qquad-\langle T(\theta_i) Z_1(-N_{i-1}^-), T(\theta_{i-1}) Q(N_{i-1}^+) \rangle + R_i = 0 && \qquad i = 2, \dots, m-1 \\
\xi_m &= -\langle T(\theta_m) Z_1(-N_{m-1}^-), T(\theta_{m-1}) Q(N_{m-1}^+) \rangle + R_m = 0
\end{aligned}
\end{equation}
are satisfied, where the remainder terms have uniform bound
\[
|R_i| \leq C r^{-3N}.
\]
$Q_m(n)$ can be written piecewise in the form \cref{Upiecewise}, and the following estimates \cref{Westimates} hold:
\begin{equation}\label{Westimates}
\begin{aligned}
\|\tilde{Q}_i^\pm\| &\leq C r^{-N} \\
\tilde{Q}_i^+(N_i^+) &= T(\theta_{i+1}) Q(-N_i^-) + \mathcal{O}(r^{-2N}) \\
\tilde{Q}_{i+1}^-(-N_i^-) &= T(\theta_i) Q(N_i^+) + \mathcal{O}(r^{-2N}) .
\end{aligned}
\end{equation}
\end{theorem}

\begin{remark}If equation \cref{diffeq} has a conserved quantity, i.e. a function $E: \R^{4d} \rightarrow \R$ such that $E(F(U)) = E(U)$, we can remove one of the jump conditions in \cref{jumpcondexist} as is done in \cite{SandstedeStrut}.
\end{remark}

\subsection{Eigenvalue problem}

We will now turn to the spectral stability of multi-pulses. In particular, we will locate the interaction eigenvalues. Let $Q_m(n) = (q_m(n), q_m(n-1))$ be an $m-$pulse solution to \cref{diffeq} constructed according to \cref{ntmulti}. By \cref{ntmulti}, $Q_m(n)$ can be written piecewise in the form \cref{Upiecewise}. The eigenvalue problem is
\begin{equation}\label{multiEVP}
V(n+1) = DF(q_m(n)) V(n) + \lambda B V(n),
\end{equation}
where $DF(q_m(n))$ and $B$ are given by \cref{latticeDF} and \cref{latticeB}. Since $q_m(n)$ decays exponentially to 0 and $F$ is smooth, $DF(q_m(n))$ is exponentially asymptotic to the constant coefficient matrix $DF(0)$, which is hyperbolic.

We will also assume a Melnikov sum condition holds. Since we have a Hamiltonian system, the standard Melnikov sum $M_1$ is 0 since $J$ is skew-symmetric.
\begin{align}\label{MelnikovM1zero}
M_1 &= \sum_{n=-\infty}^\infty \langle Z_1(n+1), B T'(0)Q(n) \rangle
= \sum_{n=-\infty}^\infty \langle R'(0)q(n), J R'(0)q(n) \rangle
= 0.
\end{align}
We note if $M_1 \neq 0$, which can occur in non-Hamiltonian systems, the analysis is much simpler and in fact is the discrete analogue of \cite{Sandstede1998}. We will assume that the following higher order Melnikov sum is nonzero.

\begin{hypothesis}\label{melnikovhyp}
The following Melnikov-like condition holds.
\begin{align}\label{M2cond}
M_2 &= \sum_{n=-\infty}^\infty \langle Z_1(n+1), B \partial_\omega Q(n) \rangle = \sum_{n=-\infty}^\infty \langle R'(0) q(n), \partial_\omega q(n) \rangle \neq 0 .
\end{align}
\end{hypothesis}
\noindent In general, the Melnikov condition \cref{M2cond} can only be verified numerically.

We can now state the following theorem, in which we locate the eigenvalues of \cref{latticeEVP} resulting from interactions between neighboring pulses.

\begin{theorem}\label{stabilitytheorem}
Assume \cref{symmetryhyp}, \cref{boundstatehyp}, \cref{intersectionhyp}, and \cref{melnikovhyp}. Let $Q_m(n)$ be an $m-$pulse solution to \cref{diffeq} constructed according to \cref{ntmulti} with pulse distances $\{ N_1, \dots, N_{m-1}\}$ and symmetry parameters $\{\theta_1, \dots, \theta_m\}$. Then there exists $\delta > 0$ small with the following property. There exists a bounded, nonzero solution $V(n)$ of the eigenvalue problem \cref{multiEVP} for $|\lambda| < \delta$ if and only if $E(\lambda) = 0$, where
\begin{equation}\label{Elambda}
E(\lambda) = \det(A - M_2 \lambda^2 I + R(\lambda)).
\end{equation}
$M_2$ is defined in \cref{melnikovhyp}, and $A$ is the tridiagonal $m \times m$ matrix
\begin{align}\label{matrixA}
A &= \begin{pmatrix}
-a_1 & a_1 & & & \\
-\tilde{a}_1 & \tilde{a}_1 - a_2 & a_2 \\
& -\tilde{a}_2 & \tilde{a}_2 - a_3 & a_3 \\
& \ddots & & \ddots \\
& & & -\tilde{a}_{m-1} & \tilde{a}_{m-1}  \\
\end{pmatrix},
\end{align}
where
\begin{align*}
a_i &= \langle T(\theta_i) Z_1(N_i^+), T(\theta_{i+1}) T'(0)Q(-N_i^-) \rangle \\
\tilde{a}_i &= \langle T(\theta_{i+1}) Z_1(-N_i^-), T(\theta_i) T'(0)Q(N_i^+) \rangle \:.
\end{align*}
The remainder term has uniform bound
\begin{align}\label{Rbound2}
|R(\lambda)| \leq C\left( (r^{-N} + |\lambda|)^3 \right),
\end{align}
where $N = \frac{1}{2}\min\{N_1, \dots, N_{m-1} \}$.
\end{theorem}

\subsection{Transverse intersection}

We present one more result, which concerns the existence of multi-pulse solutions in the case where the stable manifold $W^s(0)$ and unstable manifold $W^u(0)$ intersect transversely, as opposed to the one-dimensional intersection in \cref{intersectionhyp}. This is particularly useful for DNLS, as this occurs when we consider its real-valued solutions. In the transverse intersection case, we have a much more general result. Consider the difference equation
\begin{equation}\label{diffeqtransv}
U(n+1) = F(u(n)),
\end{equation}
where $F: \R^k \rightarrow \R^k$ is smooth. We make the following assumptions about $F$.

\begin{hypothesis}\label{transversehyp}
The following hold concerning the function $F$.
\begin{enumerate}[(i)]
\item There exists a finite group $G$ (which may be the trivial group) for which the group action is a unitary group of symmetries $T(\theta)$ on $\R^k$ such that 
\begin{equation}\label{symmetrytransverse}
F(T(\theta)U) = F(U)T(\theta)
\end{equation}
for all $\theta \in G$ and all $U \in \R^k$. 
\item 0 is a hyperbolic equilibrium for $F$, thus there exists a radius $r > 1$ such that for all eigenvalues $\nu$ of $DF(0)$, $|\nu| \leq 1/r$ or $|\nu| > r$. Furthermore, $\dim E^s, \dim E^u \geq 1$, where $E^s$ and $E^u$ are the stable and unstable eigenspaces of $DF(0)$.
\item There exists a primary pulse homoclinic orbit solution $Q(n)$ to \cref{diffeqtransv} which connects the equilibrium at 0 to itself.
\item The stable and unstable manifolds $W^s(0)$ and $W^u(0)$ intersect transversely.
\end{enumerate}
\end{hypothesis}
We note that for DNLS, the group $G$ is $( \{\pm 1\}, \cdot)$. In this case, Lin's method yields a unique $m-$pulse solution to \cref{diffeqtransv}.

\begin{theorem}\label{transversemulti}
Assume \cref{transversehyp}, and let $Q(n)$ be the primary pulse solution to \cref{diffeq}. Then there exists a positive integer $N_0$ with the following property. For all $m > 1$, pulse distances $N_i \geq N_0$ and symmetry parameters $\theta_i$, there exists a unique $m-$pulse solution $Q_m(n)$ to \cref{diffeq} which can be written in the form \cref{Upiecewise}. The remainder terms $\tilde{Q}_i^\pm(n)$ have the same estimates as in \cref{ntmulti}.
\end{theorem}

\section{Discrete NLS equation}\label{sec:DNLS}

\subsection{Background}

We will now apply the results the previous section to the 
DNLS to illustrate the impact of the discrete Lin's method. Before we do that, we will give a brief overview what is already known. Many more details can be found in \cite{Kevrekidis2009,pelinovsky_2011}. 

At the anti-continuum limit, equation \cref{DNLSequilib} reduces to a system of decoupled algebraic equations. Any $u_n$ with $u_n \in \{ 0, \pm \sqrt{\omega}\}$ is a solution. For $d > 0$, the DNLS possesses two real-valued, symmetric, single pulse solutions (up to rotation): on-site solutions, which are centered on a single lattice point; and off-site solutions, which are centered between two adjacent lattice points \cite{Kevrekidis2009}. The on-site solution has a single eigenvalue at 0 from rotational symmetry. The off-site solution has an additional pair of real eigenvalues; since the off-site solution is spectrally unstable, we will only consider the on-site solution from here on as the foundation for the 
single pulse state. 

For sufficiently small $d$, $m$-pulse solutions exist to equation \cref{DNLSequilib} for any pulse distances as long as the phase differences satisfy $\Delta \theta_i \in \{0, \pi\}$ \cite[Proposition 2.1]{Pelinovsky2005}. For sufficiently small $d$, this $m-$pulse is spectrally unstable unless all of the phase differences $\Delta \theta_i$ are $\pi$; in that case there are $m-1$ pairs of purely imaginary eigenvalues with negative Krein signature \cite[Theorem 3.6]{Pelinovsky2005}.
This means that these eigendirections, although neutrally
stable, are prone to instabilities when parameters 
(such as $d$) are varied upon collision with other
eigenvalues. However, they may also lead to instabilities
at a purely nonlinear level (despite potential spectral
stability) as a result of the mechanism explored, e.g.,
in~\cite{CUCCAGNA200938,PRL_2015}.
For any $d$ for which the $m-$pulse exists, if one or more phase differences $\Delta \theta_i$ is 0, it follows from Sturm-Liouville theory that there is at least one positive, real eigenvalue \cite{Kapitula2001a}.

\subsection{Main results}

Let $q(n)$ be the on-site, real-valued soliton solution to \cref{DNLSequilib}. We will characterize an $m-$pulse solution to \cref{DNLSequilib} in terms of the $m-1$ pulse distances $\{ N_1, \dots, N_{m-1} \}$ and phase differences $\{ \Delta\theta_1, \dots, \Delta\theta_{m-1} \}$ between consecutive copies of $q(n)$. We have the following theorem regarding the existence of $m-$pulse solutions.

\begin{theorem}\label{DNLSexisttheorem}
There exists a positive integer $N_0$ (which depends on $\omega$ and $d$), with the following property. For any $m \geq 2$, pulse distances $N_i \geq N_0$, and phase differences $\Delta\theta_i \in \{0, \pi\}$, there exists a unique $m-$pulse solution $q_m(n)$ to \cref{DNLSequilib} which resembles $m$ consecutive copies of the on-site pulse $q(n)$. No other phase differences are possible.
\end{theorem}

By \cref{Lkernel1} and \cref{Lkernel2}, the linearization about $q_m(n)$ has a kernel with algebraic multiplicity 2 and geometric multiplicity 1 which is a result of rotational invariance. The following theorem locates the small eigenvalues of the linearization about $q_m(n)$ resulting from interaction between consecutive copies of $q(n)$. 

\begin{theorem}\label{DNLSeigtheorem}
Let $q_m(n)$ be an $m-$pulse solution to \cref{DNLSequilib} with pulse distances $N_i$ and phase differences $\Delta\theta_i$. Assume that $M > 0$, where
\[
M = \sum_{n=-\infty}^\infty q_n \partial_\omega q_n = \partial_\omega \left( \frac{1}{2} \sum_{n=-\infty}^\infty q_n^2 \right).
\]
Let $N = \frac{1}{2} \min\{ N_1, \dots, N_{m-1}\}$. Then for $N$ sufficiently large, there exist $m-1$ pairs of interaction eigenvalues $\{\pm \lambda_1, \dots, \pm \lambda_{m-1}\}$, which can be grouped as follows. There are $k_\pi$ pairs of purely imaginary eigenvalues and $k_0$ pairs of real eigenvalues, where $k_\pi$ is the number of phase differences $\Delta\theta_i$ which are $\pi$, and $k_0$ is the number of phase differences $\Delta\theta_i$ which are $0$. The $\lambda_j$ are close to 0 and are given by the following formula
\begin{align}\label{eigsDNLS}
\lambda_j &= \sqrt{\frac{d \mu_j}{M}} + \mathcal{O}(r^{-2N}) && j = 1, \dots, m-1,
\end{align}
where $r$ is defined in \cref{eigr}, $d$ is the coupling
constant and $\{ \mu_1, \dots, \mu_{m-1} \}$ are the distinct, real, nonzero eigenvalues of the symmetric, tridiagonal matrix
\begin{align}\label{DNLSmatrixA}
A &= \begin{pmatrix}
-\cos(\Delta\theta_1) b_1 & \cos(\Delta\theta_1) b_1 & & &  \\
\cos(\Delta\theta_1) b_1 & -\cos(\Delta\theta_1) b_1 - \cos(\Delta\theta_2) b_2 & \cos(\Delta\theta_2) b_2 \\
& \ddots & \ddots \\
& &  \cos(\Delta\theta_{m-1}) b_{m-1} & -\cos(\Delta\theta_{m-1}) b_{m-1}  \\
\end{pmatrix},
\end{align}
where
\begin{align}\label{bieq}
b_i = \begin{cases}
q\left(\frac{N_i}{2}\right) \left[ q\left(\frac{N_i}{2} + 1\right) - q\left(\frac{N_i}{2} - 1\right) \right] & N_i \text{ even} \\
q\left(\frac{N_i-1}{2}\right)\left(\frac{N_i+3}{2}\right) 
- \left(\frac{N_i+1}{2}\right)q\left(\frac{N_i-3}{2}\right) & N_i \text{ odd}
\end{cases} \:.
\end{align}
\end{theorem}

\begin{remark}
There is strong numerical evidence that $M > 0$, i.e. the Melnikov condition is satisfied.
\end{remark}

\begin{remark}
If all the nonzero eigenvalues $\mu_j$ of $A$ are larger than $\mathcal{O}(r^{-4N})$, then the formula \cref{eigsDNLS} is the sum of a leading order term and a small remainder term. A good approximation for the eigenvalues $\lambda_j$ can be obtained by computing the eigenvalues of $A$. A sufficient condition for this is $N_{\mathrm{max}} < 2 N$, where $N_{\mathrm{max}} = \frac{1}{2} \max\{ N_1, \dots, N_{m-1}\}$.
\end{remark}

\noindent In addition, we remark that if $b_1 = \dots = b_{m-1} = b$, $A = -b \mathcal{M}_1$, where the matrix $\mathcal{M}_1$ is defined in \cite[(2.84)]{Kevrekidis2009} and represents interactions between neighboring sites.

We can compute the nonzero eigenvalues of \cref{DNLSmatrixA} in several special cases. In the first corollary, we consider the case where the pulse distances $N_i$ are equal.

\begin{corollary}\label{DNLSeigcorr}Let $q_m(n)$ be an $m-$pulse solution to \cref{DNLSequilib} with pulse distances $N_i = 2N$ and phase differences $\Delta\theta_i$. Then the $\lambda_j$ are as follows.
\begin{enumerate}[(i)]
\item For $m = 2$, we have
\begin{align}\label{2pulseeigs}
\lambda_1 &= 
\begin{cases}
\sqrt{2}\nu  + \mathcal{O}(r^{-2N}) & \Delta\theta_1 = 0 \\
\sqrt{2}\nu i + \mathcal{O}(r^{-2N}) & \Delta\theta_1 = \pi
\end{cases} \:.
\end{align}
\item For $m = 3$, we have
\begin{align}\label{3pulseequaleigs}
\lambda_{1, 2} &= \begin{cases}
\nu, \sqrt{3} \nu + \mathcal{O}(r^{-2N}) & (\Delta\theta_1, \Delta\theta_2) = (0, 0) \\
3^{1/4}\nu, 3^{1/4}\nu i + \mathcal{O}(r^{-2N}) & (\Delta\theta_1, \Delta\theta_2) = (0, \pi) \\
\nu i, \sqrt{3} \nu i + \mathcal{O}(r^{-2N}) & (\Delta\theta_1, \Delta\theta_2) = (\pi, \pi)
\end{cases}\:.
\end{align}
\item For $m > 3$, if $\Delta\theta_i = \Delta\theta$ for all $i$,
\begin{align*}
\lambda_j = \begin{cases}
\sqrt{2\left( \cos\frac{\pi j}{m} - 1 \right)}\nu + \mathcal{O}(r^{-2N}) & \Delta\theta = 0 \\
\sqrt{2\left( \cos\frac{\pi j}{m} - 1 \right)}\nu i + \mathcal{O}(r^{-2N}) & \Delta\theta = \pi
\end{cases}
\end{align*}
for $j = 1, \dots, m-1$.
\end{enumerate}
where $\nu = \sqrt{\frac{|b|d}{M}} = \mathcal{O}(r^{-N})$, and $b$ is given by equation \cref{bieq}.
\end{corollary}

In the second corollary, we give a general formula for the eigenvalues for a 3-pulse.

\begin{corollary}\label{DNLSeigcorr2}
Let $q_3(n)$ be an $3-$pulse solution to \cref{DNLSequilib} with pulse distances $N_1, N_2$ and phase differences $\Delta \theta_1, \Delta \theta_2$. Then $\lambda_1, \lambda_2$ are given by
\begin{equation}\label{3pulseeigs}
\begin{aligned}
\lambda_{1,2} = \sqrt{\frac{d}{M}}
&\Big( -b_1\cos\Delta\theta_1 - b_2\cos\Delta\theta_2  \\
&\pm \sqrt{b_1^2 + b_2^2 - b_1 b_2\cos\Delta\theta_1 \cos\Delta\theta_2} \Big)^{1/2} + \mathcal{O}(r^{-2N}).
\end{aligned}
\end{equation}
If $N_1 < N_2 < 2 N_1$, then, to leading order, these have magnitude
\begin{equation}\label{3pulsemag}
\begin{aligned}
|\lambda_1| &= \sqrt{\frac{2 |b_1| d}{M}} = \mathcal{O}(r^{-N_1/2}) \\
|\lambda_2| &= \sqrt{\frac{3 |b_2| d}{2 M}} = \mathcal{O}(r^{-N_2/2}) \:,
\end{aligned}
\end{equation}
where $b_1$ and $b_2$ are given by equation \cref{bieq}.
\end{corollary}

\subsection{Numerical results}

In this section, we provide numerical verification for the results in the previous section. We first construct multi-pulse solutions to the steady state DNLS problem by using Matlab for parameter continuation in the coupling constant $d$ from the anti-continuum limit. We then find the eigenvalues of the linearization about this solution using Matlab's \texttt{eig} function. 

\begin{figure}
\centering
\includegraphics[width=5cm]{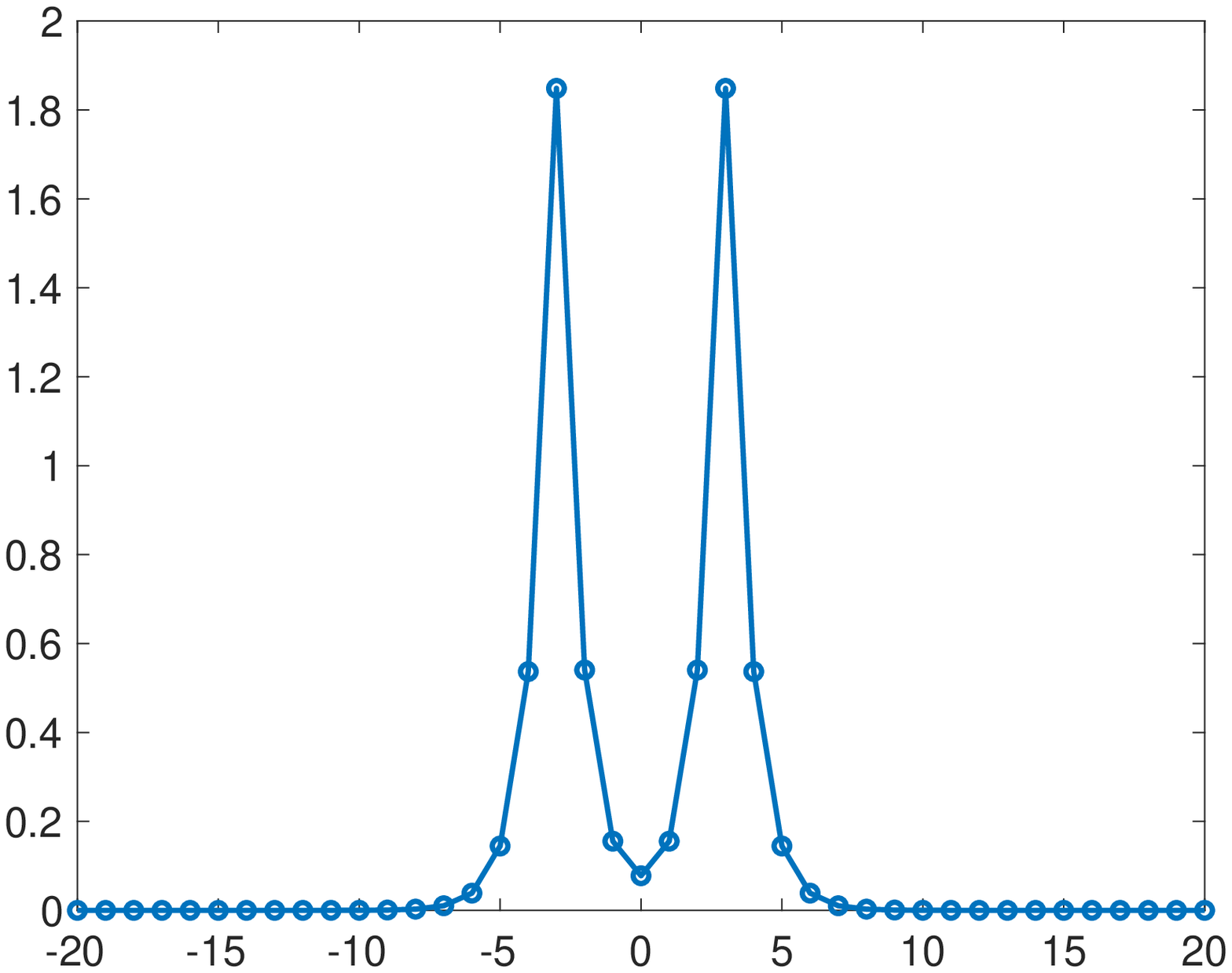}
\includegraphics[width=5cm]{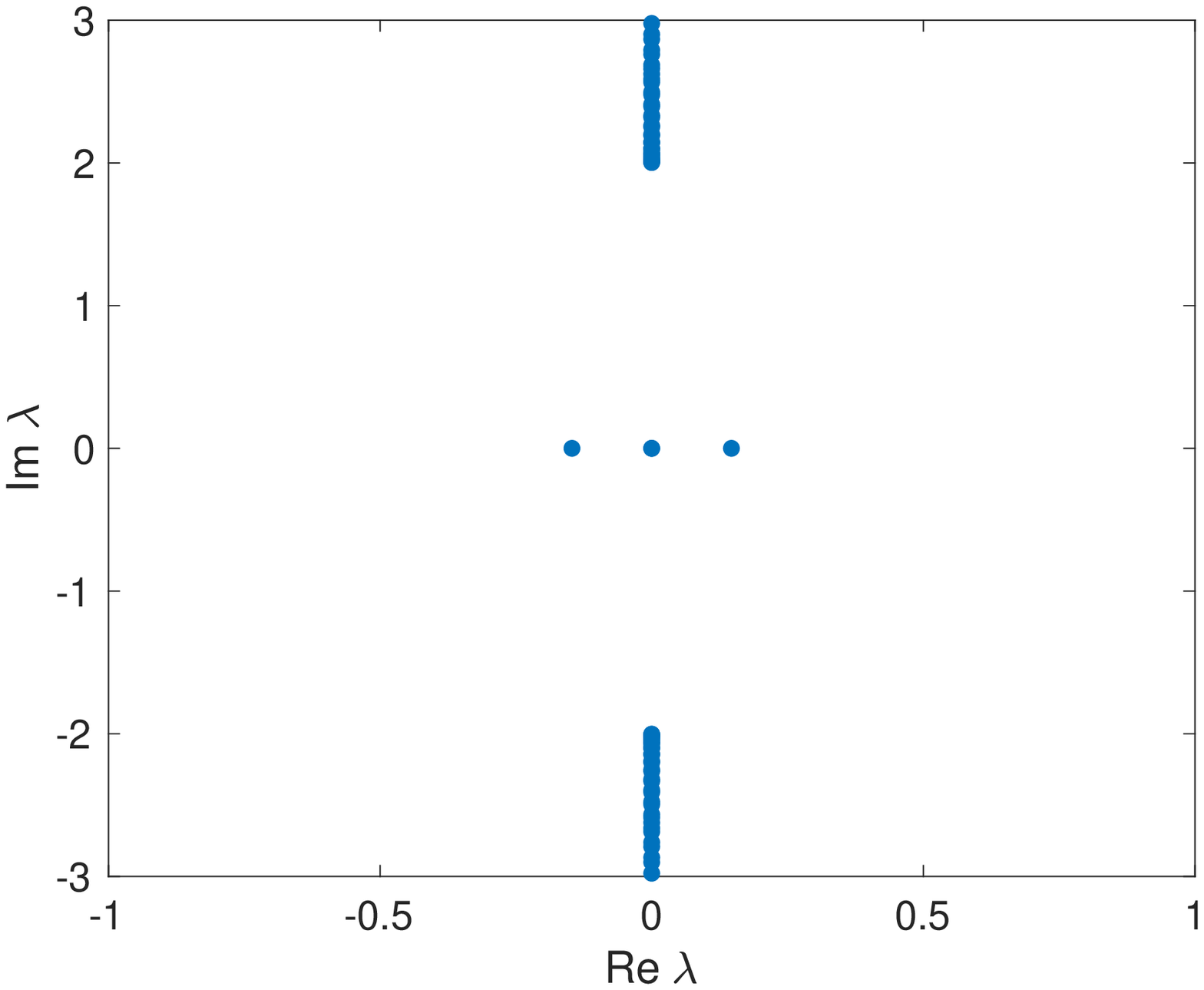}
\includegraphics[width=5cm]{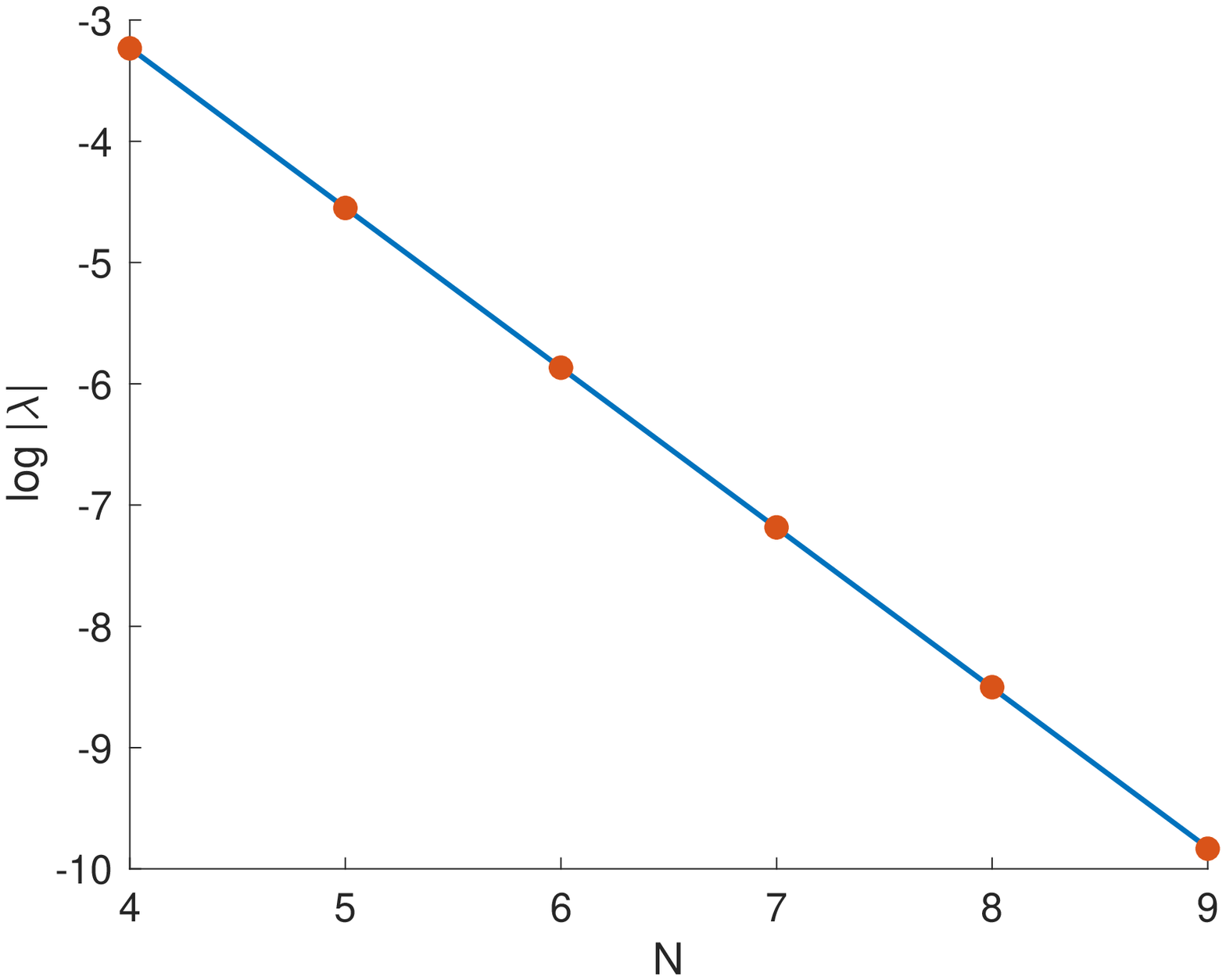}
\includegraphics[width=5cm]{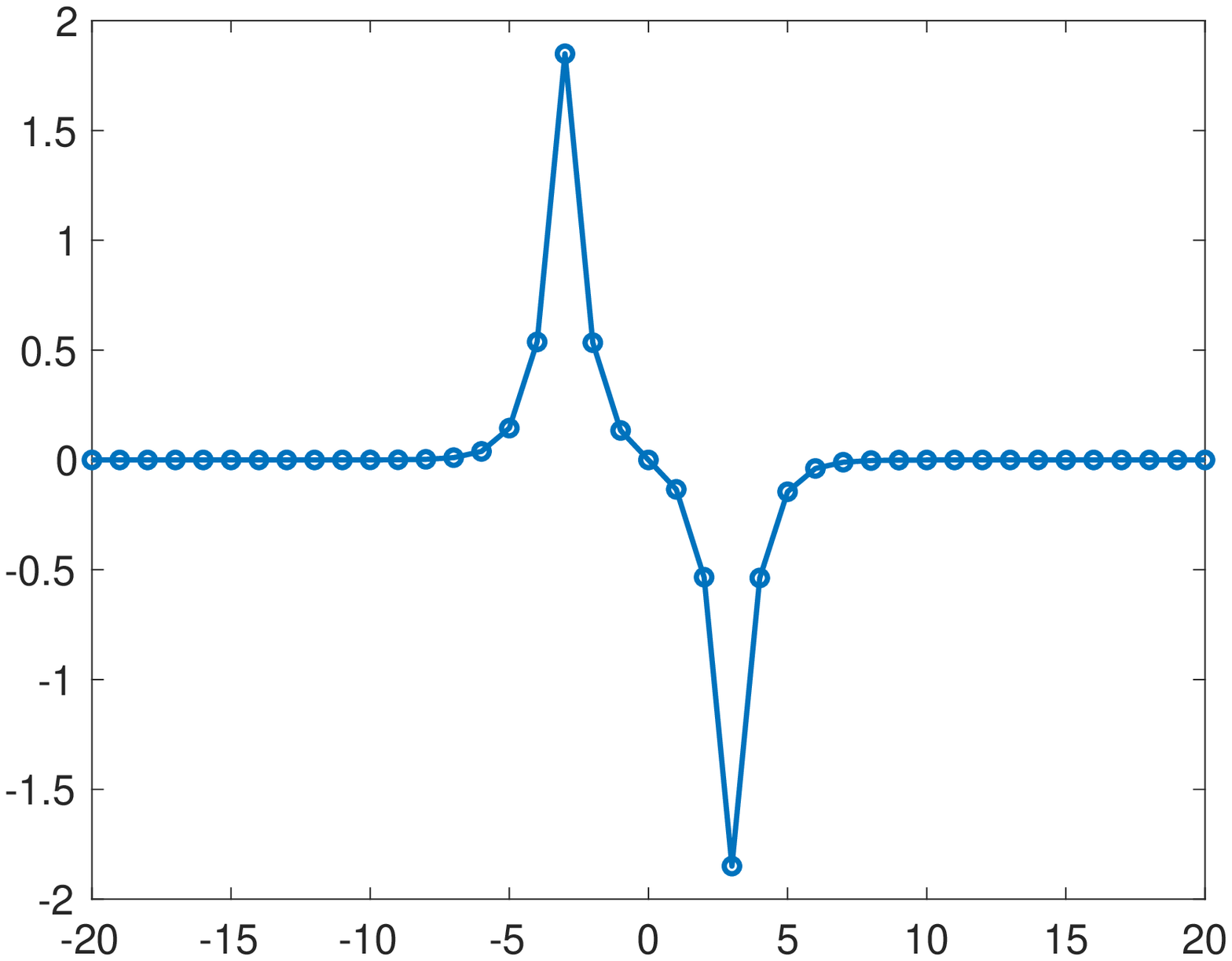}
\includegraphics[width=5cm]{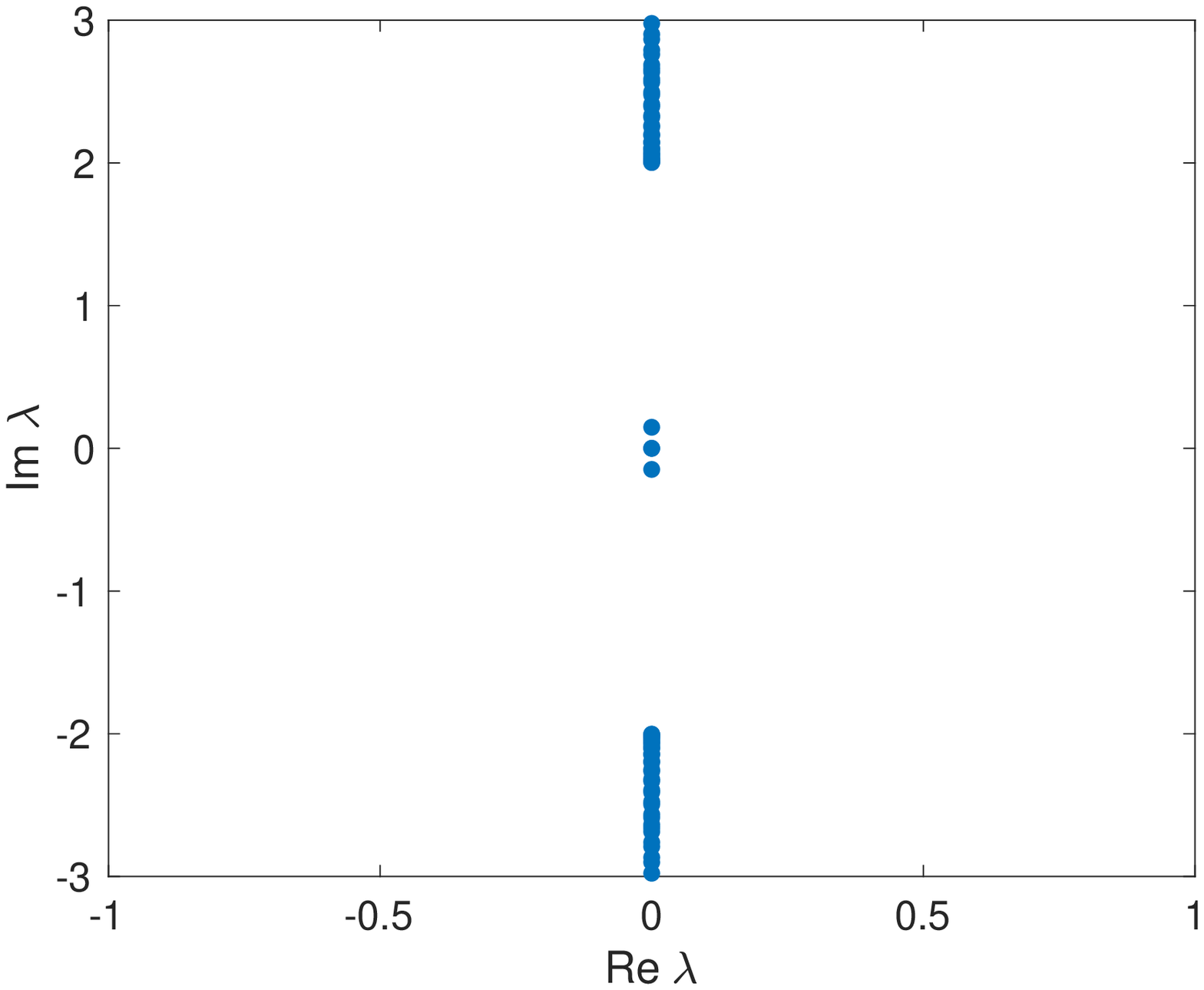}
\includegraphics[width=5cm]{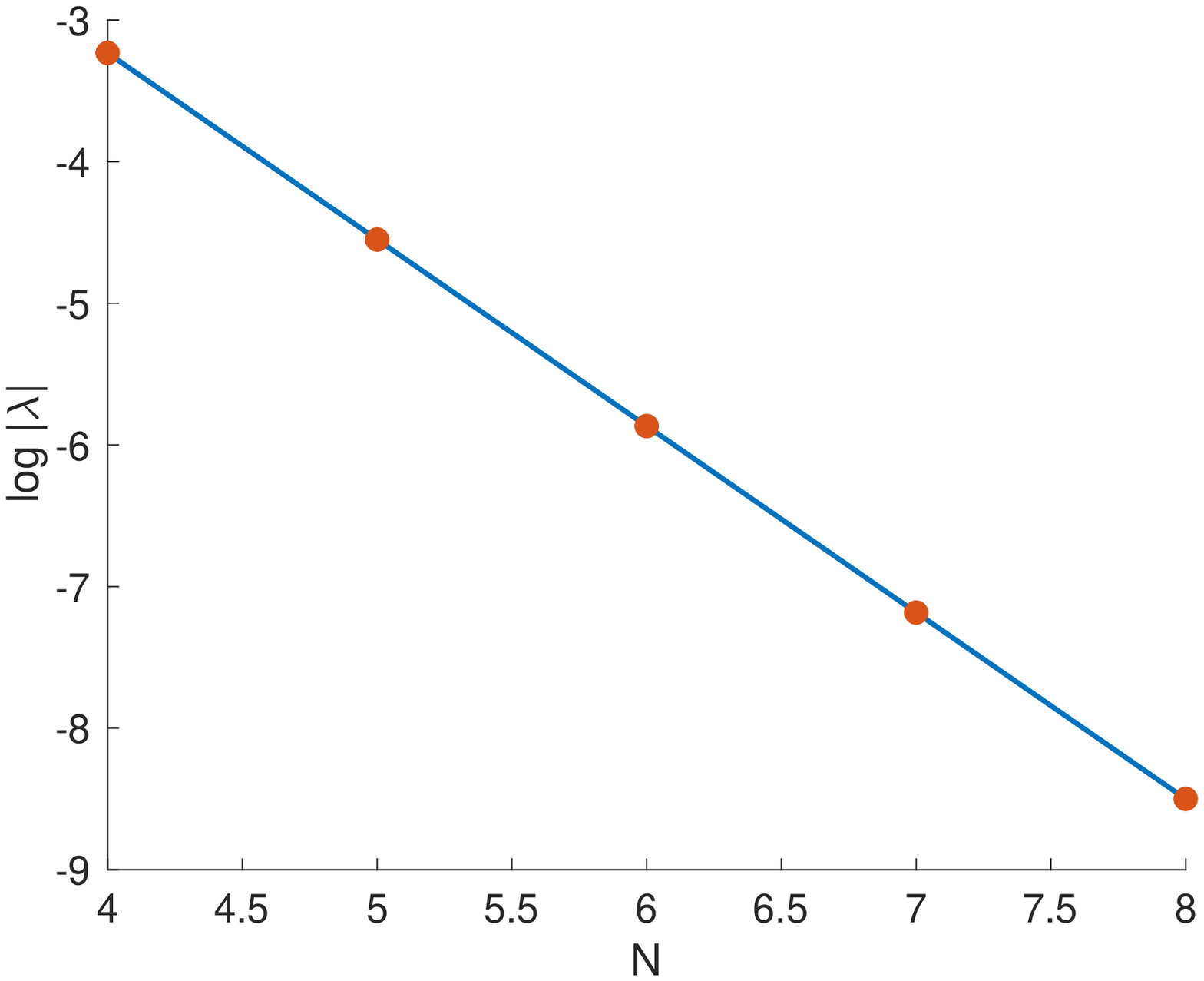}
\caption{Solution profile (left panel), spectral plane eigenvalue pattern (center panel), and plot of $\log(\lambda)$ vs. $N$ with least squares linear regression line (right panel) for $++$ (top) and $+-$ (bottom) pulses. The symbolic notation here and 
below follows that of~\cite{alfimov}, referring
with a symbolic sign representation to the positive
or negative value of the peak of the pulse.
Parameters $\omega = 2$ and $d = 1.0$.}
\label{fig:eigendecay1}
\end{figure}

First, we look at multi-pulses where the pulse distances are equal.  The left and center panels of \cref{fig:eigendecay1} show the pulse profile and eigenvalue pattern for the two double pulses (of relative
phase $0$ and $\pi$). Equation \cref{2pulseeigs} from \cref{DNLSeigcorr} states that for fixed $\omega$ and $d$, the interaction eigenvalues decay as $r^{-N}$. In the right panel of \cref{fig:eigendecay1}, we plot $\log \lambda$ vs. $N$ for the two possible double pulses and construct a least-squares linear regression line. In both cases, the relative error in the slope of this line (which is predicted to be $-\log r$) is order $10^{-4}$. This result provides theoretical and numerical support to the earlier observations of~\cite{Kapitula2001a}.

\begin{figure}
\centering
\includegraphics[width=5cm]{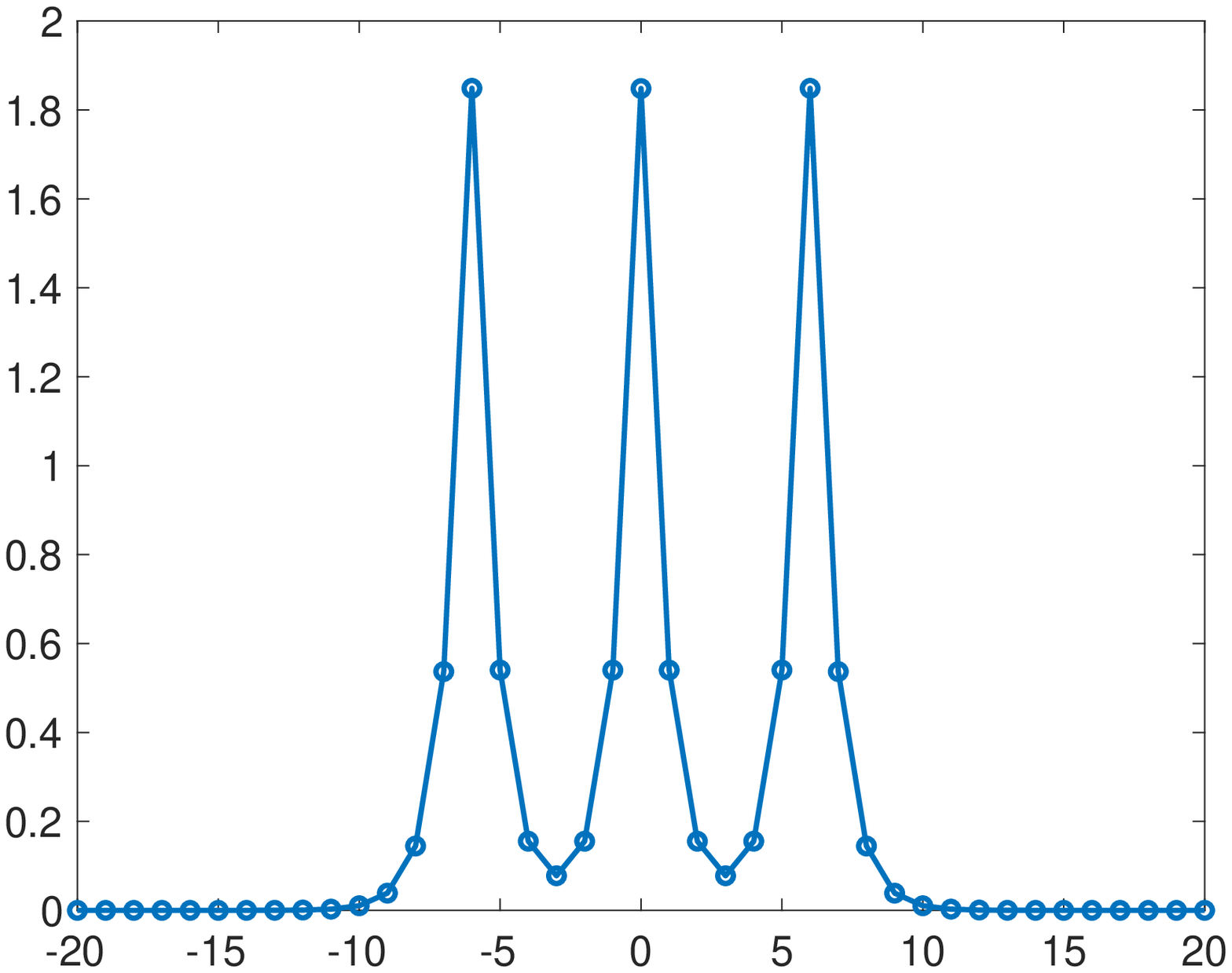}
\includegraphics[width=5cm]{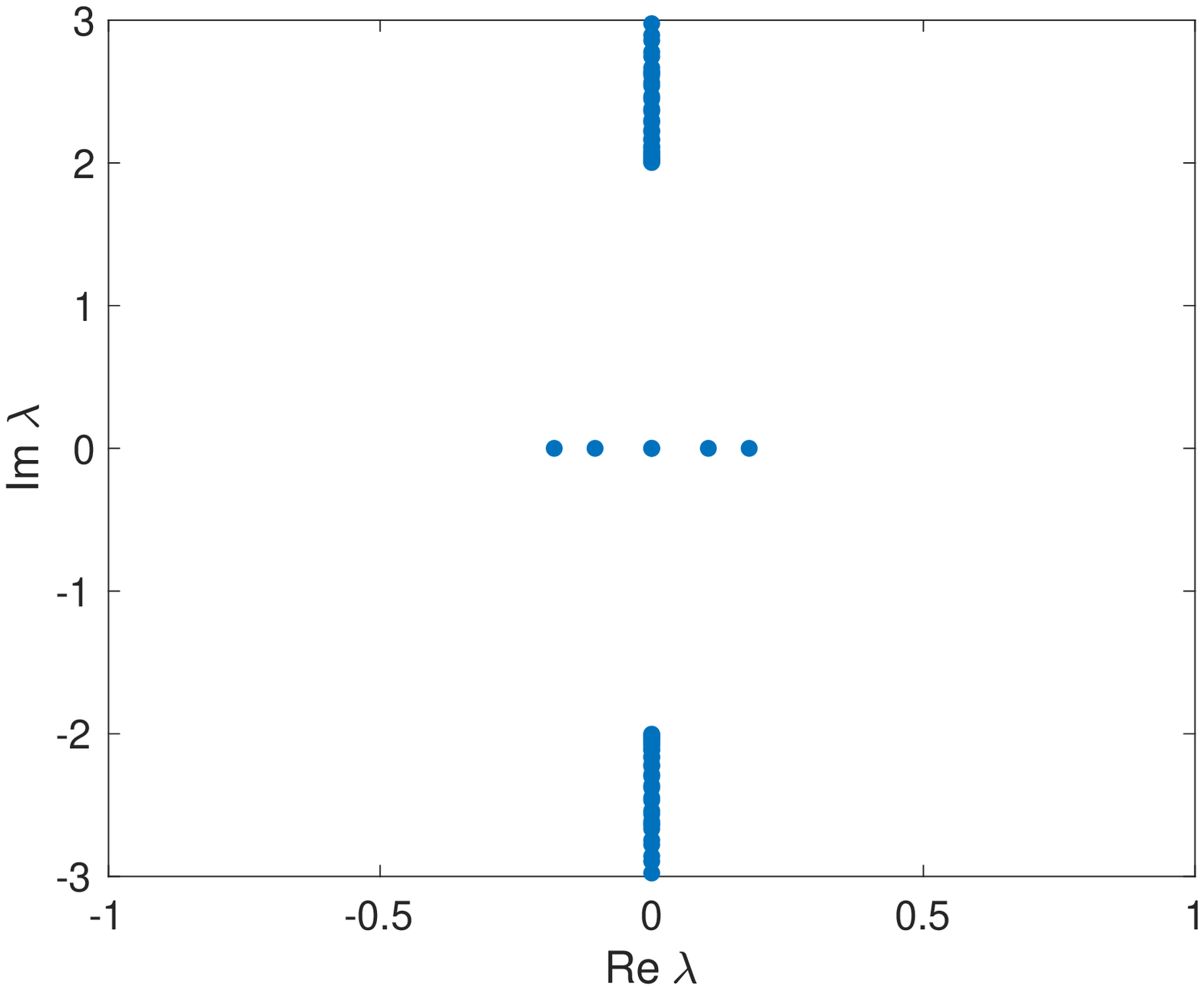}
\includegraphics[width=5cm]{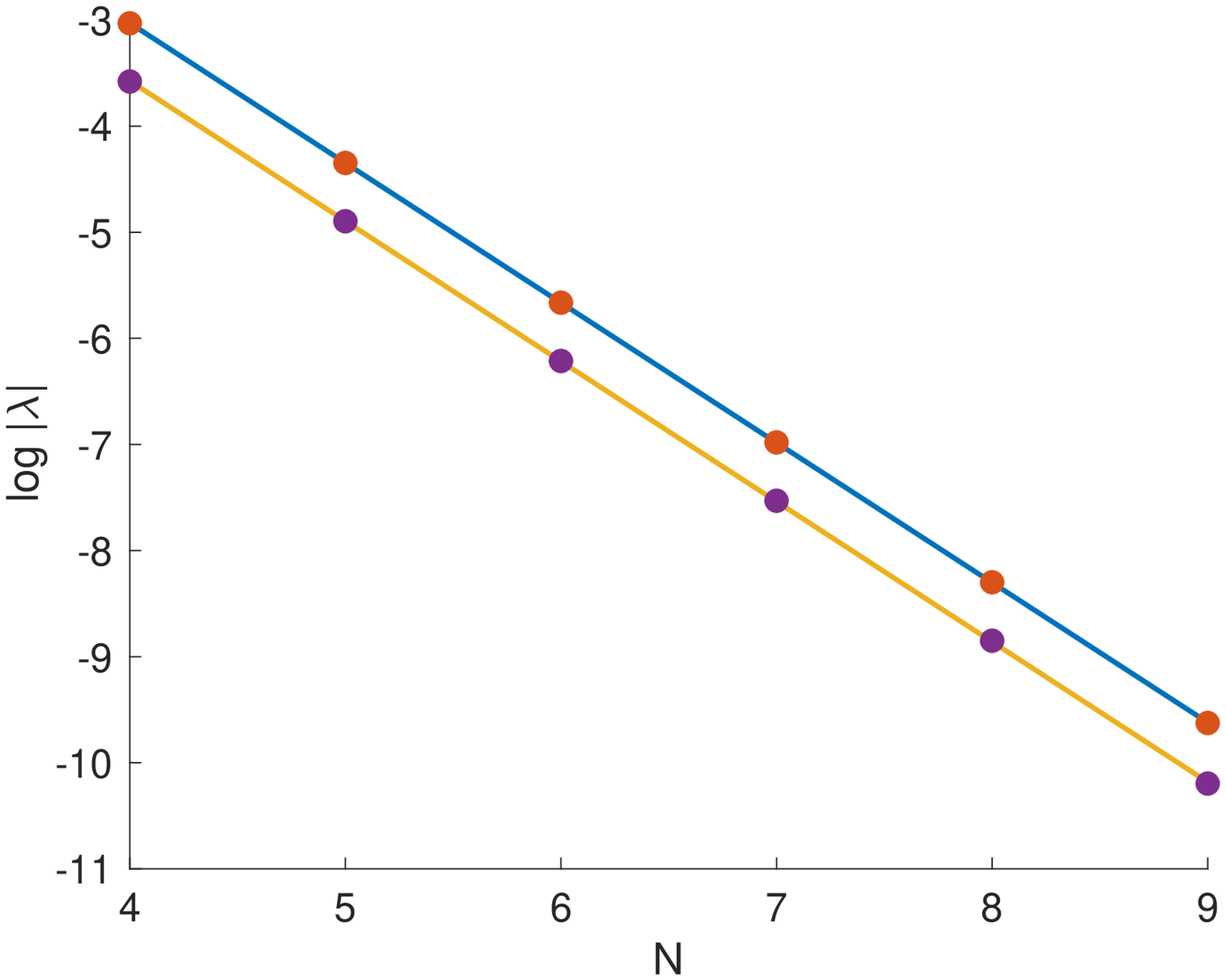}
\includegraphics[width=5cm]{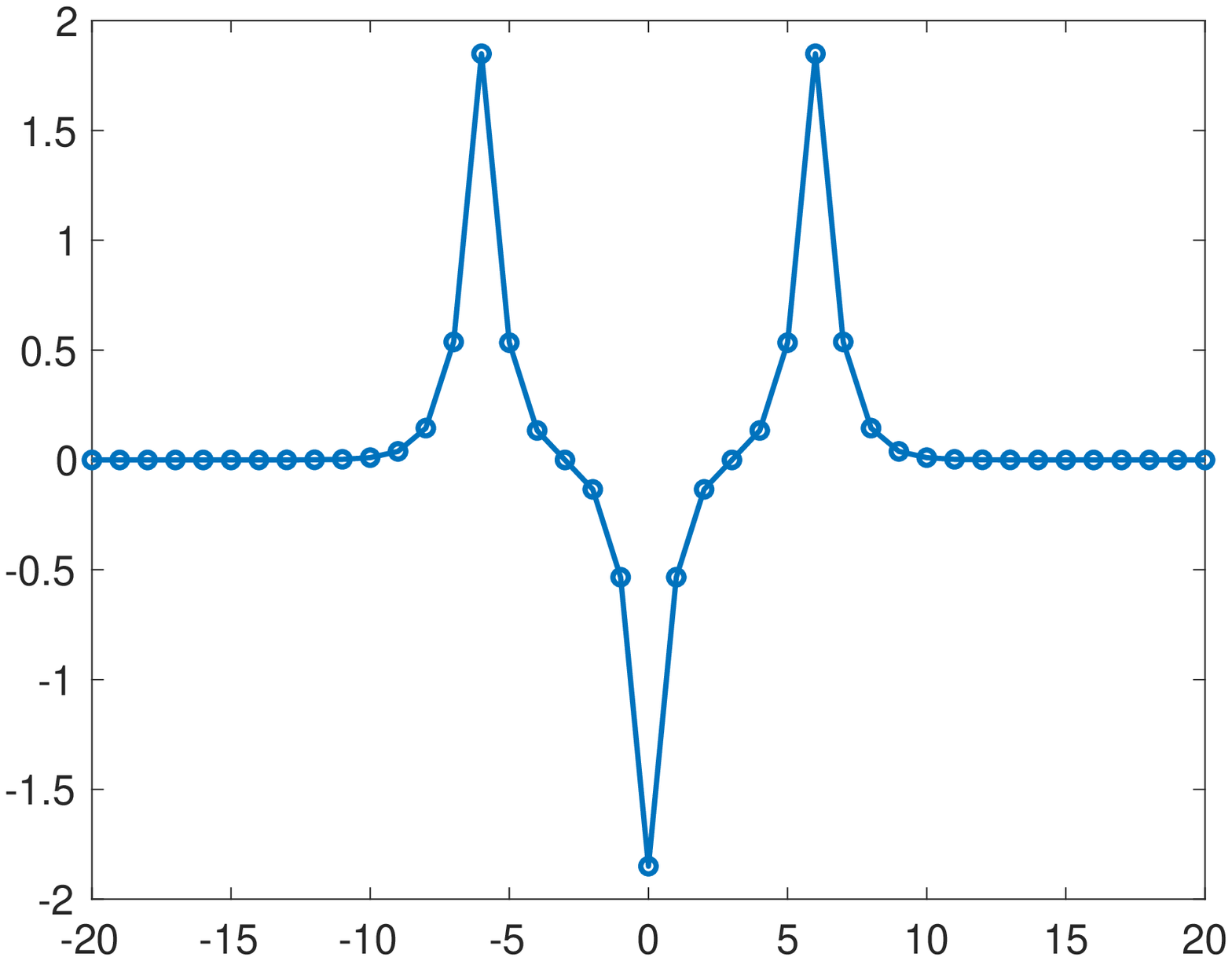}
\includegraphics[width=5cm]{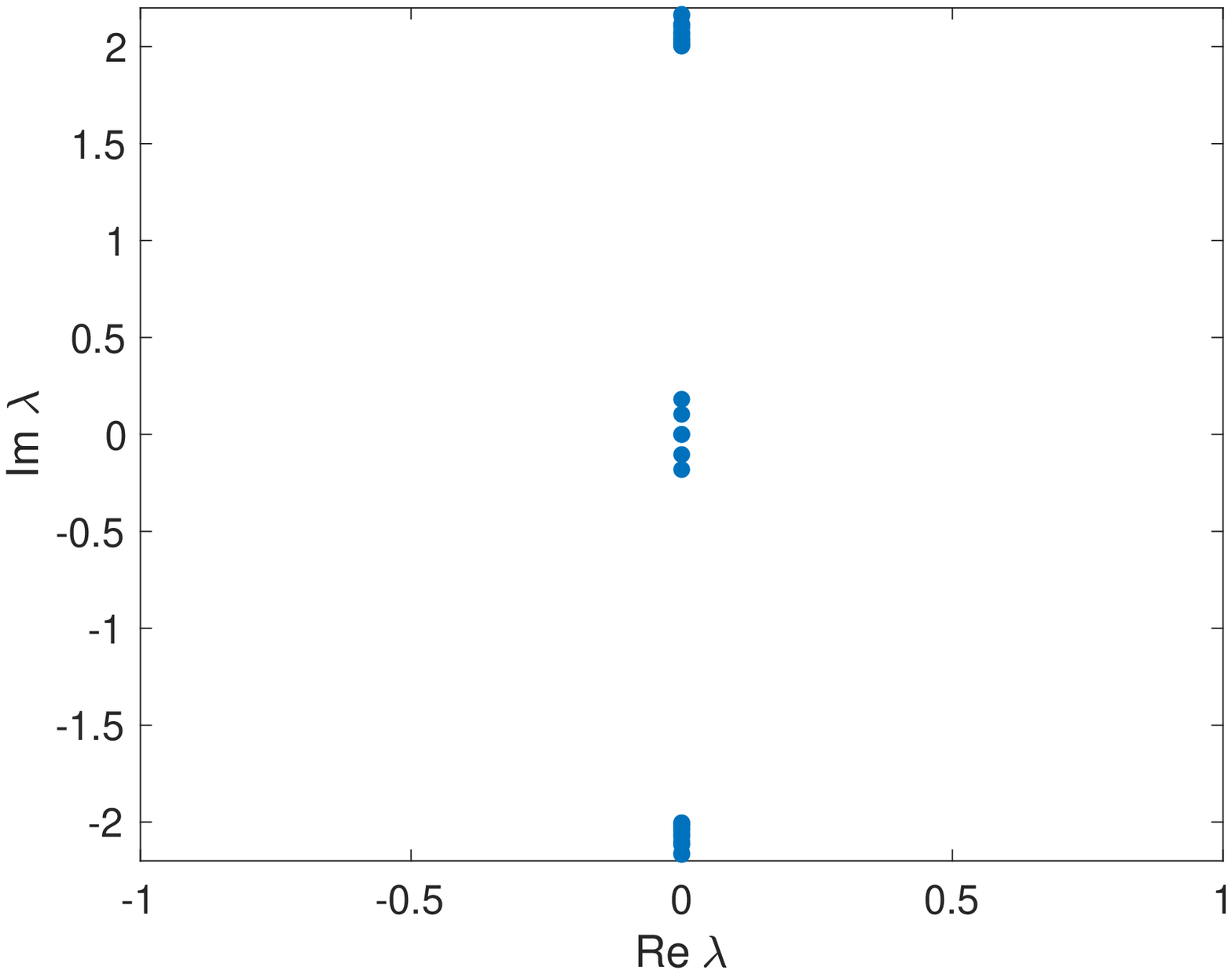}
\includegraphics[width=5cm]{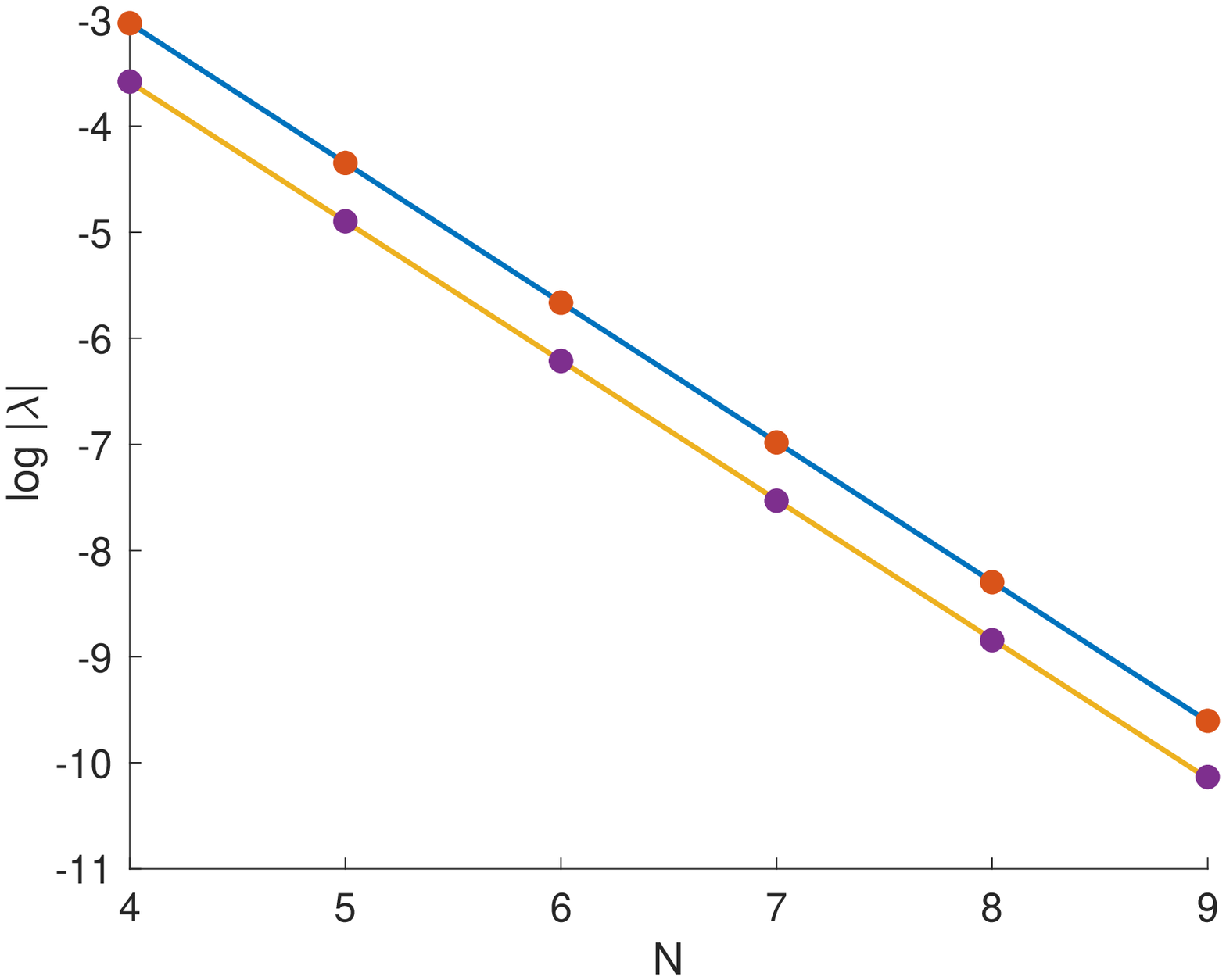}
\includegraphics[width=5cm]{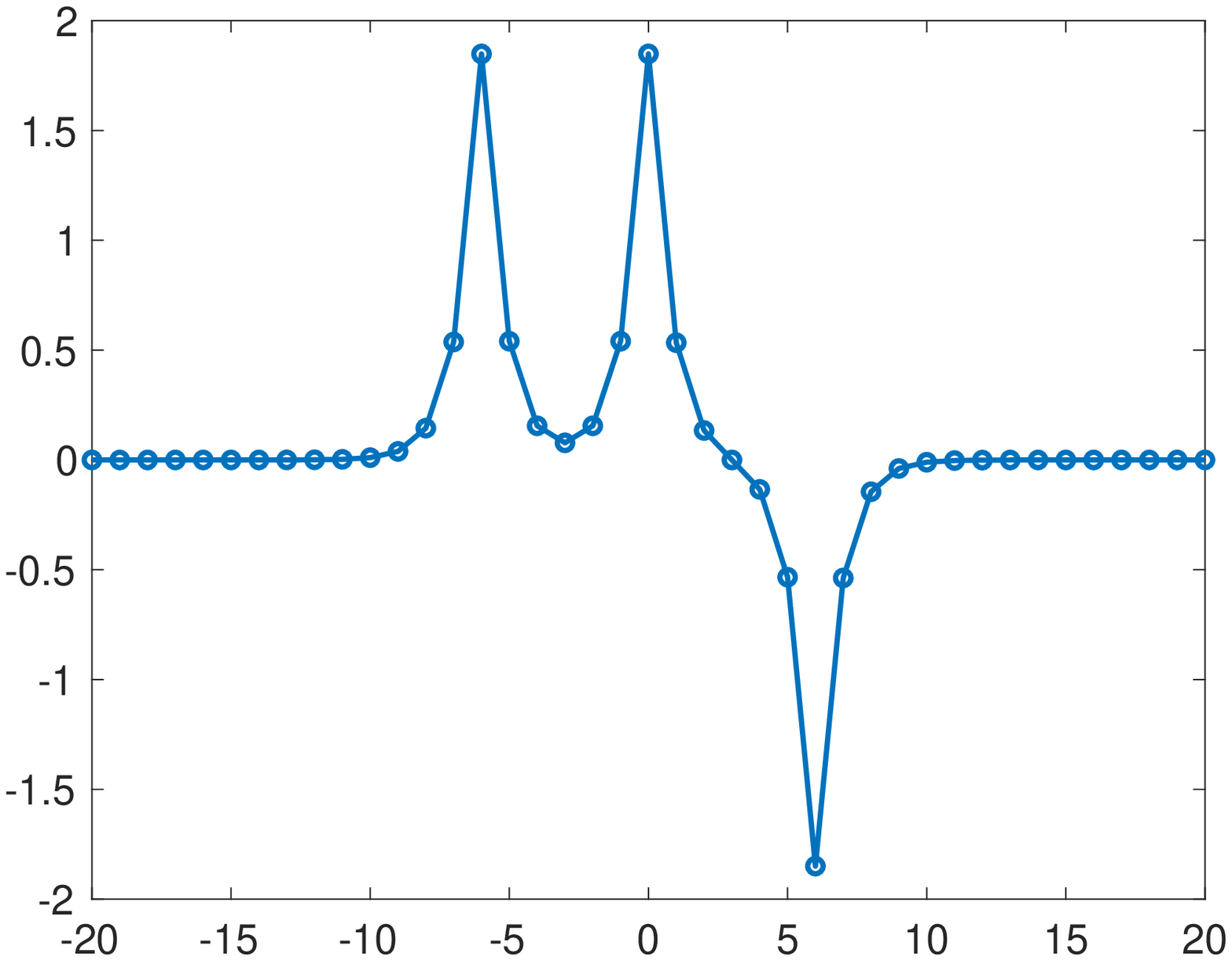}
\includegraphics[width=5cm]{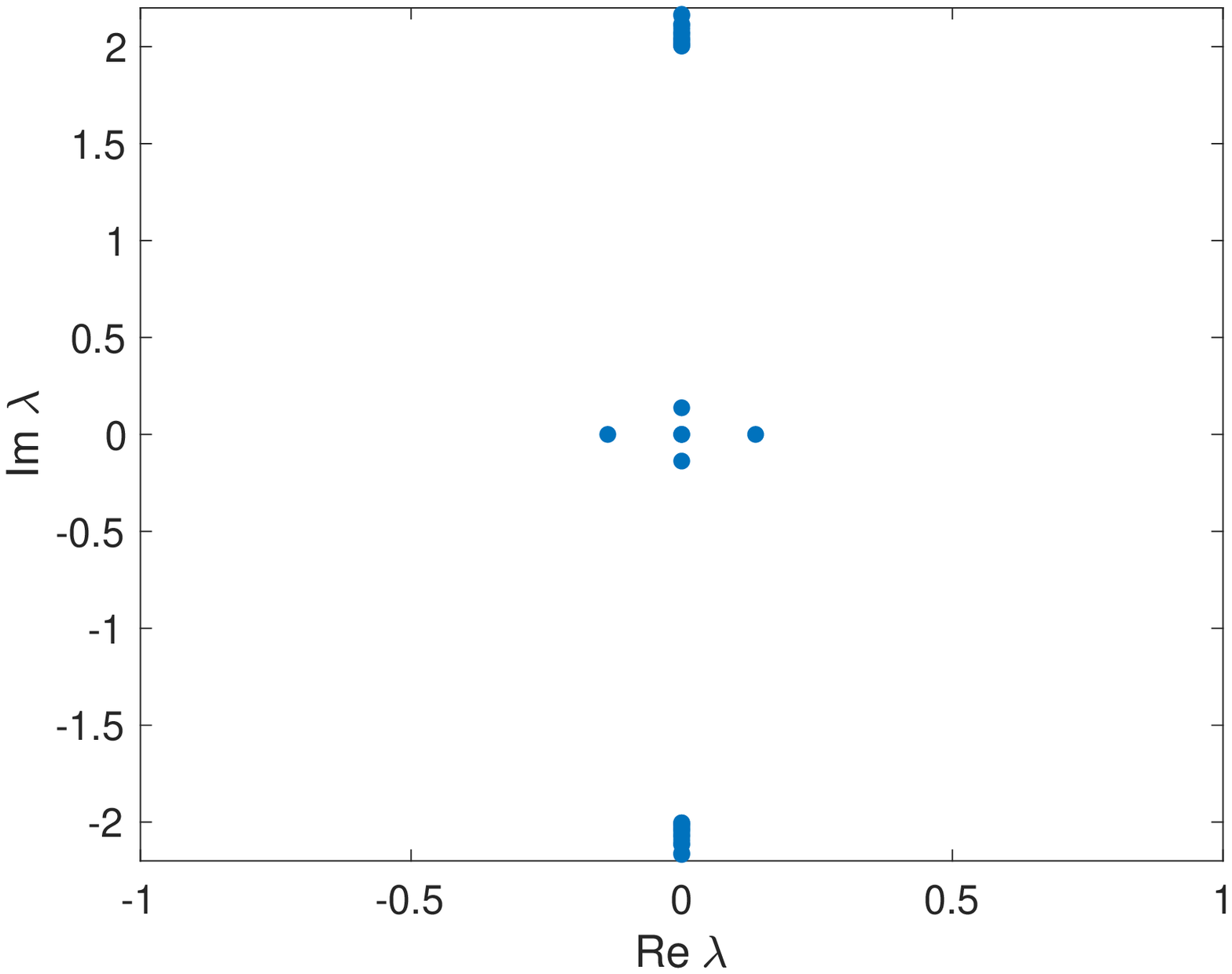}
\includegraphics[width=5cm]{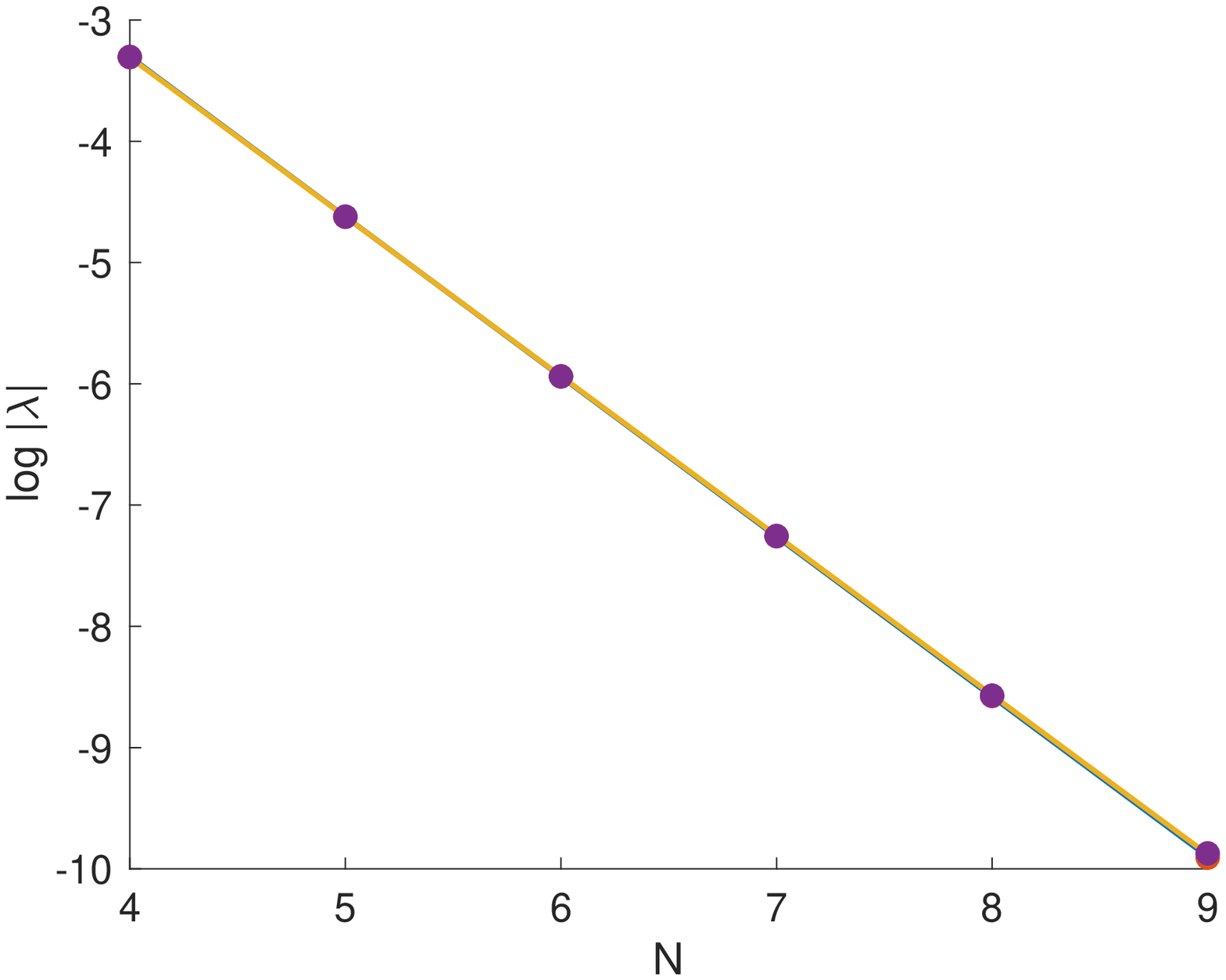}
\caption{Solution profile (left panel), spectral plane eigenvalue pattern (center panel), and plot of $\log(\lambda)$ vs. $N$ with least squares linear regression line (right panel) for the three triple pulse cases: $+++$ (top), $+-+$ (middle), and $++-$ pulses. Parameters $\omega = 2$ and $d = 1.0$.}
\label{fig:eigendecay2}
\end{figure}

We do the same for triple pulses with equal pulse distances in \cref{fig:eigendecay2}. Since the pulse distances are equal, both sets of interaction eigenvalues decay as $r^{-N}$ by equation \cref{3pulseequaleigs} from \cref{DNLSeigcorr}. In the right panel of \cref{fig:eigendecay2}, we plot $\log \lambda$ vs. $N$ for the three triple pulses and construct a least-squares linear regression line. In all three cases, 
namely the in-phase (or $+++$) pulse, the 
out-of-phase (or $+-+$) and finally the 
intermediate/mixed phase case
(or $++-$), the relative error of the slope of the least squares linear regression line is of order $10^{-4}$.  

We can also look at triple pulses with unequal pulse distances $N_1$ and $N_2$. If $N_1 < N_2 < 2 N_1$, then by \cref{DNLSeigcorr2}, there are two pairs of eigenvalues of order $r^{-N_1/2}$ and $r^{-N_2/2}$. We can similarly verify these decay rates numerically.

Finally, we can compute the leading order term in equation \cref{eigsDNLS} and compare that to the numerical result. A value for $\omega$ is chosen, and the single pulse solution $q(n; \omega)$ is constructed numerically using parameter continuation from the anti-continuum limit until the desired coupling parameter $d$ is reached. The terms $b_i$ from the matrix $A$ are computed by using equation \cref{bieq} with the numerically constructed solution $q(n; \omega)$. For the derivative $\partial_\omega q(n; \omega)$, solutions $q(n; \omega + \epsilon)$ and $q(n; \omega - \epsilon)$ are constructed numerically for small $\epsilon$ by parameter continuation from the anti-continuum limit to the same value of $d$. The derivative $\partial_\omega q(n; \omega)$ is computed from these via a centered finite difference method; this is used together with $q(n; \omega)$ to calculate the Melnikov sum $M$. 

\begin{figure}
\centering
\includegraphics[width=7cm]{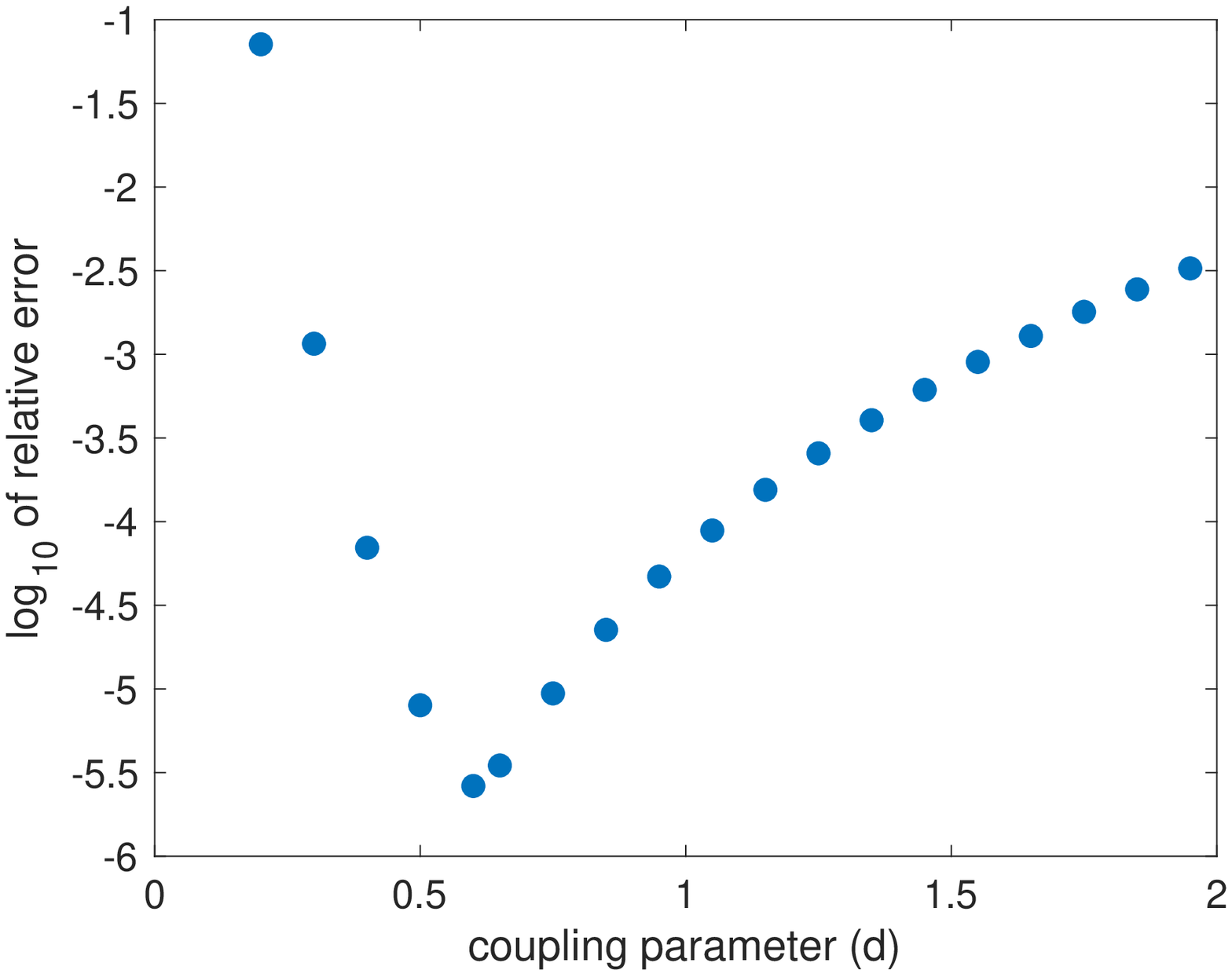}
\includegraphics[width=7cm]{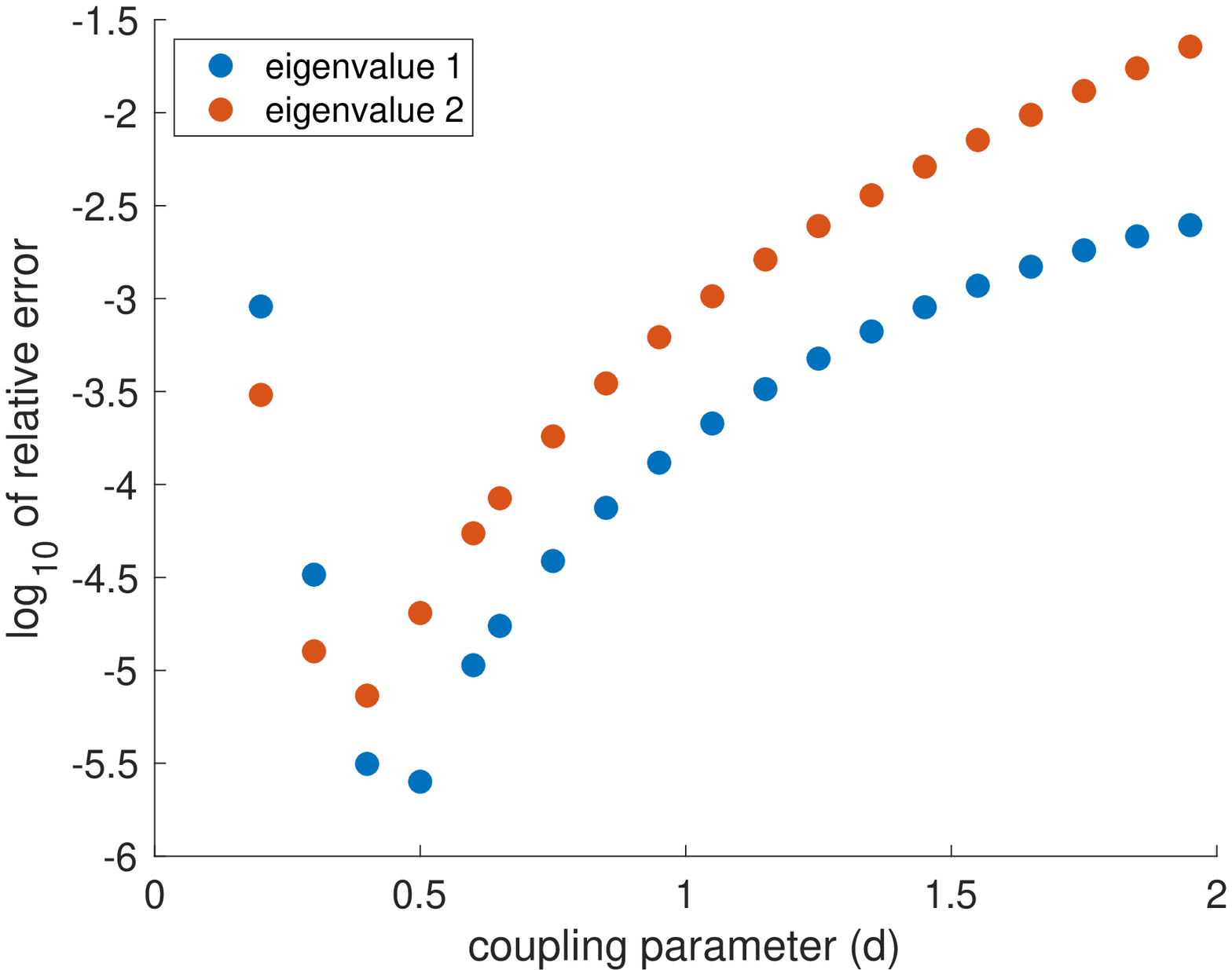}
\caption{Log of relative error of eigenvalues vs. coupling parameter $d$ for double (in phase) 
pulse $++$ ($N_1 = 10$) and triple (out-of-phase) pulse $+-+$ ($N_1 = N_2 = 8)$. $\omega = 2$ in both cases.}
\label{fig:error1}
\end{figure}

First, we consider the case of equal pulse distances. We use the expressions from \cref{DNLSeigcorr} to compute the leading order term for the interaction eigenvalues, and we compare this to the results from Matlab's \texttt{eig} function. In \cref{fig:error1} we fix the inter-pulse distances and plot the log of the relative error of the eigenvalues versus the coupling parameter $d$. For intermediate values of $d$, the relative error is less than $10^{-3}$. Since the results of Theorem \cref{stabilitytheorem} are not uniform in $d$, i.e. they hold for sufficiently large $N$ once $d$ and $\omega$ are chosen, we do not expect to have a nice relationship between the error and $d$. This is furthermore complicated by the fact that additional sources of error arise from  
numerically approximating $b_i$ and $M$. In principle, though, the method (and the asymptotic prediction) yields
satisfactory results except for the vicinity of the 
anti-continuum limit (where the notion of the single
pulse is highly discrete) and the near-continuum limit
(where the role of discreteness is too weak). 
It is interesting to point out that at a ``middle
ground'' between these two limits, namely around $d=0.5$,
we observe the optimal performance of the theoretical
prediction. 

\begin{figure}
\centering
\includegraphics[width=7cm]{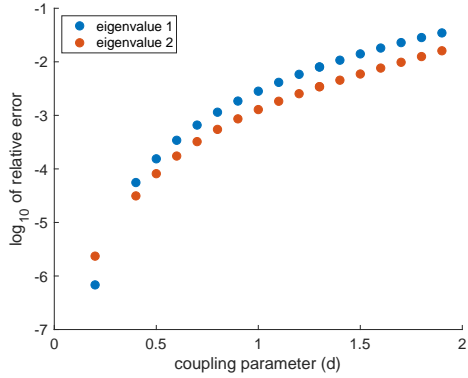}
\caption{Log of relative error of eigenvalues vs. coupling parameter $d$ for triple in-phase pulse $+++$ with unequal pulse distances ($N_1 = 8, N_2 = 6$), $\omega = 2$.}
\label{fig:error2}
\end{figure}

We can also do this for triple pulses with unequal pulse distances. In this case, we use \cref{DNLSeigcorr2} to compute the eigenvalues to leading order. \cref{fig:error2} shows the log of the relative error of the eigenvalues versus the coupling parameter $d$.
For intermediate values of $d$, the relative error is again less than $10^{-3}$. Once again this validates the relevance
of the method especially so for the case of 
intermediate ranges of the coupling parameter $d$.

\section{Conclusions \& future challenges}

In this paper we used Lin's method to construct multi-pulses in discrete systems and to find the small eigenvalues resulting from interaction between neighboring pulses in these structures. In doing so, we are able to extend known results about DNLS to parameter regimes which are further from the anti-continuum limit. In essence, we replace the requirement that the coupling parameter $d$ be small by the condition that the pulses are well separated. This method also allows us to estimate these interaction eigenvalues to a good degree of accuracy for intermediate values of $d$.

The theoretical results we obtained will apply to many other Hamiltonian systems, as long as the coupling by nodes is via the discrete second order centered difference operator $\Delta_2$. Since these restrictions were motivated partly by mathematical convenience, future work could extend these results to a broader class of Hamiltonian systems. Indeed, there exist numerous 
examples worth considering ranging from simpler
ones such as discrete multiple-kink states in 
the discrete sine-Gordon equation~\cite{peyrard},
to settings of first order PDE discretizations
related, e.g., to the Burgers model~\cite{turner}
or even discretizations of third order models
such as the Korteweg-de Vries equation~\cite{ohta}. 

Another direction for future work is characterizing the family of multi-pulse solutions which arises as the coupling parameter $d$ is varied. Recent work~\cite{Jason2019} has investigated stationary, spatially localized patterns in lattice dynamical systems which change as a parameter is varied; the coupling parameter in this case is fixed. In some cases, these patterns exist along a closed bifurcation curves known as an isola. Numerical continuation with AUTO in the coupling parameter $d$ suggests that multi-pulse solutions in DNLS exist on an isola. The parameter $d$ varies over a bounded interval which includes the origin, thus the isola contains solutions to both the focusing and defocusing equation.

A final direction for future works
would concern the consideration of higher dimensional
settings. Here, the interaction between pulses 
would involve the geometric nature of the configuration
they form and the ``line of sight'' between them.
The latter is expected (from the limited observations that
there exist~\cite{alanold}) to determine the nature
of the interaction eigenvalues. Here, however, the
scenarios can also be fundamentally richer as 
coherent states involving topological charge/vorticity
may come into play~\cite{Kevrekidis2009}. In the latter
case, it is less straightforward to identify what the
conclusions may be and considering such more complex
configurations (given also their experimental
observation~\cite{vo3a,vo3b}) may be of particular 
interest.

\section{Proof of existence theorems}

In this section, we will prove \cref{ntmulti} and \cref{transversemulti}. Since the proofs are very similar, we will prove \cref{ntmulti} then state what modifications are necessary for the proof of \cref{transversemulti}. Throughout this section, we will assume \cref{symmetryhyp}, \cref{boundstatehyp}, and \cref{intersectionhyp}. We begin with setting up the exponential dichotomy necessary for the proof. The technique of the proof is very similar to that in \cite{Sandstede1997}.

\subsection{Discrete exponential dichotomy}

First, we define the discrete evolution operator for linear difference equations.

\begin{lemma}[Discrete Evolution Operator]\label{evolop}
Consider the difference equation together with its adjoint
\begin{align}
V(n+1) &= A(n) V(n) \label{diffeqevol} \\
Z(n+1) &= [A(n)^{-1}]^* Z(n) , \label{adjeqevol}
\end{align}
where $n \in \Z$, $V(n) \in R^d$, and the matrix $A(n)$ is invertible for all $n$. Define the discrete evolution operator by
\begin{equation}\label{evol}
\Phi(m, n) = 
\begin{cases}
I & m = n \\
A(m-1) \dots A(n+1) A(n) & m > n \\
A^{-1}(m) \dots A^{-1}(n-2) A^{-1}(n-1) & m < n
\end{cases} \:.
\end{equation}
\begin{enumerate}[(i)]
\item The evolution operators $\Phi$ of \cref{diffeqevol} and $\Psi$ of \cref{adjeqevol} are related by
\begin{equation}\label{adjevol}
\Psi(m, n) = \Phi(n, m)^*.
\end{equation}
\item If $V(n)$ is a solution to \cref{diffeqevol} and $Z(n)$ is a solution to \cref{adjeqevol}, then the inner product $\langle V(n), Z(n) \rangle$ is constant in $n$.
\end{enumerate}

\begin{proof}
For (i), the result holds trivially for $m = n$. For, $m < n$ we have
\begin{align*}
\Psi(m, n) &= A(m)^* \dots A(n-2)^* A(n-1)^* \\
&= [A(n-1) A(n-2) \dots A(m)]^* \\
&= \Phi(n, m)^* .
\end{align*}
The case for $m > n$ is similar.

For (ii), we have
\begin{align*}
\langle V(n+1), Z(n+1) \rangle &= 
\langle A(n) V(n), [A(n)^{-1}]^* Z(n) \rangle \\
&= \langle A(n)^{-1} A(n) V(n), Z(n) \rangle \\
&= \langle V(n), Z(n) \rangle .
\end{align*}
\end{proof}
\end{lemma}

Next, we give a criterion for an exponential dichotomy.

\begin{lemma}[Exponential Dichotomy]\label{dichotomy}
Consider the difference equation
\begin{equation}\label{diffeqdichot}
V(n+1) = A(n) V(n).
\end{equation}
Suppose there exists a constant $r > 1$ and a constant coefficient matrix $A$ such that 
\begin{equation}\label{Anexpdecay}
|A(n) - A| \leq C r^{|n|}
\end{equation}
and $|\lambda| \geq r$ or $|\lambda| \leq 1/r$ for all eigenvalues $\lambda$ of $A$. Then \cref{diffeqdichot} has exponential dichotomies on $Z^\pm$. Specifically, there exist projections $P_\pm^s$ and $P_\pm^u$ defined on $\Z^\pm$ such that the following are true.
\begin{enumerate}[(i)]
\item Let $\Phi(m, n)$ be the evolution operator for \cref{diffeqdichot}. Then 
\begin{equation}\label{projcommute}
P_\pm^{s/u}(m) \Phi(m, n) =  \Phi(m, n) P_\pm^{s/u}(n).
\end{equation}

\item Let $\Phi_\pm^{s/u}(m, n) = \Phi(m, n) P_\pm^{s/u}(n)$ for $m, n \geq 0$ and $m, n \leq 0$ (respectively). Then we have the estimates
\begin{equation*}
\begin{aligned}
|\Phi_+^s(m, n)| \leq C r^{m - n} && 0 \leq n \leq m \\
|\Phi_+^u(m, n)| \leq C r^{n - m} && 0 \leq m \leq n \\
|\Phi_-^s(m, n)| \leq C r^{m - n} && n \leq m \leq 0 \\
|\Phi_-^u(m, n)| \leq C r^{n - m} && m \leq n \leq 0 \:,
\end{aligned}
\end{equation*}
where the evolution operator $\Phi(m, n)$ is defined in \cref{evolop}. 

\item Let $E^{s/u}$ be the stable and unstable eigenspaces of $A$, and let $Q^{s/u}$ the corresponding eigenprojections. Then we have
\begin{align*}
\dim \ran P_\pm^s(n) &= \dim E^s \\
\dim \ran P_\pm^u(n) &= \dim E^u
\end{align*}
and the exponential decay rates
\begin{align}\label{projexpdecay}
| P_\pm^{s/u}(n) - Q^{s/u} | \leq C r^{|n|}.
\end{align}
\end{enumerate}

\begin{proof}
We will consider the problem on $\Z^+$. Since $A$ is constant coefficient and hyperbolic, the difference equation $W(n+1) = A  W(n)$ has an exponential dichotomy on $\R^+$. All the results except for \cref{projexpdecay} follow directly from \cite[Proposition 2.5]{Beyn1997}. Equation \cref{projexpdecay} follows from using the estimate \cref{Anexpdecay} in the proof of \cite[Proposition 2.5]{Beyn1997}.
\end{proof}
\end{lemma}

The last thing we will need is a version of the variation of constants formula for the discrete setting.

\begin{lemma}[Discrete variation of constants]\label{VOC}
The solution $V(n)$ to the initial value problem
\begin{equation*}
\begin{aligned}
V(n+1) &= A(n) V(n) + G(V(n), n) \\
V(n_0) &= V_{n_0}
\end{aligned}
\end{equation*}
can be written in summation form as 
\begin{equation}\label{VOCformula}
V(n) = 
\begin{cases}
V_{n_0} & n = n_0 \\
\Phi(n, n_0) V_{n_0} + \sum_{j = n_0}^{n-1} \Phi(n, j+1) G(V(j), j)) & n > n_0 \\
\Phi(n, n_0) V_{n_0} - \sum_{j = n}^{n_0-1} \Phi(n, j+1) G(V(j), j)) & n < n_0 
\end{cases}\:.
\end{equation}
\begin{proof}
For $n = n_0 + 1$,
\[
V(n_0 + 1) = A(n_0) V(n_0) + G(V(n_0), n_0) = \Phi(n_0+1, n_0) V(n_0) + \Phi(n_0, n_0) G(V(n_0), n_0).
\]
Iterate this to get the result for $n > n_0$. The case for $n < n_0$ is similar.
\end{proof}
\end{lemma}

\subsection{Fixed point formulation}

To find a solution to the system of equations \cref{Usystem}, we will rewrite the system as a fixed point problem. First, we expand $F$ in a Taylor series about $T(\theta_i) Q(n)$ to get
\begin{align*}
F(U_i^\pm(n)) &= F(T(\theta_i) Q(n) + \tilde{Q}_i^\pm(n)) \\
&= F(T(\theta_i) Q(n)) + D F(T(\theta_i) Q(n)) \tilde{Q}_i^\pm(n) + G(\tilde{Q}_i^\pm(n)) \\
&= T(\theta_i)DF(Q(n))T(\theta_i)^{-1} Q(n)) \tilde{Q}_i^\pm + G(\tilde{Q}_i^\pm(n)),
\end{align*}
where $G(\tilde{Q}_i^\pm(n)) = \mathcal{O}(|\tilde{Q}_i^\pm|^2)$ with $G(0) = 0$ and $DG(0) = 0$, and we used the symmetry relation \cref{symmetryrel} in the last line. Finally, let
\begin{align} \label{didef}
d_i &= T(\theta_{i+1}) Q(-N_i^-) - T(\theta_i) Q(N_i^+).
\end{align}
Substituting these into \cref{Usystem}, we obtain the following system of equations for the remainder functions $\tilde{Q}_i^\pm$.
\begin{align}
\tilde{Q}_i^\pm(n+1) &= T(\theta)D F(Q(n))T(\theta)^{-1} \tilde{Q}_i^\pm(n) + G(\tilde{Q}_i^\pm(n)) \label{Wsystem1} \\
\tilde{Q}_i^+(N_i^+) - W_{i+1}^-(-N_i^-) &= d_i \label{Wsystem2} \\
\tilde{Q}_i^+(0) - \tilde{Q}_i^-(0) &= 0 . \label{Wsystem3}
\end{align}

Next, we look at the variational and adjoint variational equations associated with \cref{diffeq}, which are
\begin{align}
V(n+1) &= D F(Q(n)) V(n) \label{vareq} \\
Z(n+1) &= [D F(Q(n))^*]^{-1} Z(n) . \label{adjvareq} 
\end{align}
The variational equation \cref{vareq} has a bounded solution $T'(0) Q(n)$, thus we can decompose the tangent spaces to $W^s(0)$ and $W^u(0)$ at $Q(0)$ as
\begin{align*}
T_{Q(0)} W^u(0) &= Y^- \oplus \R T'(0) Q(0) \\
T_{Q(0)} W^s(0) &= Y^+ \oplus \R T'(0) Q(0) .
\end{align*}
The adjoint variational equation also has a unique bounded solution $Z_1(n)$ given by \cref{adjvarsol}. By \cref{evolop}, $Z_1(0) \perp T'(0) Q(0) \oplus Y^- \oplus Y^+$, thus we can decompose $\R^d$ as
\begin{equation}\label{nontdecomp}
\R^d = \R T'(0) Q(0) \oplus Y^+ \oplus Y^- \oplus \R Z_1(0) .
\end{equation}
Since $T(\theta)$ is unitary, we also have the decomposition
\begin{equation}\label{nontdecompT}
\R^d = \R T(\theta_i) T'(0) Q(0) \oplus T(\theta_i) Y^+ \oplus T(\theta_i) Y^- \oplus \R T(\theta_i) Z_1(0).
\end{equation}
Finally, since perturbations in the direction of $T(\theta_i) T'(0) Q(0)$ are handled by the symmetry parameter $\theta_i$, we may without loss of generality choose $\tilde{Q}_i^\pm$ so that 
\begin{equation}\label{W0loc}
\tilde{Q}_i^\pm(0) \in T(\theta_i) Y^+ \oplus T(\theta_i) Y^- \oplus \R T(\theta_i) Z_1(0).
\end{equation}

Let $\Phi(m, n; \theta)$ be the evolution operator for
\begin{equation}\label{Vtheta}
V(n+1; \theta) = T(\theta) D F(Q(n)) T(\theta)^{-1} V(n; \theta) .
\end{equation}
We note that since $T'(0)$ commutes with $T(\theta)$, $T'(0)Q(n)$ is a solution to \cref{Vtheta}. Using \cref{symmetryrel}, the evolution operators are related to those for $\theta = 0$ by
\begin{equation}\label{evoloptheta}
\Phi(m, n; \theta) = T(\theta)\Phi(m, n; 0)T(\theta)^{-1}.
\end{equation}
Since $T(\theta) D F(Q(n)) T(\theta)^{-1}$ decays exponentially to $T(\theta) D F(0) T(\theta)^{-1}$ and $D F(0)$ is hyperbolic, equation \cref{Vtheta} has exponential dichotomies on $\Z^+$ and $\Z^-$ by \cref{dichotomy}, and we note that the estimates from \cref{dichotomy} do not depend on $\theta$. Let $P_{s/u}^\pm(m; \theta)$ and $\Phi_{s/u}^\pm(m, n; \theta)$ be the projections and evolutions for this exponential dichotomy on $\Z^\pm$. The projections $P_{s/u}^\pm(m; \theta)$ are related to those for $\theta = 0$ by
\[
P_{s/u}^\pm(m; \theta) = T(\theta)P_{s/u}^\pm(m; 0)T(\theta)^{-1}
\]
Finally, let $E^s(\theta)$ and $E^u(\theta)$ be the stable and unstable eigenspaces of $T(\theta) D_U F(0) T(\theta)^{-1}$, and let $P_0^s(\theta)$ and $P_0^u(\theta)$ be the corresponding eigenprojections.

Next, as in \cite{Sandstede1997} and \cite{Knobloch2000}, we write equation \cref{Wsystem1} in fixed-point form using the discrete variation of constants formula \cref{VOCformula} together with projections on the stable and unstable subspaces of the exponential dichotomy.
\begin{equation}\label{FPeqs1}
\begin{aligned}
\tilde{Q}_i^-(n) &= 
\Phi_s^-(n, -N_{i-1}^-; \theta_i) a_{i-1}^- + \Phi_u^-(n, 0; \theta_i) b_i^-  \\
&+ \sum_{j = -N_{i-1}^-}^{n-1} \Phi_s^-(n, j+1; \theta_i) G_i^-(\tilde{Q}_i^-(j)) - \sum_{j = n}^{-1} \Phi_u^-(n, j+1; \theta_i) G_i^-(\tilde{Q}_i^-(j)) \\
\tilde{Q}_i^+(n) &= \Phi_u^+(n, N_i^+; \theta_i) a_i^+ + \Phi_s^+(n, 0; \theta_i) b_i^+ \\
&+ \sum_{j = 0}^{n-1} \Phi_s^+(n, j+1; \theta_i) G_i^+(\tilde{Q}_i^+(j)) 
- \sum_{j = n}^{N_i^+-1} \Phi_u^+(n, j+1; \theta_i) G_i^+(\tilde{Q}_i^+(j)),
\end{aligned}
\end{equation}
where $\tilde{Q}_i^-(n) \in \ell^\infty([-N_{i-1}^-, 0])$, $\tilde{Q}_i^+(n) \in \ell^\infty([0, N_i^+])$, and the sums are defined to be $0$ if the upper index is smaller than the lower index. For the initial conditions, 
\begin{enumerate}
\item $a_i^- \in E^s(\theta_i)$, $a_i^+ \in E^u(\theta_i)$, and $a_0^- = a_m^+ = 0$.
\item $b_i^+ \in T(\theta_i) Y^+$ and $b_i^- \in T(\theta_i) Y^-$.
\end{enumerate}
We note that we do not need to include a component in $T'(0) Q(0)$ in $b_i^\pm$, since that direction is handled by the symmetry parameter $\theta_i$.

Since we wish to construct a homoclinic orbit to the rest state at 0, we take the initial conditions $a_0^- = 0$ and $a_m^+ = 0$. For these cases, the fixed point equations are given by
\begin{align*}
W_1^-(n) &= \Phi_u^-(n, 0; \theta_i) b_i^- 
+ \sum_{j = -\infty}^{n-1} \Phi_s^-(n, j+1; \theta_i) G_i^-(\tilde{Q}_i^-(j)) - \sum_{j = n}^{-1} \Phi_u^-(n, j+1; \theta_i) G_i^-(\tilde{Q}_i^-(j)) \\
W_m^+(n) &= \Phi_s^+(n, 0; \theta_i) b_i^+ 
+ \sum_{j = 0}^{n-1} \Phi_s^+(n, j+1; \theta_i) G_i^+(\tilde{Q}_i^+(j)) 
- \sum_{j = n}^\infty \Phi_u^+(n, j+1; \theta_i) G_i^+(\tilde{Q}_i^+(j)),
\end{align*}
where the infinite sums converge due to the exponential dichotomy. 

\subsection{Inversion}

As in \cite{Sandstede1997}, we will solve equations \cref{Wsystem1}, \cref{Wsystem2}, and \cref{Wsystem3} in stages. In the first lemma of this section, we solve equation \cref{Wsystem1} for $\tilde{Q}_i^\pm$. 

\begin{lemma}\label{inv1}
There exist unique bounded functions $\tilde{Q}_i^\pm(n)$ such that equation \cref{Wsystem1} is satisfied. These solutions depend smoothly on the initial conditions $a_i^\pm$ and $b_i\pm$, and we have the estimates
\begin{equation}\label{Wipmest}
\begin{aligned}
\|\tilde{Q}_i^-\| &\leq C (|a_{i-1}^-| + |b_i^-|) \\
\|\tilde{Q}_i^+\| &\leq C (|a_i^+| + |b_i^+| ) \:.
\end{aligned}
\end{equation}
For the interior pieces, we have the piecewise estimates
\begin{equation}\label{Wipiecewise}
\begin{aligned}
|\tilde{Q}_i^-(n)| &\leq C (r^{-(N_{i-1}^- + n)}|a_{i-1}^-| + r^n|b_i^-|) && n \in [-N_{i-1}^-, 0] \\
|\tilde{Q}_i^+(n)| &\leq C (r^{-(N_i^+ - n)}|a_i^+| + r^{-n}|b_i^+| ) && n \in [0, N_i^+] \:.
\end{aligned}
\end{equation}

\begin{proof}
First, we show that the RHS of the fixed point equations \cref{FPeqs1} defines a smooth map from $\ell^\infty$ (on the appropriate interval) to itself. For the $\tilde{Q}_i^-$, we have
\begin{align}\label{inv1est1}
|\Phi_s^-(n, &-N_{i-1}^-; \theta_i) a_{i-1}^-| + |\Phi_u^-(n, 0; \theta_i) b_i^-| \leq C ( |a_{i-1}^-| + |b_i^-|) 
\end{align}
and
\begin{align*}
\left| \sum_{j = -N_{i-1}^-}^{n-1} \Phi_s^-(n, j+1; \theta_i) G_i^-(\tilde{Q}_i^-(j))\right| + \left|\sum_{j = n}^{-1} \Phi_u^-(n, j+1; \theta_i) G_i^-(\tilde{Q}_i^-(j))\right| 
\leq C \|\tilde{Q}_i^-\|_{\ell^\infty([-N_{i-1}, 0])}^2 \:,
\end{align*}
both of which are independent of $n$. Define the map
$K_i^-: \ell^\infty([-N_{i-1}, 0]) \times E^s \times Y^- \rightarrow \ell^\infty([-N_{i-1}, 0])$ by
\begin{align}\label{inv1map}
K_i^-(\tilde{Q}_i^-(n), &a_{i-1}^-, b_i^-) = \tilde{Q}_i^-(n) - \Phi_s^-(n, -N_{i-1}^-; \theta_i) a_{i-1}^- - \Phi_u^-(n, 0; \theta_i) b_i^-  \\
&- \sum_{j = -N_{i-1}^-}^{n-1} \Phi_s^-(n, j+1; \theta_i) G_i^-(\tilde{Q}_i^-(j)) + \sum_{j = n}^{-1} \Phi_u^-(n, j+1; \theta_i) G_i^-(\tilde{Q}_i^-(j)) . \nonumber
\end{align}
Since 0 is an equilibrium, $K(0, 0, 0) = 0$. It is straightforward to show that the Fr\'echet derivative of $K_i^-$ with respect to $\tilde{Q}_i^-$ at $(\tilde{Q}_i^-(n), a_{i-1}^-, b_i^-) = (0, 0, 0)$ is a Banach space isomorphism on $l^\infty([-N_{i-1}, 0])$. Thus we can solve for $\tilde{Q}_i^-(x)$ in terms of $(a_{i-1}^-, b_i^-)$ using the IFT. This dependence is smooth, since the map $K_i^-$ is smooth. The estimate \cref{Wipmest} on $\tilde{Q}_i^-$ comes from \cref{inv1est1}, since the terms in \cref{inv1map} involving sums are quadratic in $\tilde{Q}_i^\pm$. The case for $\tilde{Q}_i^+$ is similar. It is not hard to obtain the piecewise estimates \cref{Wipiecewise} for the interior pieces.
\end{proof}
\end{lemma}

Next, we use the center matching conditions at $N_i^\pm$ to solve equation \cref{Wsystem2}. This will give us the initial conditions $a_i^\pm$.

\begin{lemma}\label{inv2}
For $i = 1, \dots m-1$ there is a unique pair of initial conditions $(a_i^+, a_i^-) \in E^u(\theta_i) \times E^s(
\theta_i)$ such that the matching conditions \cref{Wsystem2} are satisfied. $(a_i^+, a_i^-)$ depends smoothly on $(b_i^+, b_{i+1}^-)$, and we have the following expressions for $a_i^-$ and $a_i^+$. 
\begin{equation}\label{aipmest}
\begin{aligned}
a_i^+ &= P_0^u(\theta_i) d_i + \tilde{a}_i^+ \\
a_i^- &= -P_0^s(\theta_i) d_i + \tilde{a}_i^- ,
\end{aligned}
\end{equation}
where 
\begin{equation}\label{tildeaest}
\tilde{a}_i^\pm = \mathcal{O}(r^{-N}(|b_i^+|+|b_{i+1}^-|) + |b_i^+|^2+|b_{i+1}^-|^2) .
\end{equation}
In terms of $Q(\pm N_i^\pm)$, we can write \cref{aipmest} as 
\begin{equation}\label{aipmexp}
\begin{aligned}
a_i^- &= T(\theta_i) Q(N_i^+) + \tilde{a}_i^- + \mathcal{O}(r^{-2N}) \\
a_i^+ &= T(\theta_{i+1}) Q(-N_i^-) + \tilde{a}_i^+ + \mathcal{O}(r^{-2N}).
\end{aligned}
\end{equation}

\begin{proof}
Evaluating the fixed point equations \cref{FPeqs1} at $\pm N_i^\pm$ and subtracting, solving equation \cref{Wsystem2} is equivalent to solving $K_i(a_i^+, a_i^-, b_i^+, b_{i+1}^-) = 0$, where $K_i: E^s \times E^u \times Y^+ \times Y^- \rightarrow \R^d$ is defined by
\begin{align*}
K_i(a_i^+, &a_i^-, b_i^+, b_{i+1}^-) \\
&= a_i^+ - a_i^- - d_i + (P_u^+(N_i^+; \theta_i) - P_0^u) a_i^+ - (P_s^-(-N_i^-; \theta_{i+1}) - P_0^s) a_i^- \\
&+ \Phi_s^+(N_i^+, 0; \theta_i) b_i^+ - \Phi_u^-(-N_i^-, 0; \theta_{i+1}) b_{i+1}^- \\
&+ \sum_{j = 0}^{N_i^+-1} \Phi_s^+(N_i^+, j+1; \theta_i) G_i^+(\tilde{Q}_i^+(j; a_i^+, b_i^+)) \\
&+ \sum_{j = -N_i^-}^{-1} \Phi_u^-(-N_i^-, j+1; \theta_{i+1}) G_i^-(\tilde{Q}_{i+1}^-(j; a_i^-, b_{i+1}^-)),
\end{align*}
and we substituted $W_{i+1}^-(n; a_i^-, b_{i+1}^-)$ and $\tilde{Q}_i^+(n; a_i^+, b_i^+)$ from \cref{inv1}. Next, we note that $K_i(0,0,0,0) = 0$ and that 
\begin{align*}
\frac{\partial}{\partial a_i^-} K_i(0, 0, 0, 0) &= -1 + \mathcal{O}(r^{-N_i^-}) \\
\frac{\partial}{\partial a_i^+} K_i(0, 0, 0, 0) &= 1 + \mathcal{O}(r^{-N_i^+}),
\end{align*}
since the derivatives of the terms in $K_i$ involving sums will be 0 since $G_i^\pm$ is quadratic in $\tilde{Q}_i^\pm$, thus quadratic order in $a_i^\pm$ by \cref{inv1}. For sufficiently large $N$, $D_{a_i^\pm} K(0, 0, 0, 0, 0)$ is invertible in a neighborhood of $(0, 0, 0, 0, 0)$. Thus, since $(a_i^+, a_i^-) \in E^s(\theta_i) \oplus E^u(\theta_i) = \R^d$, we can use the IFT to solve for $a_i^\pm$ in terms of $(b_i^+$, $b_{i+1}^-)$ for $(b_i^+, b_{i+1}^-)$ sufficiently small.

To get the estimates on and expressions for $a_i^\pm$, we project $K_i(a_i^+, a_i^-, b_i^+, b_{i+1}^-) = 0$ onto $E^s(\theta_i)$ and $E^u(\theta_i)$ in turn to get 
\begin{align*}
a_i^+ &= P_0^u(\theta_i) d_i + \mathcal{O}(r^{-N}(|b_i^+|+|b_{i+1}^-|) + |b_i^+|^2+|b_{i+1}^-|^2) \\
a_i^- &= -P_0^s(\theta_i) d_i + \mathcal{O}(r^{-N}(|b_i^+|+|b_{i+1}^-|) + |b_i^+|^2+|b_{i+1}^-|^2),
\end{align*}
which we can write in the form \cref{aipmest} with estimates \cref{tildeaest}. 

To write these in terms of $Q(\pm N_i^\pm)$, we note that
\begin{align*}
P_0^s(\theta_i) T(\theta_i) Q(N_i^+) &= (P_0^s(\theta_i) - P_s^+(N_i^+; \theta_i)) T(\theta_i) Q(N_i^+) + P_s^+(N_i^+; \theta_i) T(\theta_i) Q(N_i^+) \\
&= T(\theta_i)(P_0^s(0) - P_s^+(N_i^+; 0)) Q(N_i^+) + P_s^+(N_i^+; \theta_i) T(\theta_i) Q(N_i^+) \\
&= T(\theta_i) Q(N_i^+) + \mathcal{O}(r^{-2N}),
\end{align*}
where in the second line we used \cref{projexpdecay}. Similarly, we can show that 
\begin{align*}
P_0^s T(\theta_{i+1}) Q(-N_i^-) &= \mathcal{O}(r^{-2N}).
\end{align*}
Substituting these into \cref{aipmest} we obtain \cref{aipmexp}.
\end{proof}
\end{lemma}

It only remains to satisfy \cref{Wsystem3}, which is the jump condition at 0. We will not in general be able to solve equation \cref{Wsystem3}. In the next lemma, we will solve for the initial conditions $b_i^\pm$. This will give us a unique solution which will generically have $m$ jumps in the direction of $T(\theta) Z_1(0)$. We will obtain a set of $m$ jump conditions in the direction of $T(\theta) Z_1(0)$ which will depend on the symmetry parameters $\theta_i$. Satisfying the jump conditions, which solves \cref{Wsystem3}, can be accomplished by adjusting the symmetry parameters.

Recall that for all $\theta \in \R$ we have the decomposition
\[
\R^d = \R T(\theta) T'(0)Q(0) \oplus T(\theta) Y^+ \oplus T(\theta) Y^- \oplus \R T(\theta) Z_1(0).
\]
Projecting in these directions, we can write \cref{Wsystem3} as the system of equations
\begin{align}
P_{T(\theta_i)T'(0)Q(0)}\left( \tilde{Q}_i^+(0) - \tilde{Q}_i^-(0) \right) &= 0 \label{jumpS1} \\
P_{T(\theta_i)Y^+ \oplus T(\theta_i)Y^-}\left( \tilde{Q}_i^+(0) - \tilde{Q}_i^-(0) \right) &= 0 \label{jumpnonZ} \\
P_{\R T(\theta_i)Z_1(0)} \left( \tilde{Q}_i^+(0) - \tilde{Q}_i^-(0) \right) &= 0 . \label{jumpZ}
\end{align}
Since $\tilde{Q}_i^\pm(0) \in Y^+ \oplus Y^- \oplus \R Z_1(0)$, equation \cref{jumpS1} is automatically satisfied. Since $b_i^+ \in T(\theta) Y^+$ and $b_i^- \in T(\theta) Y^-$, we will be able to satisfy \cref{jumpnonZ} by solving for the $b_i^\pm$, which we do in the following lemma.

\begin{lemma}\label{inv3nt}
For $i = 1, \dots m$ there is a unique pair of initial conditions $(b_i^-, b_i^+) \in T(\theta_i) Y^- \times T(\theta_i) Y^+$ such that \cref{jumpnonZ} is satisfied. We have the uniform bound
\begin{equation}\label{bboundnt}
b = \mathcal{O}(r^{-2N}).
\end{equation}

\begin{proof}
For convenience, let $X_i = T(\theta_i)Y^+ \oplus T(\theta_i)Y^-$. Evaluating the fixed point equations \cref{FPeqs1} at 0, subtracting, and applying the projection $P_{X_i}$ to both sides, we have 
\begin{align*}
P_{X_i}( \tilde{Q}_i^+(0) &- \tilde{Q}_i^-(0)) = b_i^+ - b_i^- 
+ P_{X_i}( \Phi_u^+(0, N_i^+; \theta_i) a_i^+) - P_{X_i}( \Phi_s^-(0, -N_{i-1}^-; \theta_i) a_{i-1}^-)  \\
&- P_{X_i} \left( \sum_{j = 0}^{N_i^+-1} \Phi_u^+(0, j+1; \theta_i) G_i^+(\tilde{Q}_i^+(j)) 
- \sum_{j = -N_{i-1}^-}^{-1} \Phi_s^-(0, j+1; \theta_i) G_i^-(\tilde{Q}_i^-(j)) \right) .
\end{align*}
Next, substitute $\tilde{Q}_i^\pm$ from \cref{inv1} and $a_i^\pm$ from \cref{inv2}. Define the spaces
\begin{align}\label{spaceYt}
Y &= \bigoplus_{i=1}^m (T(\theta) Y^+ \oplus T(\theta) Y^-) = \bigoplus_{i=1}^m \R^d \\
Z &= \bigoplus_{i=1}^{m-1} \R^d .
\end{align}
Let $b = (b_1^+, b_1^-, \dots, b_m^+, b_m^-) \in Y$ and $d = (d_1, \dots, d_{m-1}) \in Z$. Define the function $K: Y \times Z \rightarrow Y$ component-wise by
\begin{align*}
K_i(b, d) &= 
 b_i^+ - b_i^- + P_{X_i} \left( \Phi_u^+(0, N_i^+; \theta_i) P_0^u d_i + \Phi_s^-(0, -N_{i-1}^-; \theta_i) P_0^s d_{i-1} \right) \\
&+ P_{X_i}\left( \Phi_u^+(0, N_i^+; \theta_i) \tilde{a}_i^+(b_i^+, b_{i+1}^-) 
- \Phi_s^-(0, -N_{i-1}^-; \theta_i) \tilde{a}_{i-1}^-(b_{i-1}^+, b_i^-) \right) \\
&- P_{X_i} \sum_{j = 0}^{N_i^+-1} \Phi_u^+(0, j+1; \theta_i) G_i^+(\tilde{Q}_i^+(j; b_i^+, b_{i+1}^-)) \\
&- P_{X_i} \sum_{j = -N_{i-1}^-}^{-1} \Phi_s^-(0, j+1; \theta_i) G_i^-(\tilde{Q}_i^-(j; b_{i-1}^+, b_i^-)),
\end{align*}
where $d_0 = d_m = 0$, and we have indicated the dependencies on the $b_i^\pm$. Using the estimates from \cref{inv1} and \cref{inv2}, $K(0, 0) = 0$. For the partial derivatives with respect to $b_i^\pm$, we have
\begin{align*}
\frac{\partial}{\partial b_i^+}K_i(0) &= 1 + \mathcal{O}(r^{-N})  \\
\frac{\partial}{\partial b_i^-}K_i(0) &= -1 + \mathcal{O}(r^{-N}) \\
\frac{\partial}{\partial b_{i-1}^+}K_i(0),
\frac{\partial}{\partial b_{i+1}^-}K_i(0) &= \mathcal{O}(r^{-N}).
\end{align*}
For all other indices,
\[
\frac{\partial}{\partial b_j^\pm}K_i(0) = 0.
\]
Thus, for sufficiently large $N$, the matrix $D_b K(0,0)$ is invertible. Using the IFT, there exists a unique smooth function $b: Z \rightarrow Y$ with $b(0) = 0$ such that $K(b(d),d) = 0$ for $d$ sufficiently small, which is the case for $N$ sufficiently large, since $d = \mathcal{O}(r^{-N})$. The bound for $b$ comes from projecting $K_i(b(d), d) = 0$ onto $T(\theta)Y^+$ and $T(\theta)Y^-$ together with the estimate $d = \mathcal{O}(r^{-N})$.
\end{proof}
\end{lemma}

Finally, we will use \cref{jumpZ} to derive the jump conditions in the direction of $T(\theta_i) Z_1$.

\begin{lemma}\label{jumpZlemma}
The jump conditions in the direction of $T(\theta_i) Z_1$ are given by
\begin{equation*}
\begin{aligned}
\xi_1 &= \langle T(\theta_1) Z_1(N_1^+), T(\theta_{2}) Q(-N_1^-) \rangle + R_1 = 0 \\
\xi_i &= \langle T(\theta_i) Z_1(N_i^+), T(\theta_{i+1}) Q(-N_i^-) \rangle \\
&\qquad-\langle T(\theta_i) Z_1(-N_{i-1}^-), T(\theta_{i-1}) Q(N_{i-1}^+) \rangle + R_i = 0 && \qquad i = 2, \dots, m-1 \\
\xi_m &= -\langle T(\theta_m) Z_1(-N_{m-1}^-), T(\theta_{m-1}) Q(N_{m-1}^+) \rangle + R_m = 0,
\end{aligned}
\end{equation*}
where the remainder term has bound
\begin{equation}\label{RZbound}
|R_i| \leq C r^{-3N}.
\end{equation}
\begin{proof}
Evaluating the fixed point equations \cref{FPeqs1} at 0 and substituting \cref{aipmest} from \cref{inv1}, we get
\begin{align*}
\tilde{Q}_i^+(0) &- \tilde{Q}_i^-(0) = \Phi_u^+(0, N_i^+; \theta_i) P_0^u(\theta_i) d_i + \Phi_s^-(0, -N_{i-1}^-; \theta_i) P_0^s(\theta_{i-1}) d_{i-1} \\
&+ b_i^+ - b_i^- 
+ \Phi_u^+(0, N_i^+; \theta_i) \tilde{a}_i^+ - \Phi_s^-(0, -N_{i-1}^-; \theta_i) \tilde{a}_{i-1}^- \\
&- \sum_{j = 0}^{N_i^+-1} \Phi_u^+(0, j+1; \theta_i) G_i^+(\tilde{Q}_i^+(j)) 
- \sum_{j = -N_{i-1}^-}^{-1} \Phi_s^-(0, j+1; \theta_i) G_i^-(\tilde{Q}_i^-(j)).
\end{align*}
Next, we project on $\R T(\theta_i) Z_1(0)$ by taking the inner product with $T(\theta_i) Z_1(0)$. Since $b_i^\pm \in T(\theta_i) Y_i^\pm$, these terms are eliminated by the projection. For the leading order terms in \cref{jumpZ}, using equation \cref{didef} and the proof of \cref{inv2}, we have
\begin{align*}
\langle T(\theta_i) Z_1(0), \Phi_u^+(0, N_i^+; \theta_i) P_0^u d_i \rangle
&= \langle T(\theta_i) Z_1(N_i^+), T(\theta_{i+1}) Q(-N_i^-) \rangle + \mathcal{O}(r^{-3N}) \\
\langle T(\theta_i) Z_1(0), \Phi_s^-(0, -N_{i-1}^-; \theta_i) P_0^s d_{i-1} \rangle
&= -\langle T(\theta_i) Z_1(-N_{i-1}^-), T(\theta_{i-1}) Q(N_{i-1}^+) \rangle + \mathcal{O}(r^{-3N}).
\end{align*}

For the higher order terms in \cref{jumpZ}, we substitute $\tilde{Q}_i^\pm$ from \cref{inv1}, $\tilde{a}_i^\pm$ from \cref{inv2}, and $b_i^\pm$ from \cref{inv3nt}. This gives us the remainder bound \cref{RZbound}. Since $N_0^- = N_m^+ = \infty$, one of the two inner product terms vanishes in the jumps $\xi_1$ and $\xi_m$. 
\end{proof}
\end{lemma}

\subsection{Proof of Theorem \ref{ntmulti}}

The existence statement follows from the jump conditions in \cref{jumpZlemma}. The uniform bound $\|\tilde{Q}_i^\pm\| \leq C r^{-N}$ in \cref{Westimates} follows from \cref{inv1} together with the estimates on $a_i^\pm$ and $b_i^\pm$. For the second estimate in \cref{Westimates}, recall that in \cref{inv2} we solved
\begin{equation}\label{Wsubdi}
\tilde{Q}_i^+(N_i^+) - \tilde{Q}_{i+1}^-(-N_i^-) = T(\theta_{i+1}) Q(-N_i^-) - T(\theta_i) Q(N_i^+).
\end{equation}
Apply the projection $P^u_-(-N_i^-; \theta_{i+1})$, noting that it acts as the identity on $T(\theta_{i+1}) Q(-N_i^-)$. We look at the three remaining terms in \cref{Wsubdi} one at a time. For $T(\theta_i) Q(N_i^+)$, we follow the proof of \cref{inv2} and use the estimate \cref{projexpdecay} to get
\begin{align*}
P^u_-(-N_i^-; \theta_{i+1})T(\theta_i) Q(N_i^+)
&= \mathcal{O}(r^{-2N}).
\end{align*}
For $\tilde{Q}_i^+(N_i^+)$, we use the fixed point equations \cref{FPeqs1} and the uniform bound on $\tilde{Q}_i^\pm$ from \cref{inv1} to get
\begin{align*}
(I - &P^u_-(-N_i^-; \theta_{i+1})) \tilde{Q}_i^+(N_i^+) = P^s_-(-N_i^-; \theta_{i+1}) \tilde{Q}_i^+(N_i^+) = \mathcal{O}(r^{-2N}),
\end{align*}
from which it follows that
\[
P^u_-(-N_i^-; \theta_{i+1}) \tilde{Q}_i^+(N_i^+) = \tilde{Q}_i^+(N_i^+) + \mathcal{O}(r^{-2N}).
\]
For $\tilde{Q}_{i+1}^-(-N_i^-)$, we follow a similar procedure to conclude that
\[
P^u_-(-N_i^-; \theta_{i+1}) \tilde{Q}_{i+1}^-(-N_i^-) = \mathcal{O}(r^{-2N}).
\]
Combining all of these gives us the second estimate in \cref{Westimates}. For the third estimate in \cref{Westimates}, we apply the projection $P^s_+(N_i^+; \theta_i)$ to \cref{Wsubdi} and follow the same procedure.

\subsection{Proof of Theorem \ref{transversemulti}}

In the transverse intersection case, we can decompose $\R^d$ as $\R^d = Y^+ \oplus Y^-$, where $Y^+ = T_{Q(0)} W^s(0)$ and $Y^- = T_{Q(0)} W^u(0)$. \cref{inv1} and \cref{inv2} are identical. To obtain a multi-pulse, all that remains to do is solve 
\[
\tilde{Q}_i^+(0) - \tilde{Q}_i^-(0) = P_{T(\theta_i)Y^+ \oplus T(\theta_i)Y^-}( \tilde{Q}_i^+(0) - \tilde{Q}_i^-(0) ) = 0,
\]
which is done in \cref{inv3nt}. There are no remaining jump conditions to satisfy.

\section{Proof of Theorem \ref{stabilitytheorem}}

In this section, we will prove \cref{stabilitytheorem}, which provides a means of locating the interaction eigenvalues associated with a multi-pulse. Throughout this section, we will assume \cref{symmetryhyp}, \cref{boundstatehyp}, \cref{intersectionhyp}, and \cref{melnikovhyp}. The technique of the proof is similar to the proof of \cite[Theorem 2]{Sandstede1998}.

\subsection{Setup}

Using \cref{ntmulti}, let $Q_m(n)$ be an $m-$pulse solution to \cref{diffeq}, constructed using \cref{ntmulti} using pulse distances $N_i$ and symmetry parameters $\theta_i$. Write $Q_m(n)$ piecewise as
\begin{equation}\label{Qmpiecewise}
\begin{aligned}
Q_i^-(n) &= T(\theta_i) Q(n) + \tilde{Q}_i^-(n) && n \in [-N_{i-1}^-, 0] \\
Q_i^+(n) &= T(\theta_i) Q(n) + \tilde{Q}_i^+(n) && n \in [0, N_i^+].
\end{aligned}
\end{equation}
From \cref{ntmulti} and \cref{Qdecay}, we have the following bounds:
\begin{equation}\label{stabbounds1}
\begin{aligned}
Q_1(n) &= \mathcal{O}(r^{-|n|}) \\
\|\tilde{Q}\| &\leq C r^{-N} \\
|\tilde{Q}_{i+1}^-(-N_i^-) - T(\theta_i) Q_1(N_i^+)| &\leq C r^{-2N} \\
|\tilde{Q}_i^+(N_i^+) - T(\theta_{i+1}) Q_1(-N_i^-)| &\leq C r^{-2N}.
\end{aligned}
\end{equation}

Recall that the eigenvalue problem is given by 
\begin{align}\label{latticeEVP2}
V(n+1) = DF(Q_m(n)) V(n) + \lambda B V(n).
\end{align}
Following \cref{DFkernel1} and \cref{DFkernel2}, we have 
\begin{equation}\label{DFQmkernel}
\begin{aligned}
T'(0)Q_m(n+1) &= DF(Q_m(n))T'(0)Q_m(n) \\
(\partial_\omega Q_m)(n+1) &= DF(Q_m(n))\partial_\omega Q_m(n) + B T'(0)Q_m(n).
\end{aligned}
\end{equation}
As in \cite{Sandstede1998}, we will take an ansatz for the eigenfunction $V(n)$ which is a piecewise perturbation of the kernel eigenfunction. If we follow \cite{Sandstede1998} and 
use an ansatz of the form 
\begin{equation*}
V_i^\pm(n) = 
d_i T'(0)Q_m(n) + W_i^\pm(n),
\end{equation*}
we will obtain a Melnikov sum of the form \cref{MelnikovM1zero}, which is 0. Instead, we will take a piecewise ansatz of the form
\begin{equation}\label{Viansatz2}
V_i^\pm(n) = 
d_i [ T'(0)Q_m(n) + \lambda \partial_\omega Q_m(n) ] + W_i^\pm(n),
\end{equation}
where $d_i \in C$. Substituting this into \cref{latticeEVP2}, and simplifying by using \cref{DFQmkernel}, the eigenvalue problem becomes
\begin{align}\label{Weq1}
W_i^\pm(n+1)
&= DF(T(\theta_i) Q(n) ) W_i^\pm(n) + G_i^\pm(n)W_i^\pm(n) + \lambda B W_i^\pm(n) + d_i \lambda^2 B T_m(n),
\end{align}
where
\begin{equation}
G_i^\pm(n) = DF(Q_m(n)) - DF(T(\theta_i) Q(n) ).
\end{equation}

In addition to solving \cref{Weq1}, the eigenfunction must satisfy matching conditions at $n = \pm N_i$ and $n = 0$. Thus the system of equations we need to solve is
\begin{equation}\label{eigWsystem1}
\begin{aligned}
W_i^\pm(n) = DF(T(\theta_i) Q(n) ) W_i^\pm(n) &+ (G_i^\pm(n) + \lambda B) W_i^\pm(n) + \lambda^2 d_i B \tilde{H}_i^\pm(n) \\
W_i^+(N_i^+) - W_{i+1}^-(-N_i^-) &= D_i d \\
W_i^\pm(0) &\in \C T(\theta_i)Y^+ \oplus T(\theta_i)Y^- \oplus T(\theta_i) Z_1(0) \\ 
W_i^+(0) - W_i^-(0) &= 0 ,
\end{aligned}
\end{equation}
where
\begin{equation}\label{defDid}
\begin{aligned}
D_i d &= [ T(\theta_{i+1}) T'(0)Q(-N_i^-) + T'(0)\tilde{Q}_{i+1}^-(-N_i^-)] d_{i+1}
- [ T(\theta_i) T'(0)Q(N_i^+) + T'(0)\tilde{Q}_i^+(N_i^+)] d_i \\
&+ \lambda[ \partial_\omega Q_{i-1}^-(-N_i^-) d_{i+1}
- \partial_\omega Q_i^+(N_i^+)] d_i \\
\end{aligned}
\end{equation}
and
\begin{equation}\label{defHtildeH}
\begin{aligned}
\tilde{H}_i^\pm(n) &= \partial_\omega Q_i^\pm(n)  \\
H_i(n) &= T(\theta_i) \partial_\omega Q(n).
\end{aligned}
\end{equation}
We can require the third condition in \cref{eigWsystem1} since perturbations in the direction of $T(\theta_i)T'(0)Q(0)$ are handled by the $d_i T'(0)Q_m(0) = d_i T'(0)Q(n) + d_i T'(0)T'(0)\tilde{Q}_i^\pm(n)$ term in \cref{Viansatz2}. 

As in \cite{Sandstede1998} and the previous section, we will generally not be able to solve \cref{eigWsystem1}. Instead, we will relax the fourth condition in \cref{eigWsystem1} to get the system
\begin{align}
W_i^\pm(n) = DF(T(\theta_i) Q(n) ) W_i^\pm(n) &+ (G_i^\pm(n) + \lambda B) W_i^\pm(n) + \lambda^2 d_i B \tilde{H}_i^\pm(n) \label{eigsystem1} \\
W_i^+(N_i^+) - W_{i+1}^-(-N_i^-) &= D_i d \label{eigsystem2} \\
W_i^\pm(0) &\in T(\theta_i) Y^+ \oplus T(\theta_i) Y^- \oplus \C T(\theta_i) Z_1(0) \label{eigsystem3a} \\
W_i^+(0) - W_i^-(0) &\in \C T(\theta_i) Z_1(0). \label{eigsystem3b} 
\end{align}
Using Lin's method, we will be able to find a unique solution to this system. This solution, however, will generically have $m-1$ jumps at $n = 0$. Thus a solution to this system is eigenfunction if and only if the $m-1$ jump conditions
\begin{equation*}
\xi_i = \langle T(\theta_i) Z_1(0), W_i^+(0) - W_i^-(0) \rangle = 0
\end{equation*}
are satisfied. Using the bounds \cref{stabbounds1}, we have the estimates
\begin{equation}\label{stabbounds2}
\begin{aligned}
\|G_i^\pm\| &\leq C r^{-N} \\
\|\tilde{H}_i^\pm - H\| &\leq C r^{-N}. \\
\end{aligned}
\end{equation}

\subsection{Fixed point formulation}

As in \cite{Sandstede1998}, we write equation \cref{eigsystem1} as a fixed point problem using the discrete variation of constants formula from \cref{VOC} together with projections on the stable and unstable subspaces of the exponential dichotomy from \cref{dichotomy}. Let $\delta > 0$ be small, and choose $N$ sufficiently large so that $r^{-N} < \delta$. Let $\Phi(m, n; \theta_i)$ be the family of evolution operators for the equations \cref{Vtheta}. Define the spaces
\begin{align*}
V_W &= \ell^\infty([-N_{i-1}, 0]) \oplus \ell^\infty([0, N_i])  \\
V_a &= \bigoplus_{i=0}^{n-1} E^u \oplus E^s \\
V_b &= \bigoplus_{i=0}^{n-1} \ran P_-^u(0; \theta_i) \oplus \ran P_+^s(0; \theta_i)\\
V_\lambda &= B_\delta(0) \subset \C \\
V_d &= \C^d.
\end{align*}
Then for
\begin{align*}
W = (W_i^-, W_i^+) &\in V_W \\
a = (a_i^-, a_i^+) &\in V_a \\
b = (b_i^-, b_i^+) &\in V_b \\
\lambda &\in V_\lambda,
\end{align*}
the fixed point equations for the eigenvalue problem are
\begin{equation}\label{fpeig}
\begin{aligned}
W_i^-(n) &= 
\Phi_s^-(n, -N_{i-1}^-; \theta_i) a_{i-1}^- + \sum_{j = -N_{i-1}^-}^{n-1} \Phi_s^-(n, j+1; \theta_i)
[(G_i^-(j) + \lambda B) W_i^-(j) + \lambda^2 d_i B \tilde{H}_i^-(j)]
 \\
&+ \Phi_u^-(n, 0; \theta_i) b_i^- - \sum_{j = n}^{-1} \Phi_u^-(n, j+1; \theta_i) 
[(G_i^-(j) + \lambda B) W_i^-(j) + \lambda^2 d_i B \tilde{H}_i^-(j)] \\
W_i^+(n) &= \Phi_s^+(n, 0; \theta_i) b_i^+ + \sum_{j = 0}^{n-1} \Phi_s^+(n, j+1; \theta_i) 
[(G_i^+(j) + \lambda B) W_i^+(j) + \lambda^2 d_i B \tilde{H}_i^+(j)] \\
&+ \Phi_u^+(n, N_i^+; \theta_i) a_i^+ - \sum_{j = n}^{N_i^+-1} \Phi_u^+(n, j+1; \theta_i) 
[(G_i^+(j) + \lambda B) W_i^+(j) + \lambda^2 d_i B \tilde{H}_i^+(j)],
\end{aligned}
\end{equation}
where $a_0^- = a_m^+ = 0$ and the sums are defined to be $0$ if the upper index is smaller than the lower index. Since we are taking $a_0^- = a_m^+ = 0$, the corresponding equations are
\begin{align*}
W_1^-(n) &= \sum_{j = -\infty}^{n-1} \Phi_s^-(n, j+1; \theta_1)
[(G_i^-(j) + \lambda B) W_i^-(j) + \lambda^2 d_i B \tilde{H}_i^-(j)]
 \\
&+ \Phi_u^-(n, 0; \theta_1) b_i^- - \sum_{j = n}^{-1} \Phi_u^-(n, j+1; \theta_1) 
[(G_i^-(j) + \lambda B) W_i^-(j) + \lambda^2 d_i B \tilde{H}_i^-(j)] \\
W_m^+(n) &= \Phi_s^+(n, 0; \theta_m) b_i^+ + \sum_{j = 0}^{n-1} \Phi_s^+(n, j+1; \theta_m) 
[(G_i^+(j) + \lambda B) W_i^+(j) + \lambda^2 d_i B \tilde{H}_i^+(j)] \\
&- \sum_{j = n}^{\infty} \Phi_u^+(n, j+1; \theta_m) 
[(G_i^+(j) + \lambda B) W_i^+(j) + \lambda^2 d_i B \tilde{H}_i^+(j)].
\end{align*}

\subsection{Inversion}

We will now solve the eigenvalue problem series of lemmas. This is very similar to the procedure in \cite{Sandstede1998}. First, we use the fixed point equations \cref{fpeig} to solve for $W_i^\pm$. 

\begin{lemma}\label{eiginv1}
There exists an operator $W_1: V_\lambda \times V_a \times V_b \times V_d \rightarrow V_W$ such that
\[
W = W_1(\lambda)(a,b,d)
\]
is a solution to \cref{eigsystem1} for $(a,b,d)$ and $\lambda$. The operator $W_1$ is analytic in $\lambda$, linear in $(a,b,d)$, and has bound
\begin{equation}\label{W1bound}
\|W_1(\lambda)(a,b,d)\| \leq C \left( |a| + |b| + |\lambda|^2 |d| \right).
\end{equation}

\begin{proof}
Rewrite the fixed point equations \cref{fpeig} as
\[
(I - L_1(\lambda))W = L_2(\lambda)(a,b,d),
\]
where $L_1(\lambda): V_W \rightarrow V_W$ is the linear operator composed of terms in the fixed point equations involving $W$
\begin{align*}
(L_1(\lambda)W)_i^-(n) &= \sum_{j = -N_{i-1}^-}^{n-1} \Phi_s^-(n, j+1; \theta_i)
(G_i^-(j) + \lambda B) W_i^-(j) \\
&\qquad- \sum_{j = n}^{-1} \Phi_u^-(n, j+1; \theta_i) 
(G_i^-(j) + \lambda B) W_i^-(j)\\
(L_1(\lambda)W)_i^+(n) &= \sum_{j = 0}^{n-1} \Phi_s^+(n, j+1; \theta_i) 
(G_i^+(j) + \lambda B) W_i^+(j) \\
&\qquad-\sum_{j = n}^{N_i^+-1} \Phi_u^+(n, j+1; \theta_i) 
(G_i^+(j) + \lambda B) W_i^+(j)
\end{align*}
and $L_2(\lambda): V_\lambda \times V_a \times V_b $ is the linear operator composed of terms in the fixed point equations not involving $W$.
\begin{align*}
(L_2(\lambda)(a,b,d))_i^-(n) &= 
\Phi_s^-(n, -N_{i-1}^-; \theta_i) a_{i-1}^- + \sum_{j = -N_{i-1}^-}^{n-1} \Phi_s^-(n, j+1; \theta_i)
\lambda d_i B \tilde{H}_i^-(j)
 \\
&\qquad+ \Phi_u^-(n, 0; \theta_i) b_i^- - \sum_{j = n}^{-1} \Phi_u^-(n, j+1; \theta_i) 
\lambda d_i B \tilde{H}_i^-(j) \\
(L_2(\lambda)(a,b,d))_i^+(n) &= \Phi_s^+(n, 0; \theta_i) b_i^+ + \sum_{j = 0}^{n-1} \Phi_s^+(n, j+1; \theta_i)\lambda^2 d_i B \tilde{H}_i^+(j) \\
&\qquad+ \Phi_u^+(n, N_i^+; \theta_i) a_i^+ - \sum_{j = n}^{N_i^+-1} \Phi_u^+(n, j+1; \theta_i)\lambda^2 d_i B \tilde{H}_i^+(j).
\end{align*}
Using the exponential dichotomy bounds from \cref{dichotomy}, we obtain the following uniform bounds for $L_1$ and $L_2$.
\begin{align*}
\|L_1(\lambda)W)\| &\leq C \left(\|G\| + |\lambda| \right)\|W\| \leq C \delta \|W\| \\
\|L_2(\lambda)(a,b,d))\| &\leq C\left( |a| + |b| + |\lambda|^2 |d| \right).
\end{align*}
For sufficiently small $\delta$, $\|(L_1(\lambda)W)\| < 1$, thus $I - L_1(\lambda)$ is invertible on $V_W$. The inverse $(I - L_1(\lambda))^{-1}$ is analytic in $\lambda$, and we obtain the solution 
\[
W = W_1(\lambda)(a,b,d) = (I - L_1(\lambda))^{-1} L_2(\lambda(a,b,d),
\]
which is analytic in $\lambda$, linear in $(a, b, d)$, and for which we have the estimate
\begin{equation*}
\|W_1(\lambda)(a,b,d)\| \leq C \left( |a| + |b| + |\lambda|^2 |d| \right).
\end{equation*}
\end{proof}
\end{lemma}

In the next lemma, we solve equation \cref{eigsystem2}, which is the matching condition at the tails of the pulses.

\begin{lemma}\label{eiginv2}
There exist operators 
\begin{align*}
A_1 : V_\lambda \times V_b \times V_d \rightarrow V_a \\
W_2 : V_\lambda \times V_b \times V_d \rightarrow V_W
\end{align*}
such that $(a, w) = (A_1(\lambda)(b,d), W_2(\lambda)(b,d)$ solves \cref{eigsystem1} and \cref{eigsystem2} for any $(b, d)$ and $\lambda$. These operators are analytic in $\lambda$, linear in $(b,d)$, and have bounds 
\begin{align}
|A_1(\lambda)(b, d)| &\leq C \left( (r^{-N} + \|G\| + |\lambda| ) |b| + (|\lambda|^2 + |D| ) |d| \right) \label{A1bound} \\
\|W_2(\lambda)(b,d)\| &\leq C \left( |b| + (|\lambda|^2 + |D|) |d| \right). \label{W2bound}
\end{align}
Furthermore, we can write
\begin{align*}
a_i^+ &= P_0^u(\theta_i) D_i d + A_2(\lambda)_i(b,d) \\
a_i^- &= -P_0^s(\theta_i) D_i d + A_2(\lambda)_i(b,d),
\end{align*}
where $A_2$ is a bounded linear operator with bound
\begin{align}\label{A2bound}
|A_2(\lambda)(b,d)| \leq 
C\left( (r^{-N} + \|G\| + |\lambda| )|b| + (r^{-N} + \|G\| + |\lambda|)|D||d| + |\lambda|^2 |d|  \right).
\end{align}

\begin{proof}
Substituting the fixed point equations \cref{fpeig} into equation \cref{eigsystem2} and recalling that $\Phi_s^-(-N_i^-, -N_i^-; \theta_{i+1}) = P_-^s(-N_i^-; \theta_{i+1}$, $\Phi_u^+(N_i^+, N_i^+; \theta_i) = P_+^u(N_i^+; \theta_{i})$, $a_i^- \in E^s(\theta_i)$, and $a_i^+ \in E^u(\theta_i)$, we have
\begin{align}
D_i d = a_i^+ &- a_i^- + (P_u^+(N_i^+; \theta_i) - P_0^u) a_i^+ - (P_s^-(-N_i^-; \theta_{i+1}) - P_0^s) a_i^- \\
&+ \Phi_s^+(N_i^+, 0; \theta_i) b_i^+ - \Phi_u^-(-N_i^-, 0; \theta_{i+1}) b_i^- \nonumber \\
&+ \sum_{j = 0}^{N_i^+-1} \Phi_s^+(N_i^+, j+1; \theta_i) 
[(G_i^+(j) + \lambda B) W_i^+(j) + \lambda^2 d_i B \tilde{H}_i^+(j)] \nonumber \\
&- \sum_{j = -N_i^-}^{-1} \Phi_u^-(-N_i^-, j+1; \theta_{i+1}) 
[(G_i^-(j) + \lambda B) W_i^-(j) + \lambda^2 d_i B \tilde{H}_i^-(j)] . \nonumber \\
\end{align}
Substituting $W = W_1(\lambda)(a, b, d)$ from \cref{eiginv1}, we obtain an equation of the form 
\begin{equation}\label{Dideq2}
D_i d = (a_i^+ - a_i^-) + L_3(\lambda)_i(a,b,d).
\end{equation}
Using \cref{dichotomy}, the bound for $W_1$ from \cref{eiginv1}, and the estimates \cref{projexpdecay} from \cref{dichotomy}, the linear operator $L_3$ has uniform bound
\begin{align}\label{L3bound}
L_3(\lambda)(a,b,d)| &\leq C\left( (r^{-N} + \|G\| + |\lambda| ) (|a| + |b|) + |\lambda|^2 |d|  \right) \\
&\leq C \delta |a| + C\left( (r^{-N} + \|G\| + |\lambda| ) |b| + |\lambda|^2 |d|  \right) . \nonumber
\end{align}

Define the map
\[
J_1: V_a \rightarrow \bigoplus_{j=1}^{m-1} \C^d
\]
by $(J_1)_i(a_i^+, a_i^-) = a_i^+ - a_i^-$. Since $E^s \oplus E^u = \C^d$, the map $J_1$ is a linear isomorphism. Let
\[
K_1(a) = J_1 (a) + L_3(\lambda)(a, 0, 0) = J_1( I + J_1^{-1} L_3(\lambda)(a, 0) ).
\]
For sufficiently small $\delta$, $\|J_1^{-1} L_3(\lambda)(a, 0, 0)\| < 1$, thus the operator $K_1(a)$ is invertible. We can then solve for $a$ to get
\[
a = A_1(\lambda)(b, d) = S_i^{-1}(-D d - L_3(\lambda)(b, d)),
\]
which has uniform bound
\begin{equation*}
|A_1(\lambda)(b, d)| \leq C \left( (r^{-N} + \|G\| + |\lambda| ) |b| + (|\lambda|^2 + |D| ) |d|  \right).
\end{equation*}
We plug this estimate into $W_1$ to get $W_2(\lambda)(b,d)$, which satisfies the bound
\begin{equation*}
\|W_2(\lambda)(b,d)\| \leq C \left( |b| + (|\lambda|^2 + |D|) |d| \right).
\end{equation*}

\noindent Finally, we project \cref{Dideq2} onto $E^s(\theta_i)$ and $E^u(\theta_i)$ to get
\begin{align*}
a_i^+ &= P_0^u(\theta_i) D_i d - P_0^u(\theta_i) L_3(\lambda)_i(a,b,d) \\
a_i^- &= -P_0^s(\theta_i) D_i d + P_0^s(\theta_i) L_3(\lambda)_i(a,b,d).
\end{align*}
Substituting $A_1(\lambda)(b,d)$ for $a$ we obtain the equations
\begin{align*}
a_i^+ &= P_0^u(\theta_i) D_i d + A_2(\lambda)_i(b,d) \\
a_i^- &= -P_0^s(\theta_i) D_i d + A_2(\lambda)_i(b,d).
\end{align*}
Substituting the bound for $A_1$ into the bound for $L_3$, we obtain the uniform bound
\begin{align*}
|A_2(\lambda)(b,d)| \leq 
C\left( (r^{-N} + \|G\| + |\lambda| )|b| + (r^{-N} + \|G\| + |\lambda|)|D||d| + |\lambda|^2 |d|  \right).
\end{align*}
\end{proof}
\end{lemma}

The last step in the inversion is to satisfy equations \cref{eigsystem3a} and \cref{eigsystem3b}. Since we have the decomposition
\begin{equation}\label{decomp1}
C^d = \C T(\theta_i) Z_1(0) \oplus \C T(\theta_i) T'(0)Q(0) \oplus T(\theta_i) Y^+ \oplus T(\theta_i) Y^-,
\end{equation}
these two equations are equivalent to the three projections
\begin{equation}\label{projeq}
\begin{aligned}
P(T(\theta_i) T'(0)Q(0)) W_i^- &= 0 \\
P(T(\theta_i) T'(0)Q(0)) W_i^+ &= 0 \\
P(T(\theta_i) Y^+ \oplus T(\theta_i) Y^-) (W_i^+ - W_i^-) &= 0,
\end{aligned}
\end{equation}
where the kernel of each projection is the remaining elements of the direct sum decomposition \cref{decomp1}. Since we have eliminated any component in $T(\theta_i) T'(0)Q(0)$ in the first two projections, we do not need it in the third projection.

We decompose $b_i^\pm$ uniquely as $b_i^\pm = x_i^\pm + y_i^\pm$, where $x_i^\pm \in \C T(\theta_i) T'(0)Q(0)$ and $y_i^\pm \in T(\theta_i) Y^\pm$. In the next lemma, we solve the equations \cref{projeq}.

\begin{lemma}\label{eiginv3}
There exist operators 
\begin{align*}
B_1 : V_\lambda \times V_d \rightarrow V_b \\
A_3 : V_\lambda \times V_d \rightarrow V_a \\
W_3 : V_\lambda \times V_d \rightarrow V_W
\end{align*}
such that $(a, b, W) = (A_3(\lambda)(d), B_1(\lambda)(d), W_2(\lambda)(d)$ solves \cref{eigsystem1}, \cref{eigsystem2}, \cref{eigsystem3a}, and \cref{eigsystem3b} for any $d$ and $\lambda$. These operators are analytic in $\lambda$, linear in $d$, and have bounds 
\begin{align}
|B_1(\lambda)(d)| &\leq C \left( (r^{-N} + \|G\| + |\lambda|)|D| |d| + |\lambda|^2 |d| \right) \label{B1bound} \\
|A_3(\lambda)(d)| &\leq C \left(|\lambda|^2 + |D|\right)|d| \label{A3bound} \\
\|W_3(\lambda)(d)\| &\leq C \left(|\lambda|^2 + |D|\right)|d|. \label{W3bound} \\
\end{align}
Furthermore, we can write
\begin{align*}
a_i^+ &= P_0^u D_i d + A_4(\lambda)_i(d) \\
a_i^- &= -P_0^s D_i d + A_4(\lambda)_i(d),
\end{align*}
where $A_4$ is a bounded linear operator with estimate
\begin{align}\label{A4bound}
|A_4(\lambda)(d)| &\leq 
C\left( (r^{-N} + \|G\| + |\lambda|)|D||d| + |\lambda|^2 |d|  \right).
\end{align}

\begin{proof}
At $n = 0$, the fixed point equations \cref{fpeig} become 
\begin{align*}
W_i^-(0) &= x_i^- + y_i^- +
\Phi_s^-(0, -N_{i-1}^-; \theta_i) a_{i-1}^- \\
&+ \sum_{j = -N_{i-1}^-}^{-1} \Phi_s^-(0, j+1; \theta_i)
[(G_i^-(j) + \lambda B) W_i^-(j) + \lambda^2 d_i B \tilde{H}_i^-(j)] \\
W_i^+(0) &= x_i^+ + y_i^+ + \Phi_u^+(0, N_i^+; \theta_i) a_i^+ \\
&- \sum_{j = 0}^{N_i^+-1} \Phi_u^+(0, j+1; \theta_i) 
[(G_i^+(j) + \lambda B) W_i^+(j) + \lambda^2 d_i B \tilde{H}_i^+(j)].
\end{align*}
The equations \cref{projeq} can thus be written as
\begin{equation}\label{projeq2}
\begin{pmatrix}
x_i^- \\ x_i^+ \\ y_i^+ - y_i^-
\end{pmatrix}
= (L_4(\lambda)(b,d))_i \:.
\end{equation}
Using the exponential dichotomy estimates from \cref{dichotomy} and $(a, W) = (A_1(\lambda)(b,d), W_2(\lambda)(b,d))$ from \cref{eiginv2}, we get the uniform bound on $L_4$
\begin{align*}
|L_4(\lambda)(b,d)| 
&\leq C \left( (r^{-2N} + \|G\| + |\lambda|)|b| + 
(r^{-N} + \|G\| + |\lambda|)|D| |d| + |\lambda|^2 |d| )
\right) \\
&\leq C \delta(|x| + |y|) + C \left( (r^{N} + \|G\| + |\lambda|)|D| |d| + |\lambda|^2 |d| \right).
\end{align*}

Define the map
\begin{align*}
J_2: &\left( \bigoplus_{j=1}^n \C T'(0)Q(0) \oplus \C T'(0)Q(0) \right) \oplus
\left( \bigoplus_{j=1}^n Y^- \oplus Y^+ \right)  \\
&\qquad\rightarrow \bigoplus_{j=1}^n \C T'(0)Q(0) \oplus \C T'(0)Q(0) \oplus (Y^- \oplus Y^+)
\end{align*}
by 
\[
J_2( (x_i^+, x_i^-),(y_i^+, y_i^-))_i = ( x_i^+, x_i^-, y_i^+ - y_i^- ).
\]
Since $\C^d = \C T(\theta_i) Z_1(0) \oplus \C T(\theta_i) T'(0)Q(0) \oplus T(\theta_i) Y^- \oplus T(\theta_i) Y^+)$, $J_2$ is an isomorphism. Using this and the fact that $b_i = (x_i^- + y_i^-, x_i^+ + y_i^+)$, we can write \cref{projeq2} as
\begin{equation}\label{projxy2}
J_2( (x_i^+, x_i^-),(y_i^+, y_i^-))_i 
+ L_4(\lambda)_i(b_i, 0) + L_4(\lambda)_i(0, d) = 0.
\end{equation}
Consider the map
\begin{align*}
K_2(b)_i &= J_2( (x_i^+, x_i^-),(y_i^+, y_i^-))_i 
+ L_4(\lambda)_i(b_i, 0). 
\end{align*}
Substituting this in \cref{projxy2}, we have
\begin{align*}
K_2(b) &= -L_4(\lambda)(0, d).
\end{align*}
For sufficiently small $\delta$, the operator $K_2(b)$ is invertible. Thus we can solve for $b$ to get
\begin{equation}
b = B_1(\lambda)(d) = -K_2^{-1} L_4(\lambda)(0, d),
\end{equation}
where we have the uniform bound on $B_1$
\begin{equation}
|B_1(\lambda)(d)| \leq C \left( (r^{-N} + \|G\| + |\lambda|)|D| |d| + |\lambda|^2 |d| \right) .
\end{equation}

We can plug this into $A_1$, $W_2$, and $A_2$ to get operators $A_3$, $W_3$, and $A_4$ with bounds
\begin{align*}
|A_3(\lambda)(d)| &\leq C \left(|\lambda|^2 + |D|\right)|d|\\
\|W_3(\lambda)(d)\| &\leq C \left(|\lambda|^2 + |D|\right)|d| \\
|A_4(\lambda)(d)| &\leq 
C\left( (r^{-N} + \|G\| + |\lambda|)|D||d| + |\lambda|^2 |d|  \right).
\end{align*}
\end{proof}
\end{lemma}

\subsection{Jump conditions}

Given $\lambda$ and $d$, we have used Lin's method to find a unique solution to equations \cref{eigsystem1}, \cref{eigsystem2}, \cref{eigsystem3a}, and \cref{eigsystem3b}, which is given by $W = W_3(\lambda)(d)$. Such a solution will  generically have $m-1$ jumps in the direction of $T(\theta_i) Z_1(0)$, which are given by
\begin{equation}\label{jumpIP}
\xi_i = \langle T(\theta_i) Z_1(0), W_i^+(0) - W_i^-(0) \rangle.
\end{equation}
In the next lemma, we derive formulas for these jumps.

\begin{lemma}\label{jumpcond}
$W_i^+(0) = W_i^-(0)$ for $i = 1, \dots, m-1$ if and only if the $m-1$ jump conditions
\begin{equation}\label{xicond}
\xi_i = \langle T(\theta_i) Z_1(0), W_i^+(0) - W_i^-(0) \rangle = 0
\end{equation}
are satisfied. The jumps $\xi_i$ can be written as 
\begin{equation}\label{xieq}
\begin{aligned}
\xi_i = \langle T(\theta_i) &Z_1(N_i^+), P_0^u(\theta_i) D_i d \rangle 
+ \langle T(\theta_i) Z_1(-N_{i-1}^-), P_0^s(\theta_{i-1}) D_{i-1} d \rangle \\ 
&- \sum_{j = -\infty}^{\infty} \langle Z_1(j+1), B \partial_\omega Q(j)\rangle + R(\lambda)_i(d),
\end{aligned}
\end{equation}
where the remainder term $R(\lambda)(d)$ has bound
\begin{align}\label{xiRbound}
|R(\lambda)(d)| \leq C\left( (r^{-N} + \|G\| + |\lambda|)( (r^{-N} + \|G\| + |\lambda|)|D| + |\lambda|^2 \right).
\end{align}

\begin{proof}
From the previous lemma, the fixed point equations at $n = 0$ are given by 
\begin{equation}\label{fpat0}
\begin{aligned}
W_i^-(0) &= b_i^- +
\Phi_s^-(0, -N_{i-1}^-; \theta_i) a_{i-1}^- \\
&+ \sum_{j = -N_{i-1}^-}^{-1} \Phi_s^-(0, j+1; \theta_i)
[(G_i^-(j) + \lambda B) W_i^-(j) + \lambda^2 d_i B \tilde{H}_i^-(j)] \\
W_i^+(0) &= b_i^+ + \Phi_u^+(0, N_i^+; \theta_i) a_i^+ \\
&- \sum_{j = 0}^{N_i^+-1} \Phi_u^+(0, j+1; \theta_i) 
[(G_i^+(j) + \lambda B) W_i^+(j) + \lambda^2 d_i B \tilde{H}_i^+(j)].
\end{aligned}
\end{equation}
To evaluate \cref{jumpIP}, we will compute the inner product of each of the terms in \cref{fpat0} with $T(\theta_i)Z_1(0)$. The $b_i^\pm$ terms will vanish since they lie in spaces orthogonal to $T(\theta_i) Z_1(0)$. We will evaluate the remaining terms in turn. For the terms involving $a$, we substitute $A_4$ from \cref{eiginv3} to get
\begin{align*}
\langle T(\theta_i) &Z_1(0), \Phi_s^-(0, -N_{i-1}^-; \theta_i) a_{i-1}^- \rangle \\
&= -\langle T(\theta_i) Z_1(-N_{i-1}^-), P_0^s(\theta_{i-1}) D_{i-1} d \rangle + \mathcal{O}\left(r^{-N}( (r^{-N} + \|G\| + |\lambda|)|D| + |\lambda|^2 )|d| \right) \\
\langle T(\theta_i) &Z_1(0), \Phi_u^+(0, N_i^+; c_i) a_i^+ \rangle \\
&= \langle T(\theta_i) Z_1(N_i^+), P_0^u(\theta_i) D_i d \rangle + \mathcal{O}\left(r^{-N}( (r^{-N} + \|G\| + |\lambda|)|D| + |\lambda|^2 )|d| \right).
\end{align*}

The sums involving $\tilde{H}$ give us the higher order Melnikov sum $M_2$.
\begin{align*}
&\langle T(\theta_i) Z_1(0), \sum_{j = -N_{i-1}^-}^{-1} \Phi_s^-(0, j+1; \theta_i) B \tilde{H}_i^-(j) + \sum_{j = 0}^{N_i^+-1} \Phi_u^+(0, j+1; \theta_i) B \tilde{H}_i^+(j) \rangle \\
&= \sum_{j = -N_{i-1}^-}^{-1} \langle T(\theta_i) Z_1(j+1), B T(\theta_i) S_1(j) \rangle + \sum_{j = 0}^{N_i^+-1} \langle T(\theta_i) Z_1(j+1), B T(\theta_i) S_1(j) \rangle + \mathcal{O}(r^{-N})\\
&= \sum_{j = -\infty}^{\infty} \langle T(\theta_i) Z_1(j+1), B T(\theta_i) S_1(j)\rangle + \mathcal{O}(r^{-N}) \\
&= \sum_{j = -\infty}^{\infty} \langle Z_1(j+1), B S_1(j)\rangle + \mathcal{O}(r^{-N}),
\end{align*}
where in the last line we used the fact that $T(\theta)$ is unitary and commutes with $B$.

Finally, we need to obtain bound for the sum involving $W$. To do this, as in \cite{Sandstede1998}, we will need an improved bound for $W$. Plugging in the bounds for $A_3$, $W_3$, and $B_1$ into the fixed point equations \cref{fpeig}, we have piecewise bounds
\begin{align*}
|W_i^-(n)| &\leq C( r^{-(N_{i-1}^- + n)}|D| +  
(r^{-N} + \|G\| + |\lambda|)|D| + |\lambda|^2 )|d| \\
|W_i^-(n)| &\leq C( r^{-(N_i^+ - n)}|D| +  
(r^{-N} + \|G\| + |\lambda|)|D| + |\lambda|^2 )|d|.
\end{align*}

Since $Z_1(n) = T'(0)Q(n)$, it follows from \cref{Qdecay} that $Z_1(n) \leq C r^{|n|}$. Since $D F(0)$ is hyperbolic, we can find a constant $\tilde{r} > r$ such that $|Z_1(n)| \leq C \tilde{r}^{-n}$. The price to pay is a larger constant $C$. Using this bounds, the sum involving $W$ becomes
\begin{align*}
&\left| \sum_{j = -N_{i-1}^-}^{-1} \langle Z_1(j+1), 
(G_i^-(j) + \lambda B) W_i^-(j) \rangle \right| \\
&\leq C (\|G\| + |\lambda|) \sum_{j = -N_{i-1}^-}^{-1} \tilde{r}^{-|j+1|} r^{-(N_{i-1}^- + j)}|D||d| + C (\|G\| + |\lambda|)(r^{-N} + \|G\| + |\lambda|)|D| + |\lambda|^2 )|d| \\
&\leq C |D| r^{-N} (\|G\| + |\lambda|)|d| \sum_{j = 1}^\infty \left( \frac{r}{\tilde{r}}\right)^j + C (\|G\| + |\lambda|)(r^{-N} + \|G\| + |\lambda|)|D| + |\lambda|^2 )|d| \\
&\leq C (\|G\| + |\lambda|)(r^{-N} + \|G\| + |\lambda|)|D| + |\lambda|^2 )|d|.
\end{align*}
The infinite sum is convergent by our choice of $\tilde{r}$. We have a similar bound for the other sum. Putting this all together, we obtain the jump equations \cref{xieq} and the remainder bound \cref{xiRbound}.
\end{proof}
\end{lemma}

\subsection{Proof of Theorem \ref{stabilitytheorem}}

Using the estimates \cref{Westimates}, we have
\begin{align*}
T'(0) \tilde{Q}_{i+1}^-(-N_i^-) &= T(\theta_i) T'(0) Q(N_i^+) + \mathcal{O}(r^{-2N}) \\
T'(0) \tilde{Q}_i^+(N_i^+) &= T(\theta_{i+1}) T'(0) Q(-N_i^-) + \mathcal{O}(r^{-2N}),
\end{align*}
since the infinitesimal generator of a group commutes with the group elements. Substituting these into \cref{defDid} and simplifying, we have
\begin{equation}\label{Did4}
\begin{aligned}
D_i d = [ T(&\theta_{i+1}) T'(0)Q(-N_i^-) + T(\theta_i) T'(0)Q(N_i^+) ] d_{i+1} \\
&- [ T(\theta_i) T'(0)Q(N_i^+) + T(\theta_{i+1}) T'(0)Q(-N_i^-) ] d_i 
+\mathcal{O}(r^{-N}( |\lambda| + r^{-N})).
\end{aligned}
\end{equation}
Next, we substitute \cref{Did4} into jump expressions $\xi_i$ from \cref{jumpcond}. For the inner product term $\langle T(\theta_i) Z_1(N_i^+), P_0^u(\theta_i) D_i d \rangle$, we use equation \cref{projexpdecay} to get
\begin{equation*}
\langle T(\theta_i) Z_1(N_i^+), P_0^u(\theta_i) D_i d \rangle 
= \langle T(\theta_i) Z_1(N_i^+), T(\theta_{i+1}) T'(0)Q(-N_i^-) \rangle (d_{i+1} - d_i)
+ \mathcal{O}(r^{-3N}).
\end{equation*}
since $T(\theta)$ is unitary and $\langle Z_1(n), T'(0)Q(n) \rangle = 0$ for all $n$. Similarly, we have
\begin{align*}
\langle T(\theta_i) Z_1(-N_{i-1}^-), P_0^s D_{i-1} d \rangle 
&= \langle T(\theta_i) Z_1(-N_{i-1}^-), T(\theta_{i-1}) T'(0)Q(N_{i-1}^+) \rangle (d_i - d_{i-1}).
\end{align*}
Substituting these into the jump equations, we obtain the jump conditions
\begin{align*}
\xi_i = \langle &T(\theta_i) Z_1(N_i^+), T(\theta_{i+1}) T'(0)Q(-N_i^-) \rangle (d_{i+1} - d_i) \\
&+ \langle T(\theta_i) Z_1(-N_{i-1}^-), T(\theta_{i-1}) T'(0)Q(N_{i-1}^+) \rangle (d_i - d_{i-1}) \\
&- \sum_{j = -\infty}^{\infty} \langle Z_1(j+1), B S_1(j)\rangle + R(\lambda)_i(d).
\end{align*}
For the remainder term, we substitute $|D|, \|G\| = \mathcal{O}(r^{-N})$ into the remainder term in \cref{jumpcond} to get
\begin{align*}
|R(\lambda)(d)| \leq C\left( (r^{-N} + |\lambda|)^3 \right).
\end{align*} 

\section{Proofs of results from section 4}

\subsection{Proof of Theorem \ref{DNLSexisttheorem}}

First, we will look for real-valued solutions to \cref{DNLSequilib}. In this case, the stationary equation \cref{DNLS} reduces to
\begin{equation*}
d(u_{n+1} - 2 u_n + u_{n-1}) - \omega u_n + u_n^3 = 0.
\end{equation*}
For $d \neq 0$, this is equivalent to the first order difference equation $U(n+1) = F(U(n))$, where $U(n) = (u_n, \tilde{u}_n) \in \R^2$, $\tilde{u}_n = u_{n-1}$, and 
\begin{equation}\label{dnlsdiffR2}
F(U) = 
\begin{pmatrix}
\frac{\omega}{d} + 2 & -1 \\
1 & 0
\end{pmatrix}
\begin{pmatrix}
u \\ \tilde{u}
\end{pmatrix}
- \frac{1}{d} 
\begin{pmatrix}
u^3 \\ 0
\end{pmatrix}.
\end{equation}
The symmetry group $G = \{ 1, -1\}$ acts on $\R^2$ via $T(\theta) = \theta I$. For $d, \omega > 0$, $DF(0)$ has a pair of real eigenvalues $\{r, 1/r \}$, where $r$ depends on both $d$ and $\omega$, and is given by \cref{eigr}. As $d \rightarrow \infty$, $r \rightarrow 1$, thus the spectral gap decreases with increasing $d$. As $d \rightarrow 0$, $r \rightarrow \infty$.

It follows that 0 is a hyperbolic equilibrium point with 1-dimensional stable and unstable manifolds. Let $q_n$ be the symmetric, real-valued, on-site soliton solution to DNLS, and let $Q(n) = (q_n, \tilde{q}_n)$ be the primary pulse solution, where $\tilde{q}_n = q_{n-1}$. Since the variational equation does not have a bounded solution, the stable and unstable manifolds intersect transversely. Thus we have satisfied \cref{transversehyp}. Using \cref{transversemulti}, for sufficiently large $N$ (which depends on $r$, thus $\omega$ and $d$) there exist $m-$pulse solutions for any $\theta_i = \pm 1$ and lengths $N_i \geq N$. These correspond to phase differences of $0$ and $\pi$.

We will now show that there are no multi-pulse solutions with phase differences other than $0$ and $\pi$. For this, we write the DNLS equation \cref{DNLSequilib2} as the first order system \cref{diffeq} in $\R^4$. In this formulation, the primary pulse solution is given by $Q(n) = (q_n, 0, \tilde{q}_n, 0)$. The unique bounded solutions to variational equation \cref{vareq} and the adjoint variational equation \cref{adjvareq} are
\begin{align*}
T'(0) Q(n) &= (0, q_n, 0, \tilde{q}_n) \\
Z_1(n) &= (0, -\tilde{q}_n, 0, q_n).
\end{align*}
Using \cref{ntmulti}, for sufficiently large $N$ (which depends on $r$, thus $\omega$ and $d$) there exist $m-$pulse solutions with lengths $N_i^\pm$ and phase parameters $\theta_i$ if any only if the jump conditions \cref{jumpcondexist} are satisfied. Since the symmetry group $T(\theta)$ is unitary, we can rewrite the jump conditions in terms of the phase differences $\Delta \theta_i = \theta_{i+1} - \theta_i$ to get the jump conditions
\begin{align}\label{jumpDNLS}
\xi_i = \langle T(-\Delta \theta_i) Z_1(N_i^+), Q(-N_i^-) \rangle
- \langle T(\Delta \theta_{i-1}) Z_1(-N_{i-1}^-), Q(N_{i-1}^+) \rangle + R_i,
\end{align}
where we take $\Delta \theta_0 = \Delta \theta_m = 0$. The inner product terms in \cref{jumpDNLS} are
\begin{equation}\label{jumpIPs}
\begin{aligned}
\langle T(-\Delta\theta_i) Z_1(N_i^+), Q(-N_i^-) \rangle 
&= -b_i \sin(\Delta\theta_i) \\
\langle T(\Delta\theta_{i-1}) Z_1(-N_{i-1}^-), Q(N_{i-1}^+) \rangle &= -b_{i-1} \sin(\Delta\theta_{i-1}),
\end{aligned}
\end{equation}
where 
\begin{align*}
b_i &= q(N_i^+ - 1)q(N_i^-) - q(N_i^+)q(N_i^- + 1).
\end{align*}
Since the single pulse $q(n)$ is an even function, the $b_i$ are given by \cref{bieq}. Since $q(n)$ is non-negative, even, unimodal, and exponentially decaying ~\cite[Theorem 1]{herrmann_2011}, $q(n)$ is strictly decreasing as $n$ moves away from 0, thus $b_i < 0$ for all $i$. 

Let $s_i = \sin{\Delta\theta_i}$. Substituting equations \cref{jumpIPs} into \cref{jumpDNLS}, the jump conditions become
\begin{equation}\label{jumpDNLS2}
\begin{aligned}
\xi_1 &= -b_1 s_1 + R_1 \\
\xi_i &= b_{i-1} s_{i-1} - b_i s_i + R_i
&& i = 2, \dots, m-1 \\
\xi_m &= b_{m-1} s_{m-1} + R_m.
\end{aligned}
\end{equation}
Since $b_i = \mathcal{O}(r^{-2N})$ and $R_i = \mathcal{O}(r^{-3N})$, the jump conditions can only be satisfied if $s_i = \mathcal{O}(r^{-N})$. Thus we only have to consider that case from here on. Since the steady state equation \cref{DNLSequilib2} has a conserved quantity \cref{DNLSE}, we can eliminate the final equation in \cref{jumpDNLS2} as is done in \cite{SandstedeStrut}. We write the $(m-1)$ remaining jump conditions in matrix form as $H s + R = 0$, where $s = (s_1, \dots, s_{m-1})$ and $H$ is the $(m-1)\times(m-1)$ matrix
\[
H = \begin{pmatrix}
-b_1 \\
b_1 & -b_2 \\
& b_2 & -b_3 \\
&& \ddots & \ddots \\
&&& b_{m-2} & -b_{m-1} \\
\end{pmatrix}.
\]
Since $H$ is lower triangular and all the $b_i$ are nonzero, $B$ is invertible, thus $s = B^{-1}R$ is the unique value of $s$ for which all the jump conditions are satisfied. 

We showed above that for sufficiently large $N$, real-valued multi-pulses exist with phase differences which are either 0 or $\pi$; in all of those cases, $s = 0$. Since $s = B^{-1}R$ is the unique solution which satisfies the jump conditions, and $s = 0$ is also a solution, we conclude that $s = 0$ must be the unique solution that satisfies jump conditions. Thus for sufficiently large $N$, the jump conditions can only be satisfied if all of the phase differences $\Delta \theta_i$ are either 0 or $\pi$. No other phase differences are possible.

\subsection{Proof of Theorem \ref{DNLSeigtheorem}}

To find the interaction eigenvalues for DNLS, we will solve the matrix equation \cref{Elambda} from \cref{stabilitytheorem}. For the higher order Melnikov sum,
\[
M_2 = \sum_{n=-\infty}^\infty \langle Z_1(n+1), B S_1(n) \rangle =
\frac{1}{d} \sum_{n=-\infty}^\infty q(n) q_\omega(n) = \frac{1}{d}M,
\]
where
\[
M = \sum_{n=-\infty}^\infty q(n) q_\omega(n).
\]
We are assuming that $M > 0$. 

For $N$ sufficiently large, we can find the eigenvalues of \cref{multiEVP} using \cref{stabilitytheorem}. The matrix $A$ is given by \cref{DNLSmatrixA}. First, we rescale equation \cref{Elambda} by taking
\begin{align*}
A = r^{-2N} \tilde{A} \\
\lambda = r^{-N} \tilde{\lambda} \\
R(\lambda) = r^{-3N} \tilde{R}(\lambda)
\end{align*}
and dividing by $r^{-2N}$ to get the equivalent equation
\begin{equation}\label{DNLStildeE}
\tilde{E}(\lambda) = 
\det(\tilde{A} - M_2 \tilde{\lambda}^2 I + r^{-N} \tilde{R}(\lambda)) = 0.
\end{equation}
To solve $\tilde{E}(\lambda) = 0$, we need to find the eigenvalues of $\tilde{A}$. Since $\tilde{A}$ is symmetric tridiagonal, its eigenvalues are real. Furthermore, $\tilde{A}$ has an eigenvalue at 0 with corresponding eigenvector $(1, 1, \dots, 1)^T$. Let $\{ \tilde{\mu}_1, \dots, \tilde{\mu}_{m-1}\}$ be the remaining $m-1$ eigenvalues of $A$. Since $b_i < 0$ for all $i$, it follows from \cite[Lemma 5.4]{Sandstede1998} that the signs of $\{ \tilde{\mu}_1, \dots, \tilde{\mu}_{m-1}\}$ are determined by the phase differences $\Delta\theta_i$. Specifically, $A$ has $k_\pi$ negative real eigenvalues (counting multiplicity), where $k_\pi$ is the number of $\Delta\theta_i$ which are $\pi$, and $A$ has $k_0$ positive real eigenvalues (counting multiplicity), where $k_0$ is the number of $\Delta\theta_i$ which are $0$. 

Next, we show that the eigenvalues of $\tilde{A}$ are distinct. The eigenvalue problem $(\tilde{A} - \mu I)v = 0$ is equivalent to the Sturm-Liouville difference equation with Dirichlet boundary conditions
\begin{equation}\label{SLdiff2}
\begin{aligned}
\nabla( p_j \Delta d_j ) &= \mu d_j && \qquad j = 1, \dots, m \\
d_0 &= 0 \\
d_{m+1} &= 0,
\end{aligned}
\end{equation}
where $p_i = \cos(\Delta\theta_i) b_i$, $\Delta$ is the forward difference operator $\Delta f(j) = f(j+1) - f(j)$ and $\nabla$ is the backward difference operator $\nabla f(j) = f(j) - f(j-1)$. It follows from \cite[Corollary 2.2.7]{Jirari1995} that the eigenvalues of \cref{SLdiff2}, thus the eigenvalues of $\tilde{A}$, are distinct.

We can now solve equation \cref{DNLStildeE} for $\lambda$. By \cref{DFQmkernel}, we will always have an eigenvalue at 0 with algebraic multiplicity 2 and geometric multiplicity 1. The remaining eigenvalues result from interaction between the pulses. Let $\eta = r^{-N}$, and rewrite equation \cref{DNLStildeE} as 
\begin{equation}\label{DNLStildeE2}
K(\tilde{\lambda}; \eta) = \det(\tilde{A} - M_2 \tilde{\lambda}^2 I + \eta \tilde{R}(\lambda)).
\end{equation}
For $j = 1, \dots, m-1$, $K(\pm \sqrt{\tilde{\mu}_j / M_2 }; 0) = 0$. Since the eigenvalues of $\tilde{A}$ are distinct, 
\[
\frac{\partial}{\partial \tilde{\lambda}} K(\tilde{\lambda}; 0)\Big|_{\tilde{\lambda} = \pm \sqrt{\tilde{\mu}_j / M_2 }} \neq 0.
\]
Using the implicit function theorem, we can solve for $\tilde{\lambda}$ as a function of $\eta$ near $(\tilde{\lambda}, \eta) = (\pm \sqrt{\tilde{\mu}_j / M_2 }; 0)$. Thus for sufficiently small $\eta$, we can find smooth functions $\tilde{\lambda}_j^\pm(\eta)$ such that $\tilde{\lambda}_j^\pm(0) = \pm \sqrt{\tilde{\mu}_j / M_2 }$ and $K(\tilde{\lambda}_j^\pm(\eta); \eta) = 0$. Expanding $\tilde{\lambda}(\eta)$ in a Taylor series about $\eta = 0$ and taking $\eta = r^{-N}$, we can write $\tilde{\lambda}_j^\pm$ as $\tilde{\lambda}_j^\pm(N) = \tilde{\lambda}_j^\pm + \mathcal{O}(r^{-N})$. Undoing the scaling and taking $M_2 = M/d$, the interaction eigenvalues are given by
\begin{align*}
\lambda^\pm_j &= \pm \sqrt{\frac{d \mu_j}{M}} + \mathcal{O}(r^{-2N}) && j = 1, \dots, m-1 .
\end{align*}

By Hamiltonian symmetry, the eigenvalues of DNLS must come in quartets $\pm \alpha \pm i \beta$. Since the $\mu_j$ are distinct and only come in pairs, the eigenvalues $\lambda_j^\pm$ must be pairs which are real or purely imaginary. Thus there are $(m - 1)$ pairs of nonzero interaction eigenvalues at $\lambda = \pm \lambda_j$, given by 
\begin{align*}
\lambda_j &= \sqrt{\frac{d \mu_j}{M}} + \mathcal{O}(r^{-2N}) && j = 1, \dots, m-1.
\end{align*}
These are either real or purely imaginary, and the remainder term cannot move these off of the real or imaginary axis. Since $M, d > 0$, we conclude that there are $k_\pi$ pairs of purely imaginary eigenvalues and $k_0$ pairs of real eigenvalues.

We note that upon variations of $d$,
these interaction eigenvalues may collide with other eigenvalues including
the ones associated with the continuous spectrum
and lead to quartets as, for example, in some of the
cases in~\cite{Pelinovsky2005}. We can ensure this will not happen by choosing $N$ sufficiently large.

\subsection{Proof of Corollaries \ref{DNLSeigcorr} and \ref{DNLSeigcorr2}}

First, we prove \cref{DNLSeigcorr}. For (i), the matrix $A$ in the case of the 2-pulse has a single eigenvalue $\mu_1 = -\cos(\Delta\theta_1) b_1$. For (ii), the matrix $A$ in the case of the symmetric 3-pulse is given by
\begin{align*}
A &= b \begin{pmatrix}
-\cos(\Delta\theta_1) & \cos(\Delta\theta_1) & 0  \\
\cos(\Delta\theta_1) & -\cos(\Delta\theta_1) - \cos(\Delta\theta_2) & \cos(\Delta\theta_2) \\ 
0 & \cos(\Delta\theta_2) & -\cos(\Delta\theta_2) \\
\end{pmatrix},
\end{align*}
which has nonzero eigenvalues
\[
\mu_{1, 2} = \left( \pm\sqrt{\cos(\Delta\theta_1)^2 - \cos(\Delta\theta_1) \cos(\Delta\theta_2) + \cos(\Delta\theta_2)^2} - \cos(\Delta\theta_1) - \cos(\Delta\theta_2) \right)b.
\]
For the three distinct 3-pulses, these eigenvalues are
\begin{align*}
\mu_{1, 2} = \begin{cases}
-3b, -b & (\Delta\theta_1, \Delta\theta_2) = (0, 0) \\
\pm \sqrt{3}b & (\Delta\theta_1, \Delta\theta_2) = (0, \pi) \\
3b, b & (\Delta\theta_1, \Delta\theta_2) = (\pi, \pi)
\end{cases}.
\end{align*}
For (iii), if $b_i = b$ and $\Delta\theta_i = \Delta\theta$ for all $i$, the eigenvalue problem $(A - \mu I)v = 0$ is equivalent to the difference equation with Neumann boundary conditions
\begin{equation*}
\begin{aligned}
v_{n-1} - 2 v_n + v_{n+1} &- \frac{\mu}{b \cos(\Delta\theta)} v_n = 0 \\
v_0 &= v_1 \\
v_{m+1} &= v_m,
\end{aligned}
\end{equation*}
which has solutions
\begin{align*}
\mu_j &= 2 b \left( \cos\frac{\pi j}{m} - 1 \right) \cos (\Delta\theta) && j = 1, \dots, m.
\end{align*}

For \cref{DNLSeigcorr2}, equation \cref{3pulseeigs} follows from computing the eigenvalues of $A$ explicitly for the 3-pulse and noting that $(\cos \Delta \theta_i)^2 = 1$ since $\Delta \theta_i \in \{0, \pi\}$. We note that for $N_1 < N_2$, $b_1 > b_2$. Thus we write
\[
\sqrt{b_1^2 + b_2^2 - b_1 b_2\cos\Delta\theta_1 \cos\Delta\theta_2} = b_1
\sqrt{1 + \frac{b_2^2}{b_1^2} - \frac{b_2}{b_1} \cos\Delta\theta_1 \cos\Delta\theta_2} 
\]
and expand in a Taylor series to obtain the estimates \cref{3pulsemag}.

\paragraph{Acknowledgments}

This material is based upon work supported by the U.S. National Science Foundation under grants DMS-1148284 (R.P.), DMS-1809074 (P.G.K.), and DMS-1714429 (B.S.). P.G.K. also gratefully acknowledges support from the Leverhulme
Trust during his stay at the University of Oxford.

\bibliography{DiscreteLinPaper.bib}

\end{document}